\documentclass[11pt,reqno,twoside]{amsart}
\usepackage{amsfonts}
\usepackage{amsmath,amssymb,amsopn}
\usepackage[dvips,bottom=1in,right=1in,top=1in, left=1in]{geometry}
\usepackage{fancyvrb}
\usepackage{graphicx}
\usepackage{multirow}
\usepackage{mathrsfs,mathtools}
\usepackage{enumitem}
\usepackage{commath}
\usepackage{hyperref}
\usepackage{orcidlink}
\usepackage{url}
\usepackage{algorithm}
\usepackage[noend]{algpseudocode}
\usepackage{dsfont}
\usepackage{color}
\usepackage{cite}
\usepackage[normalem]{ulem}

\def\RR{{\mathbb{R}}}

\newtheorem{remark}{Remark}[subsection]

\usepackage[capitalize,nameinlink]{cleveref}[0.19]
\crefname{section}{section}{sections}
\crefname{subsection}{subsection}{subsections}
\Crefname{section}{Section}{Sections}
\Crefname{subsection}{Subsection}{Subsections}
\crefname{hypothesis}{Hypothesis}{Hypotheses}
\Crefname{figure}{Figure}{Figures}
\crefformat{equation}{\textup{#2(#1)#3}}
\crefrangeformat{equation}{\textup{#3(#1)#4--#5(#2)#6}}
\crefmultiformat{equation}{\textup{#2(#1)#3}}{ and \textup{#2(#1)#3}}
{, \textup{#2(#1)#3}}{, and \textup{#2(#1)#3}}
\crefrangemultiformat{equation}{\textup{#3(#1)#4--#5(#2)#6}}%
{ and \textup{#3(#1)#4--#5(#2)#6}}{, \textup{#3(#1)#4--#5(#2)#6}}{, and \textup{#3(#1)#4--#5(#2)#6}}

\Crefformat{equation}{#2Equation~\textup{(#1)}#3}
\Crefrangeformat{equation}{Equations~\textup{#3(#1)#4--#5(#2)#6}}
\Crefmultiformat{equation}{Equations~\textup{#2(#1)#3}}{ and \textup{#2(#1)#3}}
{, \textup{#2(#1)#3}}{, and \textup{#2(#1)#3}}
\Crefrangemultiformat{equation}{Equations~\textup{#3(#1)#4--#5(#2)#6}}%
{ and \textup{#3(#1)#4--#5(#2)#6}}{, \textup{#3(#1)#4--#5(#2)#6}}{, and \textup{#3(#1)#4--#5(#2)#6}}

\crefdefaultlabelformat{#2\textup{#1}#3}

\title[OCTANE]{OCTANE - Optimal Control for Tensor-based Autoencoder Network Emergence: Explicit Case}

\thanks{\\ R. Khatri and C. Olson are supported by a U.S. Naval Research Laboratory Base Program (Work Unit: 61A1G8). H. Antil is partially supported by the Office of Naval Research (ONR) under Award NO: N00014-24-1-2147, NSF grant DMS-2408877, and the Air Force Office of Scientific Research (AFOSR) under Award NO: FA9550-22-1-0248. A. Kolshorn is partially supported by the ONR Naval Research Enterprise Internship Program (NREIP) and the Research Training Group in Computation- And Data-Enabled Science (CADES) program (NSF grant DMS-2136228) at Portland State University.}

\author[R.~Khatri]{Ratna Khatri\, \orcidlink{0000-0003-0931-4025}}
\address{R. Khatri, C.~Olson. Optical Sciences Division, Code 5664, U.S. Naval Research Laboratory, Washington, DC 20375, USA}
\email{ratna.khatri.civ@us.navy.mil, colin.c.olson.civ@us.navy.mil}

\author[A.~Kolshorn]{Anthony Kolshorn\, \orcidlink{0009-0000-1546-114X}}
\address{A. Kolshorn, Portland State University,
Portland, OR 97201, USA.}
\email{kolshorn@pdx.edu}

\author[C.~Olson]{Colin Olson}

\author[H.~Antil]{Harbir Antil\, \orcidlink{0000-0002-6641-1449}}
\address{H. Antil, The Center for Mathematics and Artificial Intelligence (CMAI) \& Department of Mathematical Sciences, George Mason University,
Fairfax, VA 22030, USA.}
\email{hantil@gmu.edu}

\date{\today.}

\begin{document}

\begin{abstract}
This paper presents a novel, mathematically rigorous framework for autoencoder-type deep neural networks that combines optimal control theory and low-rank tensor methods to yield memory-efficient training and automated architecture discovery. The learning task is formulated as an optimization problem constrained by differential equations representing the encoder and decoder components of the network and the corresponding optimality conditions are derived via a Lagrangian approach. Efficient memory compression is enabled by approximating differential equation solutions on low-rank tensor manifolds using an adaptive explicit integration scheme. These concepts are combined to form OCTANE (Optimal Control for Tensor-based Autoencoder Network Emergence)—a unified training framework that yields compact autoencoder architectures, reduces memory usage, and enables effective learning, even with limited training data. The framework's utility is illustrated with application to image denoising and deblurring tasks and recommendations regarding governing hyperparameters are provided.  
\end{abstract}


  
\maketitle

\section{Introduction}
Differential equations are powerful tools for modeling dynamic systems and have recently been used to represent deep neural networks (DNNs) \cite{chen2018neural, Ruthotto_2020_DL_OC, Khatri_2020_fdnn, Antil_2022_OCfpde}. In parallel, tensor methods are gaining traction in machine learning due to their ability to efficiently compress high-dimensional data \cite{Wilkinson2019, BacciuDavide2020, Novikov_2015}. This paper introduces a unified framework for deep autoencoder training and architecture design that leverages both optimal control theory and tensor decomposition.

Autoencoders, comprising a composition of encoder $ f: \mathcal{X} \to \mathcal{H} $ and decoder $ g: \mathcal{H} \to \mathcal{X} $ maps, are widely used in applications such as image denoising, anomaly detection \cite{Olson_2020_AE, Nick_2020_AE}, manifold learning, drug discovery \cite{Zhavoronkov2019}, and information retrieval \cite{Salakhutdinov2009}. For a general overview, see also \cite{Goodfellow2016, bank2021autoencoders}. Despite their success, the design of autoencoder architectures—particularly their depth and compression profiles—remains largely heuristic.

Motivated by \cite{Ruthotto_2020_DL_OC, Khatri_2020_fdnn}, we propose modeling autoencoders as an optimal control problem, where the encoder and decoder are governed by coupled nonlinear differential equations. The learning task is to minimize a loss functional subject to these dynamics, which evolve data forward and backward in time through time-dependent operators and biases. This continuous-time formulation offers a principled alternative to traditional heuristically crafted layer-wise design.

To enable tractable compression, we solve the differential equations using a rank-adaptive explicit Euler scheme on low-rank tensor manifolds, following \cite{Rogers_2021}. The state variables \( f \) and \( g \) are discretized using the tensor-train (TT) format \cite{Oseledets_2011}, and all numerical operations are implemented using the \texttt{TT-toolbox} \cite{tt_oseledets_2014}. This scheme dynamically adjusts tensor ranks at each time step, enabling architecture discovery and reducing memory demands.

\medskip
\noindent
\textbf{Contributions.} The main contributions of this work are:
\begin{itemize}
    \item A novel optimal control formulation for autoencoder-type DNNs, with corresponding optimality conditions derived via a Lagrangian framework \cite{Khatri_2020_fdnn}.
    \item An explicit case of the rank-adaptive, tensor-based algorithmic design that discovers the autoencoder network architecture while enforcing compression across layers.
    \item A concrete algorithm, \textbf{OCTANE} (\textbf{O}ptimal \textbf{C}ontrol for \textbf{T}ensor-based \textbf{A}utoencoder \textbf{N}etwork \textbf{E}mergence), validated on image denoising and deblurring tasks constrained by limited training data.
\end{itemize}

Our approach differs from prior works such as \cite{Peng_2020_IP_AE}, which employ autoencoders iteratively for inverse problems (e.g., compressed sensing) without explicitly exploiting the architecture. Likewise, while \cite{Burger_2021_Arch_search} use neural architecture search (NAS) to construct residual autoencoders, our method organically discovers the most compact architecture from a continuous optimization framework grounded in control theory and tensor analysis.

\medskip
\noindent 
{\bf Outline:} 
In \cref{s:Prelim}, we present the relevant notations and definitions, along with a description of the rank-adaptive Euler integration scheme on tensor manifolds. In \cref{s:AE_OC}, we formulate the deep autoencoder model as a continuous-time optimal control problem and derive the corresponding optimality conditions. The discretization of the resulting system is discussed in \cref{s:Disc_AE_OC}. Section~\ref{s:OCTANE} introduces the proposed algorithm, OCTANE, for training the deep autoencoder model. Numerical experiments and results validating the effectiveness of the method are presented in \cref{s:Numerics}. Finally, \cref{s:Discuss} provides concluding remarks and outlines potential directions for future research.

\section{Preliminaries\label{s:Prelim}}

The purpose of this section is to introduce the notations and definitions that we will use throughout the paper. We begin with \cref{t:notations}, where we state the common notations. \\

\begin{table}[h!]
\centering
\caption{Table of Notations.}
\label{t:notations}
\renewcommand{\arraystretch}{1.2}
\begin{tabular}{|c|p{11.2cm}|}
\hline
\textbf{Notation} & \textbf{Description} \\
\hline
$n \in \mathbb{N}$ & Number of distinct samples \\
$n_f \in \mathbb{N}$ & Number of sample features \\
$N_e \in \mathbb{N}$ & Number of layers in the encoder (i.e., encoder depth) \\
$N_d \in \mathbb{N}$ & Number of layers in the decoder (i.e., decoder depth) \\
$n_{f_r}$ (resp. $n_{f_c}$) & Number of rows (resp. columns) in the feature data for each sample \\
$f \in \mathbb{R}^{n_{f_r} \times n_{f_c} \times n}$ & Encoder (data) variable \\
$g \in \mathbb{R}^{n_{f_r} \times n_{f_c} \times n}$ & Decoder (data) variable \\
$K \in \mathbb{R}^{n_{f_r} \times n_{f_r}}$ & Linear operator for the encoder (distinct for each layer) \\
$\tilde{K} \in \mathbb{R}^{n_{f_r} \times n_{f_r}}$ & Linear operator for the decoder (distinct for each layer) \\
$b \in \mathbb{R}$ & Bias (distinct for each encoder layer) \\
$\tilde{b} \in \mathbb{R}$ & Bias (distinct for each decoder layer) \\
$P \in \mathbb{R}^{n_{f_r} \times n_{f_c} \times n}$ & Lagrange multiplier for encoder \\
$\tilde{P} \in \mathbb{R}^{n_{f_r} \times n_{f_c} \times n}$ & Lagrange multiplier for decoder \\
$\tau \in \mathbb{R}$ & Temporal step size for encoder and decoder \\
$\sigma(\cdot)$ & Pointwise activation function for encoder \\
$\tilde{\sigma}(\cdot)$ & Pointwise activation function for decoder \\
$\hat{\sigma}(\cdot)$ & Pointwise activation function for output layer \\
$(\cdot)'$ & Derivative with respect to argument \\
$\mathrm{tr}(\cdot)$ & Trace operator \\
$(\cdot)^{\intercal}$ & Matrix transpose \\
$\odot$ & Pointwise (Hadamard) multiplication \\
$m_1$ & Number of mini-batches used in training \\
$m_2$ & Number of iterations in gradient-based optimization \\
$\alpha_{\text{train}},\; \alpha_{\text{valid}},\; \alpha_{\text{test}}$ & Reconstruction error for training, validation, and test data \\
\hline
\end{tabular}
\end{table}

\subsection{Rank-adaptive Euler Scheme \cite{Rogers_2021}\label{RAEuler}}
Consider the initial value problem 
\begin{equation}\label{IVP}
\frac{dY(t)}{dt} = \mathcal{N}(Y(t)), 
\qquad Y(0)=Y_0.
\end{equation}
Here, $Y:[0,T]\rightarrow {\mathbb R}^{n_1\times n_2 \times \cdots \times n_d}$, with $d\geq 2$, is a multi-dimensional array of real numbers (the solution tensor), and $\mathcal{N}$ is a tensor-valued nonlinear map.

The rank-adaptive forward Euler scheme for Euler method is given by
\begin{equation}
\label{FwdEuler_adap}
Y_{j+1} = {\mathfrak{T}}_{r_j}(Y_j + \tau \: {\mathfrak{T}}_{s_j}(\mathcal{N}( Y_j))).
\end{equation}
Where ${\mathfrak{T}}_{s_j}$ and ${\mathfrak{T}}_{r_j}$ are the rank adaptive SVD step truncation operators. The ranks $s_j$ and $r_j$ are selected so that, with $M_s\in \RR_+$, and $M_r\in \RR_+$,
\begin{subequations}\label{Euler_bounds}
\begin{align}
&\left\| \mathcal{N}(Y_j)- {\mathfrak{T}}_{s_j}(\mathcal{N}( Y_j)) \right\|_2 \leq M_s \tau ,\label{eq:s_bound} \\
&\left\|Y_j + \tau\:{\mathfrak{T}}_{s_j}(\mathcal{N}(Y_j))- {\mathfrak{T}}_{r_j}(Y_j + \tau \:{\mathfrak{T}}_{s_j}(\mathcal{N}( Y_j))) \right\|_2 \leq M_r \tau^{2},\label{eq:r_bound}
\end{align}
\end{subequations}
for all $j = 1,2,\ldots$, yield an order  one local truncation error for 
\eqref{FwdEuler_adap}. These conditions ensure that the numerical scheme is stable and consistent.

We extend this to a terminal value problem with a simple change of variable $t \mapsto T-t=\hat{t}$ for the reverse direction. Then, for the terminal value problem,
\begin{equation}\label{TVP}
-\frac{dZ(\hat{t})}{d\hat{t}} = \mathcal{N}(Z(\hat{t})), 
\qquad Z(T)=Z_T,
\end{equation} 
if we let $Z_j = Z(\hat{t}_j)$, we get the following rank-adaptive forward Euler scheme with the reverse time direction,
\begin{equation}
\label{FwdEulerRev_adap}
Z_{j} =
{\mathfrak{T}}_{r_{j+1}}(Z_{j+1} + \tau \: {\mathfrak{T}}_{s_{j+1}}(\mathcal{N}( Z_{j+1}))).
\end{equation}
Where ${\mathfrak{T}}_{s_{j+1}}$ and ${\mathfrak{T}}_{r_{j+1}}$ are the truncation operators as described above, and ranks are selected so that,
\begin{subequations}\label{Euler_bounds_Rev}
\begin{align}
&\left\| \mathcal{N}(Z_{j+1})- {\mathfrak{T}}_{s_{j+1}}(\mathcal{N}(Z_{j+1})) \right\|_2 \leq M_s \tau ,\label{eq:s_bound_rev} \\
&\left\|Z_{j+1} + \tau\:{\mathfrak{T}}_{s_{j+1}}(\mathcal{N}(Z_{j+1}))- {\mathfrak{T}}_{r_{j+1}}(Z_{j+1} + \tau \:{\mathfrak{T}}_{s_{j+1}}(\mathcal{N}( Z_{j+1}))) \right\|_2 \leq M_r \tau^{2},\label{eq:r_bound_rev}
\end{align}
\end{subequations}
for all $j =N-1, N-2,\ldots,0$, yield an order one local truncation error for 
\eqref{FwdEulerRev_adap}.
  
\section{Continuous Deep Autoencoder in an Optimal Control Framework \label{s:AE_OC}}

\subsection{Classical Autoencoder}
An autoencoder is a special type of feedforward neural network where the goal is to learn the composition of functions  $g(f(\hat{x}))$ for the input data $\hat{x}\in \mathcal{X}$. The encoder function $f:\mathcal{X} \to \mathcal{H}$ projects the data from higher-dimensional space $\mathcal{X}$ to lower-dimensional space $\mathcal{H}$, and the decoder function $g:\mathcal{H} \to \mathcal{X}$ takes it back to the original space. The learning process is described by minimizing a loss function $J(x,g(f(\hat{x}))$, given reference data $x\in\mathcal{X}$, which penalizes $g(f(\hat{x}))$ for being dissimilar from $x$, such as mean squared error \cite[Chapter 14]{Goodfellow2016}.

For a single layer network, the encoder and decoder functions are defined as feedforward units, 
\[
\begin{aligned}
f &= \sigma(K\hat{x}+b),\\
g &= \tilde{\sigma}(\tilde{K}f+\tilde{b}),
\end{aligned}
\]
where the image $f$ is known as the code, latent variable, or latent representation, $K, \tilde{K}, b$ and $\tilde{b}$ are the unknown weights and biases, and $\sigma$ is the activation function. If $f$ and $g$ are computed in multiple layers, then this becomes a deep autoencoder. For learning the underlying low-dimensional structure of the input data, $g(f(\hat{x}))$ is the identity map (along with $x=\hat{x}$). In the case of image denoising (resp. deblurring), $g(f(\hat{x}))$ is the denoising (resp. deblurring) operator, $x$ is the clean image and $\hat{x}$ is the noisy (resp. blurred) image. 

However, when practitioners use autoencoders, it is not clear how to select the depth and the level of compression in the network. In the next section, we address this by considering above mentioned basic notion of an autoencoder, and leverage modeling tools from optimal control theory to formulate an autoencoder model as an optimal control problem. Note that this formulation allows us to use analysis tools from differential equation theory as well, which enable us to make informed decisions regarding the depth and compression of the autoencoder architecture. 

\subsection{Continuous Deep Autoencoder as an Optimal Control Problem}

In this work, we introduce a continuous-time deep autoencoder architecture through an optimal control framework \cite{Ruthotto_2020_DL_OC, Khatri_2020_fdnn}. Autoencoders aim to learn mappings from high-dimensional data to a low-dimensional latent space and back \cite[Chapter~14]{Goodfellow2016}, typically using an encoder-decoder structure. We model the encoder and decoder as coupled ordinary differential equations, and formulate the learning task as the minimization of a cost functional $J$, subject to these dynamics. The goal is to learn the time-dependent weights and biases governing the evolution. Compression is introduced via rank reduction at the discrete level later in the paper.

While many autoencoder variants exist—sparse, denoising, contractive, residual, etc. \cite{Goodfellow2016, Burger_2021_Arch_search} (typically differing in the objective function)—our formulation accommodates a general class of regularized deep autoencoders, allowing flexibility in the loss objective and regularization terms.

Consider the reference data $x\in\RR^{n_{f_r}\times n_{f_c} \times n}$, state variables $f,g \in \RR^{n_{f_r}\times n_{f_c} \times n}$. For simplicity, we restrict ourselves to 3D tensors. Let $J(x,g(T))$ be the objective (loss) to be minimized,
$\mathcal{R}(K(t),\tilde{K}(t),b(t),\tilde{b}(t))$ be the regularizer, and let $d_t y(t) := \frac{dy(t)}{dt}$ be the standard time derivative. Moreover, we denote by 
\[
	\Theta = \left\{ \big(K(t),b(t)\big), \, \left(\tilde{K}(t),\tilde{b}(t)\right) \right\},
\]
the unknown weight and biases. Then, for $0<\bar{t}\leq T$,
\begin{equation}\label{AE_OC}
\begin{aligned}
\min_{\Theta} \; & 
\left\{ \mathcal{J}(\Theta, g(T), x) := J(x, g(T)) + \mathcal{R}(\Theta) \right\} \\
\text{subject to} \quad &
\left\{
\begin{aligned}
d_t f(t) &= \sigma\big(K(t) f(t) + b(t)\big), && \forall\, t \in (0, \bar{t}] \\
f(0) &= \hat{x}, \\
d_t g(t) &= \tilde{\sigma}\big(\tilde{K}(t) g(t) + \tilde{b}(t)\big), && \forall\, t \in (\bar{t}, T] \\
g(\bar{t}) &= f(\bar{t}),
\end{aligned}
\right.
\end{aligned}
\end{equation}

Here, the differential equations for $f$ and $g$ corresponds to encoder and decoder, respectively.

Furthermore, $K,\tilde{K} \in \RR^{n_{f_r}\; \times\; n_{f_r}}$ are the linear operators and $b,\tilde{b} \in \RR$ are the biases, and $\sigma$ and $\tilde{\sigma}$ are nonlinear pointwise activation functions. Note that $g(\bar{t}) = f(\bar{t})$ is the latent state of the encoder and decoder variables. We denote the initial data for the encoder differential equation as $\hat{x}\in \RR^{n_{f_r}\times n_{f_c} \times n}$. 

\begin{remark}
\rm
In the identity mapping case, we have $\hat{x} = x$. For image denoising and deblurring, $\hat{x} $ denotes the noisy or blurred input, respectively, while $x$ is the corresponding clean image. In classification tasks, $\hat{x}$ represents the input data and $x$ the true label, with the loss functional $J(x,g(T))$ taken as the cross-entropy loss function as in \cite{Khatri_2020_fdnn}. 
\end{remark}

To solve the optimal control problem in \cref{AE_OC}, we reformulate it using a Lagrangian functional and derive the corresponding optimality conditions. Let $ P $ and $\tilde{P}$ denote the Lagrange multipliers associated with the encoder and decoder states $f$ and $g$, respectively. The Lagrangian is given by:

\begin{equation}\label{LagDef}
\begin{aligned}
\mathcal{L}(f, g, \Theta; P, \tilde{P}) :=\; & \mathcal{J}(\Theta, g(T), x) \\
& + \left\langle d_t f(t) - \sigma\big(K(t) f(t) + b(t)\big),\, P(t) \right\rangle_{(0,\bar{t})} \\
& + \left\langle d_t g(t) - \tilde{\sigma}\big(\tilde{K}(t) g(t) + \tilde{b}(t)\big),\, \tilde{P}(t) \right\rangle_{(\bar{t},T)}
\end{aligned}
\end{equation}
where, $\langle \cdot,\cdot\rangle_{(t_1,t_2)} := \int_{t_1}^{t_2}\langle \cdot,\cdot\rangle_F \;dt$ is the $L^2$-inner product and $\langle \cdot,\cdot\rangle_F$ is the Frobenius inner product. Next, we expand the inner product and then apply integration by parts. Recall, $f(0) = \hat{x}$ and $g(\bar{t}) = f(\bar{t})$. This yields the following simplified Lagrangian functional,

\begin{equation} \label{Lag}
\begin{aligned}
\mathcal{L}(f, g, \Theta; P, \tilde{P}) =\; & \mathcal{J}(\Theta, g(T), x) 
- \left\langle d_t P(t), f(t) \right\rangle_{(0,\bar{t})} 
+ \left\langle P(\bar{t}), f(\bar{t}) \right\rangle_F \\
& - \left\langle P(0), \hat{x} \right\rangle_F 
- \left\langle P(t), \sigma\big(K(t) f(t) + b(t)\big) \right\rangle_{(0,\bar{t})} \\
& - \left\langle d_t \tilde{P}(t), g(t) \right\rangle_{(\bar{t},T)} 
+ \left\langle \tilde{P}(T), g(T) \right\rangle_F 
- \left\langle \tilde{P}(\bar{t}), f(\bar{t}) \right\rangle_F \\
& - \left\langle \tilde{P}(t), \tilde{\sigma}\big(\tilde{K}(t) g(t) + \tilde{b}(t)\big) \right\rangle_{(\bar{t},T)}.
\end{aligned}
\end{equation}

Let $(f^*,g^*,\Theta^*;P^*,\tilde{P}^*)$ denote a stationary point, then the first order necessary optimality conditions are given by the following state, adjoint and design equations: 
\begin{enumerate}
\smallskip
\item[(A)] \textbf{System of state equations.} The gradient of $\mathcal{L}(\cdot)$ with respect to $(P,\tilde{P})$ evaluated at the stationary point $(f^*,g^*,\Theta^*;P^*,\tilde{P}^*)$ yields $\nabla_{(P,\tilde{P})}\mathcal{L}(f^*,g^*,\Theta^*;P^*,\tilde{P}^*) = 0$,  
equivalently,

\begin{equation}\label{stateEq}
\left\{
\begin{aligned}
& d_t f^*(t) &= \; &\sigma(K^*(t)f^*(t)+b^*(t)), \quad &\forall \; t \in (0,\bar{t}]\\
&f^*(0) &= \; &\hat{x} & \\
& d_t g^*(t) &=\;&\tilde{\sigma}(\tilde{K}^*(t)g^*(t)+\tilde{b}^*(t)), \quad &\forall \; t \in (\bar{t},T]\\
&g^*(\bar{t}) & =\; &f^*(\bar{t}) &
\end{aligned}
\right.
\end{equation}
For the state variables $f^*$ and $g^*$, we solve \cref{stateEq} forward in time, therefore we call the system of equations \cref{stateEq} \textit{forward propagation} in the proposed autoencoder.

\smallskip
\item[(B)] \textbf{System of adjoint equations.} Next, the gradient of $\mathcal{L}(\cdot)$ with respect to $(f, g)$ evaluated at the stationary point $(f^*,g^*,\Theta^*;P^*,\tilde{P}^*)$ yields $\nabla_{(f,g)} \mathcal{L}(f^*,g^*,\Theta^*;P^*,\tilde{P}^*) = 0$, equivalently, 

\begin{equation}\label{AdjEq}
\left\{
\begin{aligned}
-&d_tP^*(t) &=\; &K^*(t)^\intercal \left(P^*(t) \odot \sigma^\prime(K^*(t)f^*(t)+b^*(t))\right),\quad &\forall \; t \in [0,\bar{t})\\
& P^*(\bar{t}) &= \; &\tilde{P}^*(\bar{t})& \\
-&d_t\tilde{P}^*(t) &=\;& \tilde{K}^*(t)^\intercal\left(\tilde{P}^*(t) \odot (\tilde{\sigma}^\prime(\tilde{K}^*(t)g^*(t)+\tilde{b}^*(t))\right), &\forall \; t \in [\bar{t},T)\\
&\tilde{P}^*(T)&= \; &- \nabla_{g} J(x,g^*(T))
\end{aligned}
\right.
\end{equation}
Notice that the adjoint variables $P^*$ and $\tilde{P}^*$, with their terminal conditions, are obtained by marching backward in time. As a result, the system of equations \cref{AdjEq} is called \textit{backward propagation} in the proposed autoencoder.

\item[(C)] \textbf{Gradient.} 

Finally, evaluating $\nabla_{\Theta}\mathcal{L}(\cdot)$ at $(f^*,g^*,\Theta^*;P^*,\tilde{P}^*)$ yield the following gradients w.r.t. the design variables,

\begin{equation}\label{DesEq}
\begin{aligned}
\nabla_{K} \mathcal{L}(f^*, g^*, \Theta^*; P^*, \tilde{P}^*) &= 
- f^*(t) \left(P^*(t) \odot \sigma'\big(K^*(t) f^*(t) + b^*(t)\big)\right)^{\intercal} 
+ \nabla_{K} \mathcal{R}(\Theta^*), \quad && \forall\, t \in (0, \bar{t}) \\[0.5em]
\nabla_{b} \mathcal{L}(f^*, g^*, \Theta^*; P^*, \tilde{P}^*) &= 
- \left\langle \sigma'\big(K^*(t) f^*(t) + b^*(t)\big), P^*(t) \right\rangle_F 
+ \nabla_{b} \mathcal{R}(\Theta^*), \quad && \forall\, t \in (0, \bar{t}) \\[0.5em]
\nabla_{\tilde{K}} \mathcal{L}(f^*, g^*, \Theta^*; P^*, \tilde{P}^*) &= 
- g^*(t) \left(\tilde{P}^*(t) \odot \tilde{\sigma}'\big(\tilde{K}^*(t) g^*(t) + \tilde{b}^*(t)\big)\right)^{\intercal} 
+ \nabla_{\tilde{K}} \mathcal{R}(\Theta^*), \quad && \forall\, t \in (\bar{t}, T) \\[0.5em]
\nabla_{\tilde{b}} \mathcal{L}(f^*, g^*, \Theta^*; P^*, \tilde{P}^*) &= 
- \left\langle \tilde{\sigma}'\big(\tilde{K}^*(t) g^*(t) + \tilde{b}^*(t)\big), \tilde{P}^*(t) \right\rangle_F 
+ \nabla_{\tilde{b}} \mathcal{R}(\Theta^*), \quad && \forall\, t \in (\bar{t}, T).
\end{aligned}
\end{equation}

\end{enumerate}
In view of $(\text{A})-(\text{C})$, we can use a gradient-based solver to find a stationary point of \cref{AE_OC}.


\section{Discrete Deep Autoencoder as an Optimal Control Problem \label{s:Disc_AE_OC}}

We adopt the \emph{optimize-then-discretize} paradigm for solving the optimal control problem. The first-order optimality conditions for the continuous formulation in \cref{AE_OC} are given by the state equations \cref{stateEq}, the adjoint equations \cref{AdjEq}, and the design equations \cref{DesEq}. Our goal is to discretize these conditions, which necessitates a numerical scheme for the differential systems \cref{stateEq}–\cref{AdjEq} and a time-discretization for the design equations \cref{DesEq}.

To this end, we employ an explicit rank-adaptive Euler scheme, as described in \cref{RAEuler}, for the forward and adjoint dynamics. The use of a rank-adaptive integration method is motivated by the need to dynamically capture both compression and expansion in the data representation. This mechanism underlies the resulting characteristic butterfly-shaped structure of the proposed autoencoder.

\subsection{Discrete Optimality Conditions}

We begin by establishing the correspondence between discrete time-stepping in the differential equation framework and layer-wise propagation in deep neural networks. To this end, we uniformly discretize the time interval $[0, T]$ with step size $\tau = T/N$, yielding the partition
\[
0 = t_0 < t_1 < t_2 < \cdots < t_{N_e} = \bar{t} < t_{N_e+1} < \cdots < t_N = T,
\]
where each time point is given by \( t_j = j \tau \) for \( j = 0, \ldots, N \). The index \( j \) thus serves as the layer index in the network, with \( N_e \) and \( N_d \) denoting the number of encoder and decoder layers, respectively. Assuming a symmetric architecture, we take \( N_d = N - N_e = N/2 \), where \( N \) is assumed even.

To apply the tensor-based numerical scheme, we represent the state and adjoint variables as tensors, denoted by \( \mathfrak{f}, \mathfrak{g}, \mathfrak{P}, \tilde{\mathfrak{P}} \), respectively. The design variables (e.g., weight matrices and biases) remain in their original matrix or scalar form. This formulation transforms the differential equations into tensor-valued nonlinear dynamical systems.

For temporal discretization of the state and adjoint equations, we employ the rank-adaptive explicit Euler schemes introduced in \cref{RAEuler}. Specifically, we use the forward Euler method for initial value problems \cref{FwdEuler_adap} to discretize the state equations, and its reverse-time counterpart for terminal value problems \cref{FwdEulerRev_adap} to discretize the adjoint equations. This results in the following discrete optimality system, starting from the tensor-valued initial condition \( \hat{\mathfrak{x}} \).

\begin{enumerate}
\smallskip
\item[(a)] \textbf{Discrete State Equations.} Forward propagation:

\begin{equation}\label{stateEq_disc}
\left\{
\begin{aligned}
\mathfrak{f}^*(t_j) &= 
\mathfrak{T}_{r_j} \left( 
\mathfrak{f}^*(t_{j-1}) 
+ \tau\, \mathfrak{T}_{s_j} \left( 
\sigma\big(K^*(t_{j-1})\, \mathfrak{f}^*(t_{j-1}) + b^*(t_{j-1})\big) 
\right) 
\right), 
&& j = 1, \dots, N_e \\[0.5em]
\mathfrak{f}^*(t_0) &= \hat{\mathfrak{x}}, \\[0.5em]
\mathfrak{g}^*(t_j) &= 
\mathfrak{T}_{r_j} \left( 
\mathfrak{g}^*(t_{j-1}) 
+ \tau\, \mathfrak{T}_{s_j} \left( 
\tilde{\sigma}\big(\tilde{K}^*(t_{j-1})\, \mathfrak{g}^*(t_{j-1}) + \tilde{b}^*(t_{j-1})\big) 
\right) 
\right), 
&& j = N_e + 1, \dots, N \\[0.5em]
\mathfrak{g}^*(t_{N_e}) &= \mathfrak{f}^*(t_{N_e}).
\end{aligned}
\right.
\end{equation}
\smallskip
\item[(b)] \textbf{Discrete Adjoint Equations.} Backward propagation:
\begin{equation}\label{AdjEq_disc} 
\left\{
\begin{aligned}
&\mathfrak{P}^*(t_j) &=\;&{\mathfrak{T}}_{r_{j+1}}\left(\mathfrak{P}^*(t_{j+1}) +\tau {\mathfrak{T}}_{s_{j+1}}\left(K^*(t_{j})^\intercal \left(\mathfrak{P}^*(t_{j+1}) \odot \sigma^\prime(K^*(t_{j})\mathfrak{f}^*(t_{j+1})+b^*(t_{j}))\right)\right)\right),\\
& & &\hspace{10cm} j = N_{e}-1,\dots,0\\
& \mathfrak{P}^*(t_{N_e}) &= \; &\tilde{\mathfrak{P}}^*(t_{N_e})& \\
& \tilde{\mathfrak{P}}^*(t_j) &=\;&{\mathfrak{T}}_{r_{j+1}}\left(\tilde{\mathfrak{P}}^*(t_{j+1}) + \tau {\mathfrak{T}}_{s_{j+1}}\left( \tilde{K}^*(t_j)^\intercal\left(\tilde{\mathfrak{P}}(t_{j+1}) \odot (\tilde{\sigma}^\prime(\tilde{K}^*(t_j)\mathfrak{g}^*(t_{j+1})+\tilde{b}^*(t_j))\right)\right)\right),\\
& & &\hspace{10cm} j = N-1,\dots,N_e\\
&\tilde{\mathfrak{P}}^*(t_N)&= \; &- \nabla_{\mathfrak{g}} J(\mathfrak{x},\mathfrak{g}^*(t_N))
\end{aligned}
\right.
\end{equation}
\smallskip
\item[(c)] \textbf{Discrete Gradient.} 

\begin{equation}\label{DesEq_disc}
\begin{aligned}
\nabla_{K} \mathcal{L}(\mathfrak{f}^*, \mathfrak{g}^*, \Theta^*; \mathfrak{P}^*, \tilde{\mathfrak{P}}^*) 
&= - \mathfrak{f}^*(t_j) 
\left( \mathfrak{P}^*(t_{j+1}) \odot \sigma'\big(K^*(t_j) \mathfrak{f}^*(t_j) + b^*(t_j)\big) \right)^{\intercal} 
+ \nabla_{K} \mathcal{R}(\Theta^*(t_j)), \\
&\hspace{8cm} j = 1, \dots, N_e \\[0.75em]
\nabla_{b} \mathcal{L}(\mathfrak{f}^*, \mathfrak{g}^*, \Theta^*; \mathfrak{P}^*, \tilde{\mathfrak{P}}^*) 
&= - \left\langle \sigma'\big(K^*(t_j) \mathfrak{f}^*(t_j) + b^*(t_j)\big), \mathfrak{P}^*(t_{j+1}) \right\rangle_F 
+ \nabla_{b} \mathcal{R}(\Theta^*(t_j)), \\
&\hspace{8cm} j = 1, \dots, N_e \\[0.75em]
\nabla_{\tilde{K}} \mathcal{L}(\mathfrak{f}^*, \mathfrak{g}^*, \Theta^*; \mathfrak{P}^*, \tilde{\mathfrak{P}}^*) 
&= - \mathfrak{g}^*(t_j) 
\left( \tilde{\mathfrak{P}}^*(t_{j+1}) \odot \tilde{\sigma}'\big(\tilde{K}^*(t_j) \mathfrak{g}^*(t_j) + \tilde{b}^*(t_j)\big) \right)^{\intercal} 
+ \nabla_{\tilde{K}} \mathcal{R}(\Theta^*(t_j)), \\
&\hspace{8cm} j = N_e + 1, \dots, N \\[0.75em]
\nabla_{\tilde{b}} \mathcal{L}(\mathfrak{f}^*, \mathfrak{g}^*, \Theta^*; \mathfrak{P}^*, \tilde{\mathfrak{P}}^*) 
&= - \left\langle \tilde{\sigma}'\big(\tilde{K}^*(t_j) \mathfrak{g}^*(t_j) + \tilde{b}^*(t_j)\big), \tilde{\mathfrak{P}}^*(t_{j+1}) \right\rangle_F  + \nabla_{\tilde{b}} \mathcal{R}(\Theta^*(t_j)), \\
&\hspace{8cm} j = N_e + 1, \dots, N.
\end{aligned}
\end{equation}
\end{enumerate}

Based on the optimality conditions $(\text{a})$–$(\text{c})$, we proceed to develop a gradient-based method to solve the original problem \cref{AE_OC}. We emphasize that, under the optimal control formulation, each evaluation of the gradient in \cref{DesEq_disc} requires solving one forward (state) and one backward (adjoint) system.

For notational simplicity, in the algorithms presented below, we omit the asterisk $(\cdot)^*$ that denotes stationary points. Furthermore, we adopt the convention $ u_i := u(t_i) $ for all $ i $, where variables are indexed by discrete time steps.


\section{OCTANE\label{s:OCTANE} }
In this section, we describe the algorithmic structure of our gradient-based learning method (OCTANE), curated using optimality conditions  \cref{stateEq_disc,AdjEq_disc,DesEq_disc}. We begin by presenting the strategy we have designed to handle the truncation operators and corresponding rank selections in \cref{rank_select}, followed by the algorithmic architecture of OCTANE in \cref{octane_alg}. \Cref{comp_cost} discusses the computational cost and memory needs.

\subsection{Truncation operator and rank selection strategy}\label{rank_select}

We outline our pseudo code strategy, incorporated into the main algorithm in \cref{alg:OCTANE_train}.
\begin{itemize}
  \item[(a)] \textbf{Rank-adaptive truncation:}
  \begin{itemize}[label={\scriptsize$\bullet$}]
    \item Recall from \cref{RAEuler}, \( \mathfrak{T}_{r_j} \) and \( \mathfrak{T}_{s_j} \) are adaptive SVD-based truncation operators.
    \item Truncation ranks \( r_j \) and \( s_j \) are chosen to satisfy the conditions in \cref{Euler_bounds} (for forward integration) and \cref{Euler_bounds_Rev} (for reverse-time adjoints).
    \item Tensor rounding routines in libraries such as \texttt{TT-toolbox} typically require:
    \begin{enumerate}
      \item[(i)] a maximum rank threshold, and
      \item[(ii)] an optional approximation tolerance (defaulting to machine precision).
    \end{enumerate}
    \item These inputs guide low-rank approximation during truncation.
  \end{itemize}

 \item[(b)] \textbf{Encoder rank selection:}
  \begin{itemize}[label={\scriptsize$\bullet$}]
    \item During forward integration of the encoder state equations:
    \begin{itemize}
      \item[(i)] Apply \( \mathfrak{T}_{s_j} \) to the inner term in \cref{FwdEuler_adap}, using the full rank of the input as the prescribed maximum threshold.
      \item[(ii)] Apply \( \mathfrak{T}_{r_j} \) to the full update, again using the full rank as the upper bound.
    \end{itemize}
    \item Store the resulting encoder ranks in a vector \( \mathbf{r}_e \) (not necessarily equal to the prescribed ranks).
  \end{itemize}

  \item[(c)] \textbf{Decoder rank selection:}
  \begin{itemize}[label={\scriptsize$\bullet$}]
    \item Define \( \mathbf{r}_d := \text{flip}(\mathbf{r}_e) \) to impose symmetry across the encoder-decoder interface.
    \item Use \( \mathbf{r}_d \) as the prescribed threshold rank vector for decoder forward integration, ensuring that \cref{Euler_bounds} holds. Toward the last layer, this procedure may return to the full-rank representation of the original input.
  \end{itemize}

  \item[(d)] \textbf{Adjoint equations:}
  \begin{itemize}[label={\scriptsize$\bullet$}]
    \item Reuse \( \mathbf{r}_e \) and \( \mathbf{r}_d \) as the maximum thresholds for the encoder and decoder adjoint equations, respectively.
    \item Ensure that truncation during reverse-time integration satisfies \cref{Euler_bounds_Rev}.
  \end{itemize}

  \item[(e)] \textbf{Autoencoder shape and rank profile:}
  \begin{itemize}[label={\scriptsize$\bullet$}]
    \item The resulting vectors of encoder and decoder ranks in forward integration, selected by the truncation operators, define the discovered autoencoder architecture.
    \item A symmetric structure (\( N_e = N_d \)) ensures one-to-one correspondence of rank profiles across the encoder and decoder.
    \item The resultant ranks are not learned via optimization; they emerge naturally from the dynamics of the problem and the input data.
  \end{itemize}

\end{itemize}
 This procedure is implemented in \cref{alg:OCTANE_train}. Note that we do not track the inner ranks \( s_j \) explicitly. We further remark that tensor rounding is the key enabler for tractable training in this framework.
 Tensor decompositions (e.g., TT, Tucker, HOSVD) are not unique and induce inherent truncation 
 error bounds. The numerical scheme in \cref{RAEuler} is agnostic to the chosen tensor format.
 Once a decomposition is fixed, its truncation error can be incorporated into the overall numerical error control strategy, as illustrated in \cref{Error_merge}.

\subsection {OCTANE Algorithm}\label{octane_alg}
OCTANE is an autoencoder training framework designed to uncover the low-rank structure of data while solving the associated learning problem in an optimal control framework. It naturally fits within unsupervised learning, particularly when learning the identity map. However, it can be readily adapted to supervised tasks (e.g., classification or regression) by replacing the reference data $x$ with labeled targets and selecting an appropriate loss function.

\begin{figure}
\includegraphics[width=.47\linewidth]{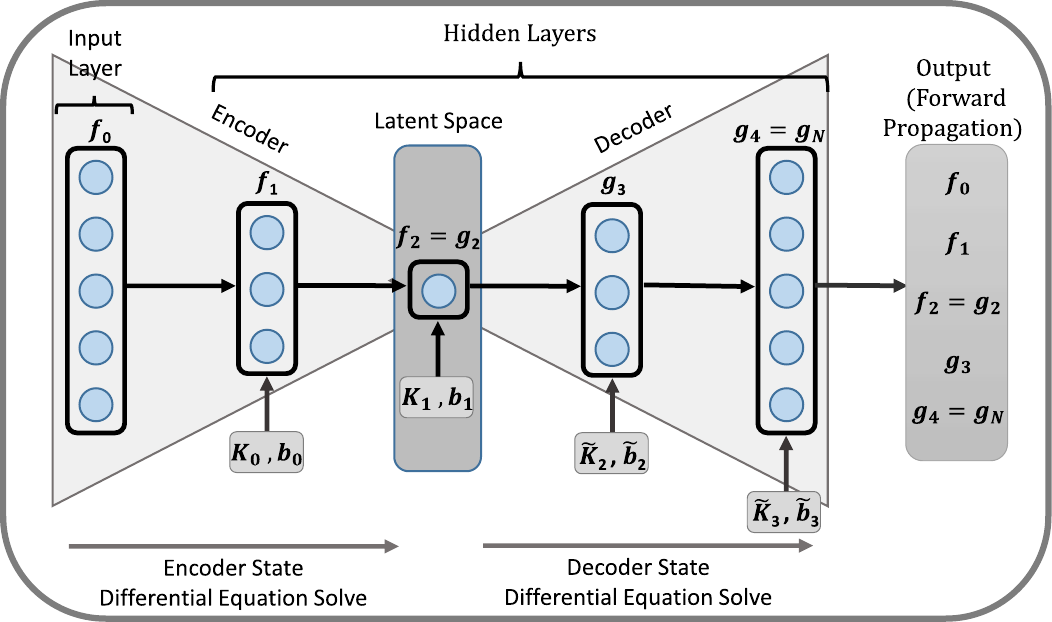}
\includegraphics[width=.47\linewidth]{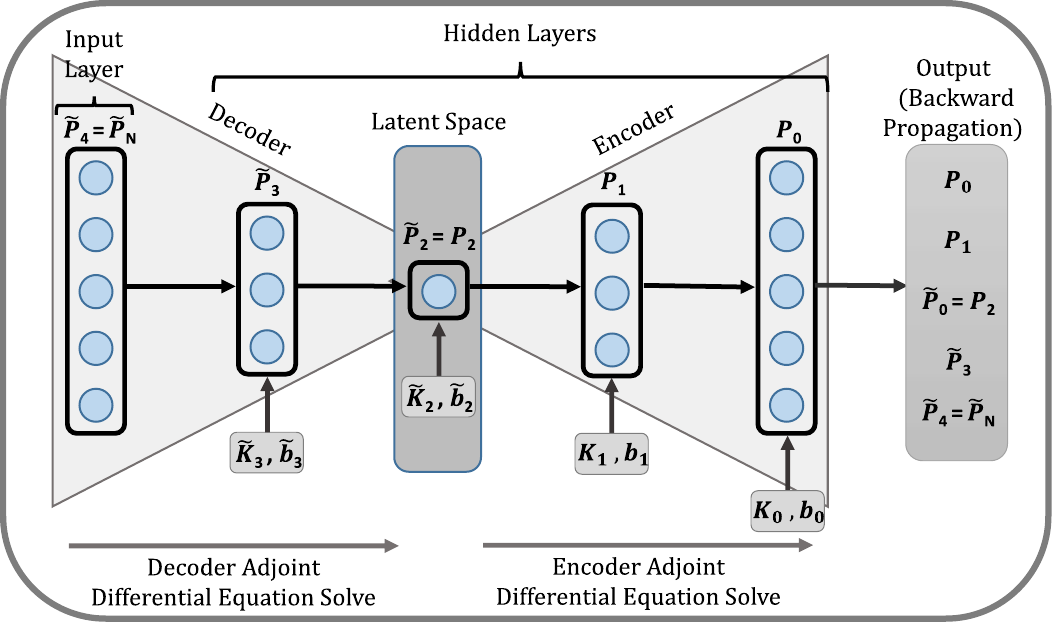}
\caption{Illustration of forward (left) and backward (right) propagation in a 4-layer OCTANE (i.e. 2 encoder layers and 2 decoder layers). The differential equations are solved using tensor-based, rank-adaptive explicit forward Euler scheme. Note that the reduction in layer size is by virtue of layer rank.}\label{f:OCTANEnet}
\end{figure}

OCTANE, as a deep learning model, includes training and testing phases. The training involves:
\begin{itemize}
\item \textbf{Forward propagation:} solving the encoder and decoder state equations;
\item \textbf{Backward propagation:} solving the corresponding adjoint equations;
\item \textbf{Gradient update:} evaluating the design equations.
\end{itemize}
These steps collectively learn the network parameters (design variables), yielding the trained autoencoder. During testing, only forward propagation is performed to reconstruct unseen input data. 
A schematic of forward and backward propagation in a sample 4-layer OCTANE network (2 encoder + 2 decoder layers) is shown in \cref{f:OCTANEnet}.

The algorithmic structure begins with forward propagation, comprising the encoder and decoder routines in \cref{alg:Bnet_fwd_encoder} and \cref{alg:Bnet_fwd_decoder}. This is followed by backward propagation, defined by the decoder and encoder adjoint routines in \cref{alg:Bnet_bckwd_decoder} and \cref{alg:Bnet_bckwd_encoder}. The complete training and testing procedures are presented in \cref{alg:OCTANE_train} and \cref{alg:OCTANE_test}, respectively.

\begin{algorithm}[!htb]
\caption{Forward Encoder}
\label{alg:Bnet_fwd_encoder}
\begin{algorithmic}[1]
\Require $\{\text{lower index } li, \text{upper index } ui\},\; \mathfrak{f}_{li-1},\; \{K_j, b_j\}_{j=li-1}^{ui-1},\; \tau,\; M_s,\; M_r$
\Ensure $\{\mathfrak{f}_j\}_{j=li-1}^{ui},\; \{r_j\}_{j=li-1}^{ui}$

\State Compute and store the actual rank of $\mathfrak{f}_{li-1}$: \quad $r_{li-1} = \text{rank}(\mathfrak{f}_{li-1})$
\For{$j = li,\dots,ui$}
    \State Compute the term $\mathfrak{u}$: \quad $\mathfrak{u} = \sigma(K_{j-1} \mathfrak{f}_{j-1} + b_{j-1})$
    
    \State Truncate the tensor $\mathfrak{u}$ with tolerance $\varepsilon$ and maximum rank $r_{li-1}$:
    \Statex \hspace{\algorithmicindent}
    $\mathfrak{u} = \mathfrak{T}_{s_j}(\mathfrak{u})$
    \Comment{Choose smallest $s_j > 0$ s.t. \eqref{eq:s_bound} holds with $M_s$. If not found, \textbf{break}}
    
    \State Update $\mathfrak{f}_j$: \quad $\mathfrak{f}_j = \mathfrak{f}_{j-1} + \tau \mathfrak{u}$
    
    \State Truncate the tensor $\mathfrak{f}_j$ with tolerance $\varepsilon$ and maximum rank $r_{li-1}$:
    \Statex \hspace{\algorithmicindent}
    $\mathfrak{f}_j = \mathfrak{T}_{r_j}(\mathfrak{f}_j)$
    \Comment{Choose smallest $r_j > 0$ s.t. \eqref{eq:r_bound} holds with $M_r$. If not found, \textbf{break}}

    \State Compute and store actual rank of $\mathfrak{f}_j$: \quad $r_j = \text{rank}(\mathfrak{f}_j)$
\EndFor
\end{algorithmic}
\end{algorithm}

\begin{algorithm}[!htb]
\caption{Forward Decoder}
\label{alg:Bnet_fwd_decoder}
\begin{algorithmic}[1]
\Require $\{\text{lower index } li, \text{upper index } ui\},\; \mathfrak{g}_{li-1},\; \{\tilde{K}_j, \tilde{b}_j\}_{j=li-1}^{ui-1},\; \{r_q\}_{q=0}^{li},\; \tau,\; M_s,\; M_r$
\Ensure $\{\mathfrak{g}_j\}_{j=li-1}^{ui},\; \{\tilde{r}_j\}_{j=li-1}^{ui}$

\State Compute and store the actual rank of $\mathfrak{g}_{li-1}$: \quad $\tilde{r}_{li-1} = \text{rank}(\mathfrak{g}_{li-1})$
\State Let prescribed ranks $\mathbf{s} = \mathbf{r}$ and set counter $q = 1$
\For{$j = li,\dots,ui$}
    \State Compute the term $\mathfrak{u}$: \quad $\mathfrak{u} = \sigma(\tilde{K}_{j-1} \mathfrak{g}_{j-1} + \tilde{b}_{j-1})$
    
    \State Truncate the tensor $\mathfrak{u}$ to the prescribed rank $\mathbf{s}_j$ with tolerance $\varepsilon$:
    \Statex \hspace{\algorithmicindent}
    $\mathfrak{u} = \mathfrak{T}_{\mathbf{s}_j}(\mathfrak{u})$
    \Comment{If \eqref{eq:s_bound} with $M_s$ fails, \textbf{break}}
    
    \State Update $\mathfrak{g}_j$: \quad
    $\mathfrak{g}_j = \mathfrak{g}_{j-1} + \tau \mathfrak{u}$
    
    \State Truncate the tensor $\mathfrak{g}_j$ to the prescribed rank $\mathbf{r}_j$ with tolerance $\varepsilon$:
    \Statex \hspace{\algorithmicindent}
    $\mathfrak{g}_j = \mathfrak{T}_{\mathbf{r}_j}(\mathfrak{g}_j)$
    \Comment{If \eqref{eq:r_bound} with $M_r$ fails, \textbf{break}}
    
    \State Compute and store actual rank of $\mathfrak{g}_j$: \quad $\tilde{r}_j = \text{rank}(\mathfrak{g}_j)$
    
    \State Increment counter: \quad $q = q + 1$
\EndFor
\end{algorithmic}
\end{algorithm}

\begin{algorithm}[!htb]
\caption{Backward Decoder}
\label{alg:Bnet_bckwd_decoder}
\begin{algorithmic}[1]
\Require $\{\text{lower index } li, \text{upper index } ui\},\; \tilde{\mathfrak{P}}_{ui+1},\; \{\mathfrak{g}_j\}_{j=li+1}^{ui+1},\; \{\tilde{K}_j, \tilde{b}_j\}_{j=li}^{ui},\; \{r_q\}_{q=0}^{li},\; \tau,\; M_s,\; M_r$
\Ensure $\{\tilde{\mathfrak{P}}_j\}_{j=li}^{ui+1},\; \{\tilde{r}_j\}_{j=li}^{ui+1}$

\State Compute and store actual rank of $\tilde{\mathfrak{P}}_{ui+1}$: \quad $\tilde{r}_{ui+1} = \text{rank}(\tilde{\mathfrak{P}}_{ui+1})$
\State Let prescribed ranks $\mathbf{s} = \mathbf{r}$ and initialize counter $q = li - 1$
\For{$j = ui,\dots,li$}
    \State Compute the update term: \quad
    $
    \mathfrak{u} = \tilde{K}_j^\intercal \left( \tilde{\mathfrak{P}}_{j+1} \odot \tilde{\sigma}^\prime(\tilde{K}_j \mathfrak{g}_{j+1} + \tilde{b}_j) \right)
    $
    
    \State Truncate $\mathfrak{u}$ to prescribed rank $s_q$ with tolerance $\varepsilon$: 
    \Statex \hspace{\algorithmicindent}
    $
    \mathfrak{u} = \mathfrak{T}_{\mathbf{s}_j}(\mathfrak{u})
    $
    \Comment{If \eqref{eq:s_bound} with $M_s$ fails, \textbf{break}}

    \State Update solution: \quad
    $
    \tilde{\mathfrak{P}}_j = \tilde{\mathfrak{P}}_{j+1} + \tau \mathfrak{u}
    $

    \State Truncate $\tilde{\mathfrak{P}}_j$ to prescribed rank $r_q$ with tolerance $\varepsilon$:
    \Statex \hspace{\algorithmicindent}
    $
    \tilde{\mathfrak{P}}_j = \mathfrak{T}_{\mathbf{r}_j}(\tilde{\mathfrak{P}}_j)
    $
    \Comment{If \eqref{eq:r_bound} with $M_r$ fails, \textbf{break}}

    \State Compute and store actual rank: \quad $\tilde{r}_j = \text{rank}(\tilde{\mathfrak{P}}_j)$
    \State Decrease counter: \quad $q = q - 1$
\EndFor
\end{algorithmic}
\end{algorithm}

\begin{algorithm}[!htb]
\caption{Backward Encoder}
\label{alg:Bnet_bckwd_encoder}
\begin{algorithmic}[1]
\Require $\{\text{lower index } li,\ \text{upper index } ui\},\; \mathfrak{P}_{ui+1},\; \{\mathfrak{f}_j\}_{j=li+1}^{ui+1},\; \{K_j, b_j\}_{j=li}^{ui},\; \{r_q\}_{q=0}^{ui},\; \tau,\; M_s,\; M_r$
\Ensure $\{\mathfrak{P}_j\}_{j=li}^{ui+1},\; \{\tilde{r}_j\}_{j=li}^{ui+1}$

\State Compute and store actual rank of $\mathfrak{P}_{ui+1}$: \quad $\tilde{r}_{ui+1} = \text{rank}(\mathfrak{P}_{ui+1})$
\State Let prescribed ranks $\mathbf{s} = \mathbf{r}$ and initialize counter $q = ui - 1$
\For{$j = ui,\dots,li$}
    \State Compute update term: \quad
    $
    \mathfrak{u} = K_j^\intercal \left( \mathfrak{P}_{j+1} \odot \sigma^\prime(K_j \mathfrak{f}_{j+1} + b_j) \right)
    $

    \State Truncate $\mathfrak{u}$ to prescribed rank $s_q$ with tolerance $\varepsilon$:
    \Statex \hspace{\algorithmicindent}
    $
    \mathfrak{u} = \mathfrak{T}_{\mathbf{s}_j}(\mathfrak{u})
    $
    \Comment{If \eqref{eq:s_bound} with $M_s$ fails, \textbf{break}}

    \State Update solution: \quad
    $
    \mathfrak{P}_j = \mathfrak{P}_{j+1} + \tau \mathfrak{u}
    $

    \State Truncate $\mathfrak{P}_j$ to prescribed rank $r_q$ with tolerance $\varepsilon$:
    \Statex \hspace{\algorithmicindent}
    $
    \mathfrak{P}_j = \mathfrak{T}_{\mathbf{r}_j}(\mathfrak{P}_j)
    $
    \Comment{If \eqref{eq:r_bound} with $M_r$ fails, \textbf{break}}

    \State Compute and store actual rank: \quad $\tilde{r}_j = \text{rank}(\mathfrak{P}_j)$
    \State Decrease counter: \quad $q = q - 1$
\EndFor
\end{algorithmic}
\end{algorithm}


\begin{algorithm}[!htb]
\caption{Training Phase of OCTANE}
\label{alg:OCTANE_train}
\begin{algorithmic}[1]

\Require Input and target data $(\hat{x}, x)$, total layers $N$, error bound constants $(M_s, M_r)$, final time $T$, batch iterations $m_1$, optimization solver iterations $m_2$
\Ensure Trained weights $\{K_j, b_j\}_{j=0}^{N_e-1}$, $\{\tilde{K}_j, \tilde{b}_j\}_{j=N_e}^{N-1}$; rank profiles $\mathbf{r}^{fe}$, $\mathbf{r}^{fd}$, $\mathbf{r}^{be}$, $\mathbf{r}^{bd}$; output $\mathfrak{g}_N$; training loss $\alpha_{\text{train}}$

\State Compute time step: $\tau \gets \frac{T}{N}$, and encoder depth: $N_e \gets \frac{N}{2}$
\State Initialize weights $\{K_j,b_j\}_{j=0}^{N_e-1}$ and $\{\tilde{K}_j,\tilde{b}_j\}_{j=N_e}^{N-1}$

\For{$i = 1$ to $m_1$}
    \State Randomly select mini-batch $(\hat{x}_i, x_i) \subset (\hat{x}, x)$
    \State Tensorize: $\hat{\mathfrak{x}} \gets \texttt{tensor}(\hat{x}_i)$, $\mathfrak{x} \gets \texttt{tensor}(x_i)$

    \State \textbf{Forward Encode:} Use \cref{alg:Bnet_fwd_encoder} with $\mathfrak{f}_0 = \hat{\mathfrak{x}}$, $li = 1$, and $ui = N_e$, to get $\{\mathfrak{f}_j\}_{j=0}^{N_e}$, and ranks $\mathbf{r}^{fe}$

    \State Define encoder ranks: $\mathbf{r}_e \gets \mathbf{r}^{fe}$, decoder ranks: $\mathbf{r}_d \gets \texttt{flip}(\mathbf{r}^{fe})$

    \State \textbf{Forward Decode:} Use \cref{alg:Bnet_fwd_decoder} with $\mathfrak{g}_{N_e} = \mathfrak{f}_{N_e}$, ranks $\mathbf{r}_d$, $li = N_e+1$, and $ui = N$, to get $\{\mathfrak{g}_j\}_{j=N_e}^{N}$, and $\mathbf{r}^{fd}$

    \State Compute terminal adjoint: $\tilde{\mathfrak{P}}_N \gets -\nabla_g \mathcal{J}(\mathfrak{x}, \mathfrak{g}_N)$

    \State \textbf{Backward Decode:} Use \cref{alg:Bnet_bckwd_decoder} with $\tilde{\mathfrak{P}}_N$, ranks $\mathbf{r}_d$, $li = N_e$, and $ui = N-1$, to get $\{\tilde{\mathfrak{P}}_j\}_{j=N_e}^{N}$, and $\mathbf{r}^{bd}$

    \State \textbf{Backward Encode:} Use \cref{alg:Bnet_bckwd_encoder} with $\mathfrak{P}_{N_e} = \tilde{\mathfrak{P}}_{N_e}$, ranks $\mathbf{r}_e$, $li = 0$, and $ui = N_e-1$, to get $\{\mathfrak{P}_j\}_{j=0}^{N_e}$, and $\mathbf{r}^{be}$

    \State \textbf{Gradient Computation:} Compute gradients using \cref{DesEq_disc}
    \[
    \begin{aligned}
    \nabla_{K_j}\mathcal{L} &= - \mathfrak{f}_j \left( \mathfrak{P}_{j+1} \odot \sigma'(K_j \mathfrak{f}_j + b_j) \right)^\intercal + \nabla_K \mathcal{R}_j, &\qquad \forall\;\; j = 0,...,N_e-1 \\
    \nabla_{b_j}\mathcal{L} &= - \langle \sigma'(K_j \mathfrak{f}_j + b_j), \mathfrak{P}_{j+1} \rangle_F + \nabla_b \mathcal{R}_j, &\forall \;\;j = 0,...,N_e-1 \\
    \nabla_{\tilde{K}_j}\mathcal{L} &= - \mathfrak{g}_j \left( \tilde{\mathfrak{P}}_{j+1} \odot \tilde{\sigma}'(\tilde{K}_j \mathfrak{g}_j + \tilde{b}_j) \right)^\intercal + \nabla_{\tilde{K}} \mathcal{R}_j, &\forall\;\; j = N_e,...,N-1 \\
    \nabla_{\tilde{b}_j}\mathcal{L} &= - \langle \tilde{\sigma}'(\tilde{K}_j \mathfrak{g}_j + \tilde{b}_j), \tilde{\mathfrak{P}}_{j+1} \rangle_F + \nabla_{\tilde{b}} \mathcal{R}_j, &\forall \;\;j = N_e,...,N-1
    \end{aligned}
    \]

    \State \textbf{Update Weights:} Use gradient-based optimizer with tolerance $\eta$ and max iterations $m_2$ to update all parameters

    \State \textbf{Compute Training Loss:} $\alpha_{\text{train}} \gets \mathcal{J}(\cdot, \mathfrak{g}_N, \mathfrak{x})$
\EndFor

\State \textbf{Post-processing:} Plot reconstruction $\mathfrak{g}_N$ and rank profiles $\mathbf{r}^{fe}, \mathbf{r}^{fd}$
\end{algorithmic}
\end{algorithm}

\begin{algorithm}[!htb]
\caption{Testing Phase of OCTANE}
\label{alg:OCTANE_test}
\begin{algorithmic}[1]

\Require Test data $(\hat{x}_{\text{test}}, x_{\text{test}})$, trained weights $\{K_j, b_j\}_{j=0}^{N_e-1}$ and $\{\tilde{K}_j, \tilde{b}_j\}_{j=N_e}^{N-1}$, encoder ranks $\mathbf{r}^{fe}$, layer count $N$, bounds $(M_s, M_r)$, final time $T$
\Ensure Test encoder and decoder ranks $\mathbf{r}^{fe}_{\text{test}}, \mathbf{r}^{fd}_{\text{test}}$; reconstruction $\mathfrak{g}_N$; test loss $\alpha_{\text{test}}$

\State Compute step size and depth: $\tau \gets \frac{T}{N}$, $N_e \gets \frac{N}{2}$
\State Define prescribed ranks: $\mathbf{r}_e \gets \mathbf{r}^{fe}$, $\mathbf{r}_d \gets \texttt{flip}(\mathbf{r}^{fe})$
\State Tensorize test data: $\hat{\mathfrak{x}} \gets \texttt{tensor}(\hat{x}_{\text{test}})$, $\mathfrak{x} \gets \texttt{tensor}(x_{\text{test}})$

\State \textbf{Forward Encode:} Use \cref{alg:Bnet_fwd_encoder} with initial condition $\mathfrak{f}_0 = \hat{\mathfrak{x}}$ and ranks $\mathbf{r}_e$, to compute $\{\mathfrak{f}_j\}_{j=0}^{N_e}$ and test encoder ranks $\mathbf{r}^{fe}_{\text{test}}$

\State \textbf{Forward Decode:} Use \cref{alg:Bnet_fwd_decoder} with initial condition $\mathfrak{g}_{N_e} = \mathfrak{f}_{N_e}$ and ranks $\mathbf{r}_d$, to compute $\{\mathfrak{g}_j\}_{j=N_e}^{N}$ and test decoder ranks $\mathbf{r}^{fd}_{\text{test}}$

\State Compute reconstruction error: $\alpha_{\text{test}} \gets \mathcal{J}(\cdot, \mathfrak{g}_N, \mathfrak{x})$

\State \textbf{Post-processing:} Plot reconstruction $\mathfrak{g}_N$ and ranks $\mathbf{r}^{fe}_{\text{test}}, \mathbf{r}^{fd}_{\text{test}}$

\end{algorithmic}
\end{algorithm}

\subsection{Computational Cost and Memory}\label{comp_cost}
Recall that each optimization step in training (\textit{line 13}, \cref{alg:OCTANE_train}) involves solving two state and two adjoint differential equations via the forward/backward encoder and decoder routines (\cref{alg:Bnet_fwd_encoder,alg:Bnet_fwd_decoder,alg:Bnet_bckwd_encoder,alg:Bnet_bckwd_decoder}). Thus, the computational cost of \textsc{OCTANE} is comparable to a standard RNN \cite{Ruthotto_2020_DL_OC}, but with dual ODE solves in each pass. However, our low-rank tensor formulation significantly reduces memory usage and storage overhead by retaining only compressed solution trajectories. This not only lowers memory requirements but also reduces overall runtime.

\section{Numerical Experiments\label{s:Numerics}}
Until now, we have presented a general formulation suitable for a broad class of problems. We now outline specific choices made for our numerical experiments, designed to illustrate the effectiveness of the proposed framework.

\subsection{Preliminaries for Numerical Experiments}
\subsubsection{Data Fidelity, Reconstruction Error and Regularization.}
We adopt the generalized Euclidean distance in $L^2$ as the 
data fidelity term, defined by 
\begin{equation}\label{Jdefl2}
J(x, g(T)) := \frac{1}{2n} \;\| \hat{\sigma}(g(T)) - x \|^2_{L^2},
\end{equation}
where $\hat{\sigma}(\cdot)$ is a nonlinear, pointwise activation function applied to the output layer $g(T)$.

To promote regularity in the network parameters, we introduce the following regularizer:
\begin{equation}\label{Regdef}
\begin{aligned}
\mathcal{R}(\Theta) = \;\;& \frac{\lambda_1}{2N_e} \sum_{j=0}^{N_e-1} \| K(t_j) \|_F^2 + \frac{\lambda_2}{2N_d} \sum_{j=N_e}^{N-1} \| \tilde{K}(t_j) \|_F^2 \\
&+ \frac{\lambda_3}{2N_e} \sum_{j=0}^{N_e-1} |b(t_j)|^2 + \frac{\lambda_4}{2N_d} \sum_{j=N_e}^{N-1} |\tilde{b}(t_j)|^2,
\end{aligned}
\end{equation}
where $\{\lambda_i\}_{i=1}^{4}$ are regularization parameters selected heuristically.

The overall objective function, defining the average reconstruction error (i.e., the regularized mean squared error) over all $n$ samples, is given by:
\begin{equation}\label{recon_err}
\alpha := \mathcal{J}(\Theta, g(T), x) 
= J(x, g(T)) + \mathcal{R}(\Theta).
\end{equation}
Here, $\| \cdot \|_F$ denotes the Frobenius norm, and $N_d = N - N_e$. The $L^2$-norm in \cref{Jdefl2} is evaluated using piecewise linear interpolation.

Furthermore, recall, \cref{AdjEq} requires the gradient $\nabla_g J(x,g(T)$, which would be given by,
\begin{equation}\label{Jgradg_mse}
\nabla_{g} J(x,g^*(T)) = \frac{1}{n}(\hat{\sigma}^{\prime}(g^*(T))^{\intercal}\left(\hat{\sigma}(g^*(T)) - x\right),
\end{equation}
with the tensorized version for \cref{AdjEq_disc} written as,
\begin{equation}\label{Jgradg}
\nabla_{\mathfrak{g}} J(\mathfrak{x},\mathfrak{g}^*(t_N)) =  \frac{1}{n}(\hat{\sigma}^{\prime}(\mathfrak{g}^*(t_N))^{\intercal}\left(\hat{\sigma}(\mathfrak{g}^*(t_N)) - \mathfrak{x}\right).
\end{equation}

\subsubsection{Activation Functions}
In our experiments, we use the hyperbolic tangent as the activation for both encoder and decoder:
\[
\sigma(x) = \tilde{\sigma}(x) = \tanh(x), \quad \sigma'(x) = \tilde{\sigma}'(x) = 1 - \tanh^2(x).
\]
For image reconstruction tasks, where pixel intensities are non-negative, we employ the smoothed ReLU from \cite[eq.~(3.2)]{antil2021_DNNBayes} as the output activation $\hat{\sigma}$, serving as a soft constraint to promote non-negativity in the output.

\subsubsection{Gradient Tests}
\begin{figure}[!htb]
\begin{center}
\includegraphics[width=0.3\textwidth, height=0.25\textwidth]{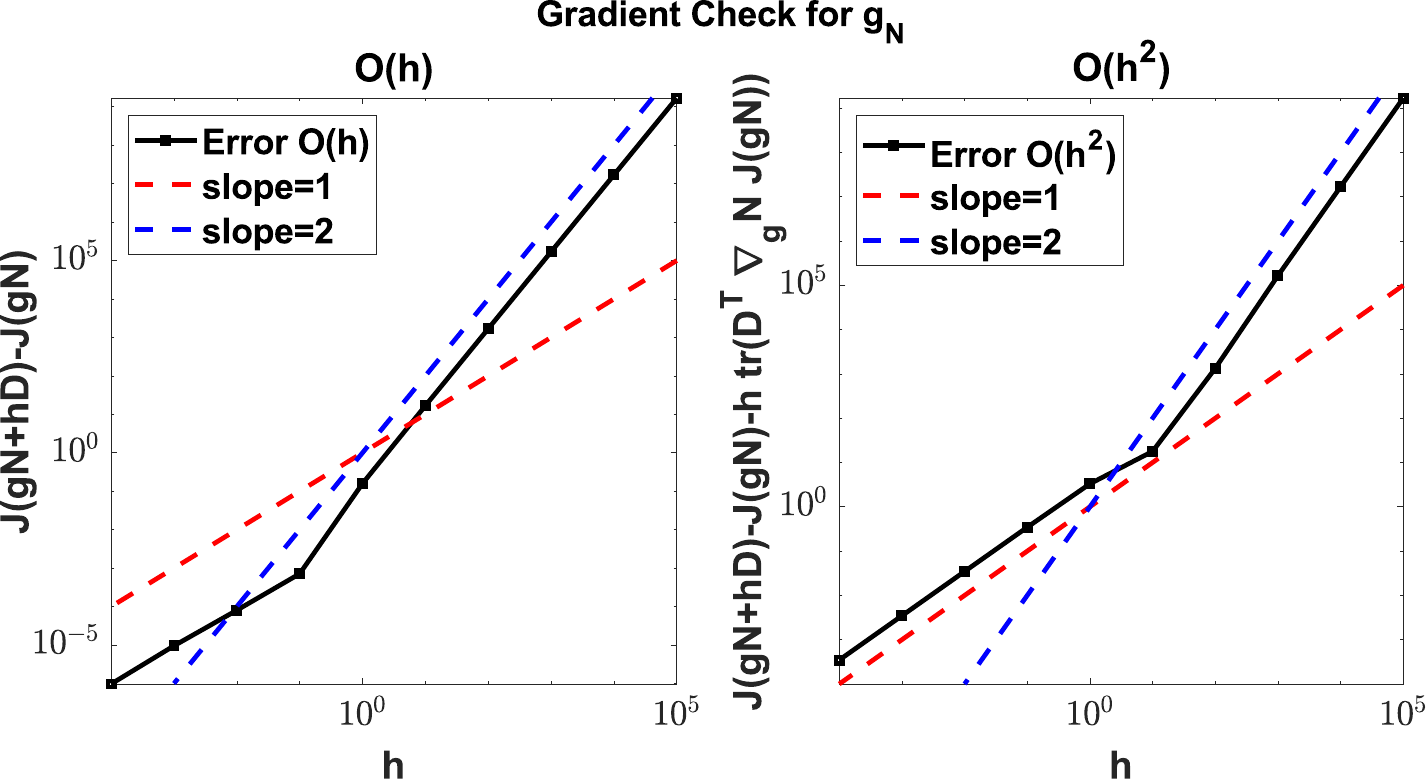} \quad
\includegraphics[width=0.3\textwidth, height=0.25\textwidth]{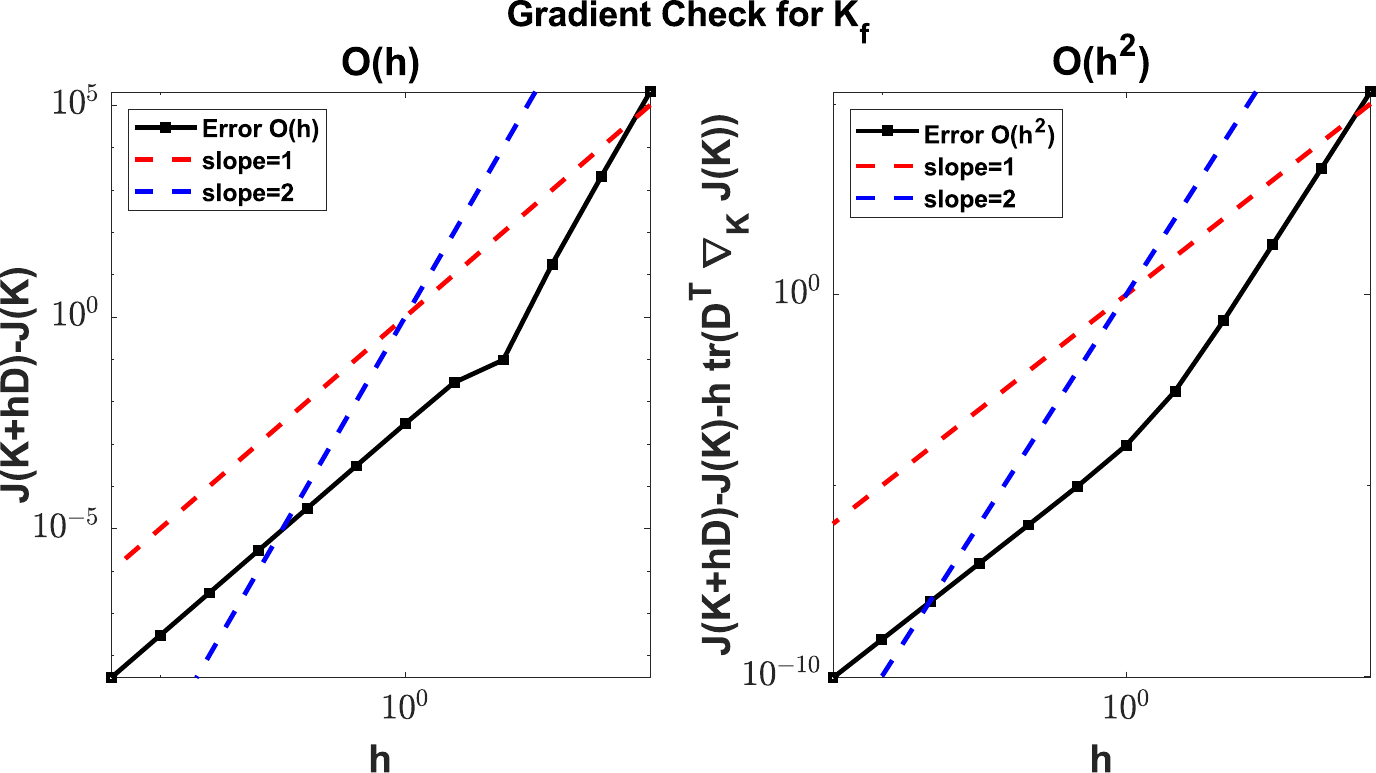} \quad 
\includegraphics[width=0.3\textwidth, height=0.25\textwidth]{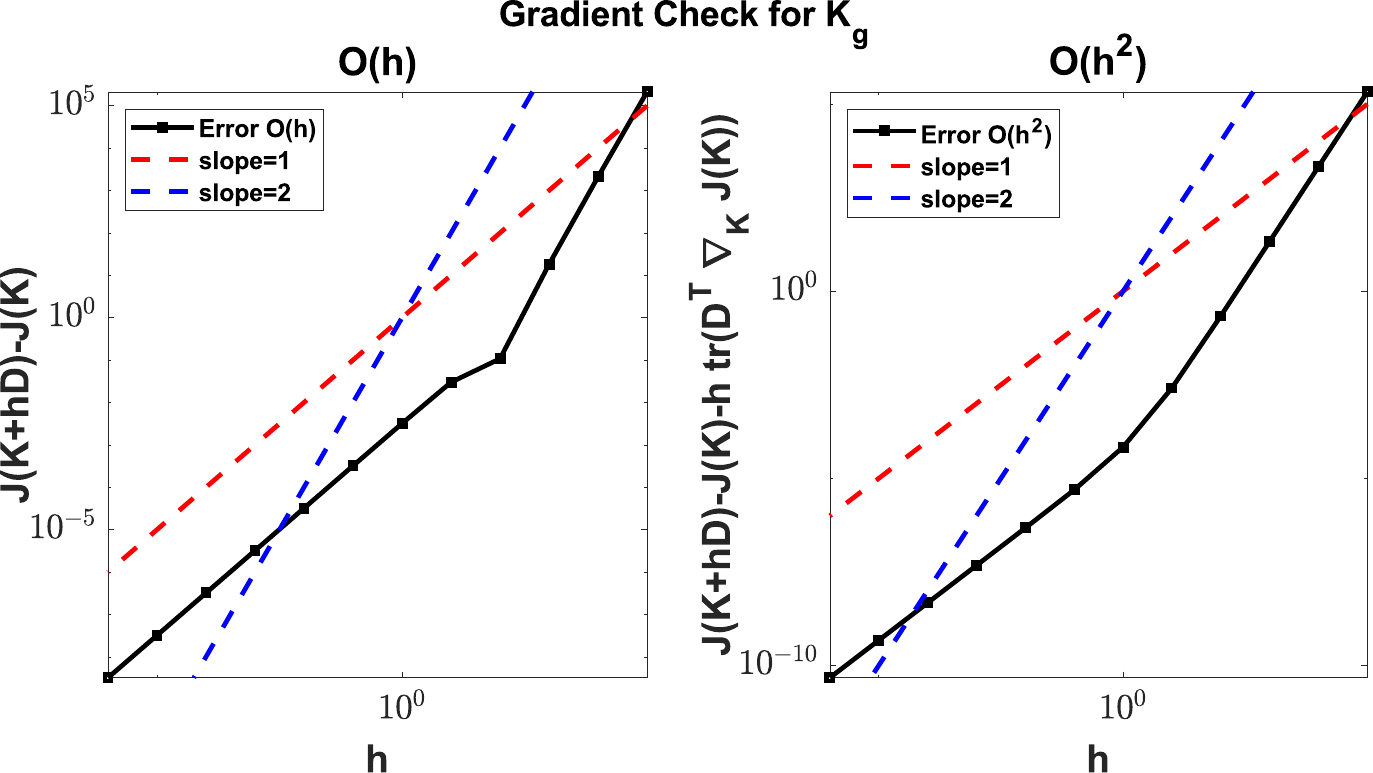}\\
\includegraphics[width=0.3\textwidth, height=0.25\textwidth]{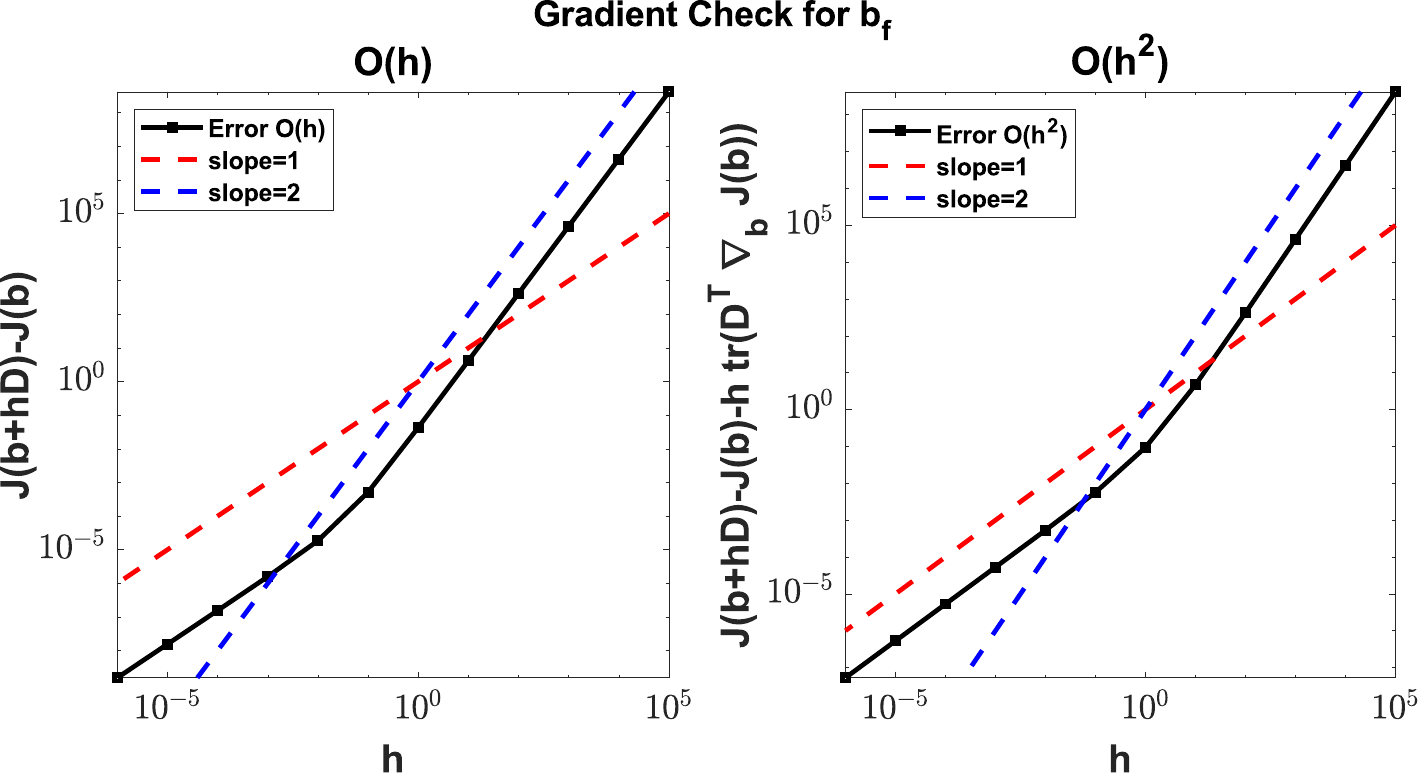} \qquad \qquad \qquad
\includegraphics[width=0.3\textwidth, height=0.25\textwidth]{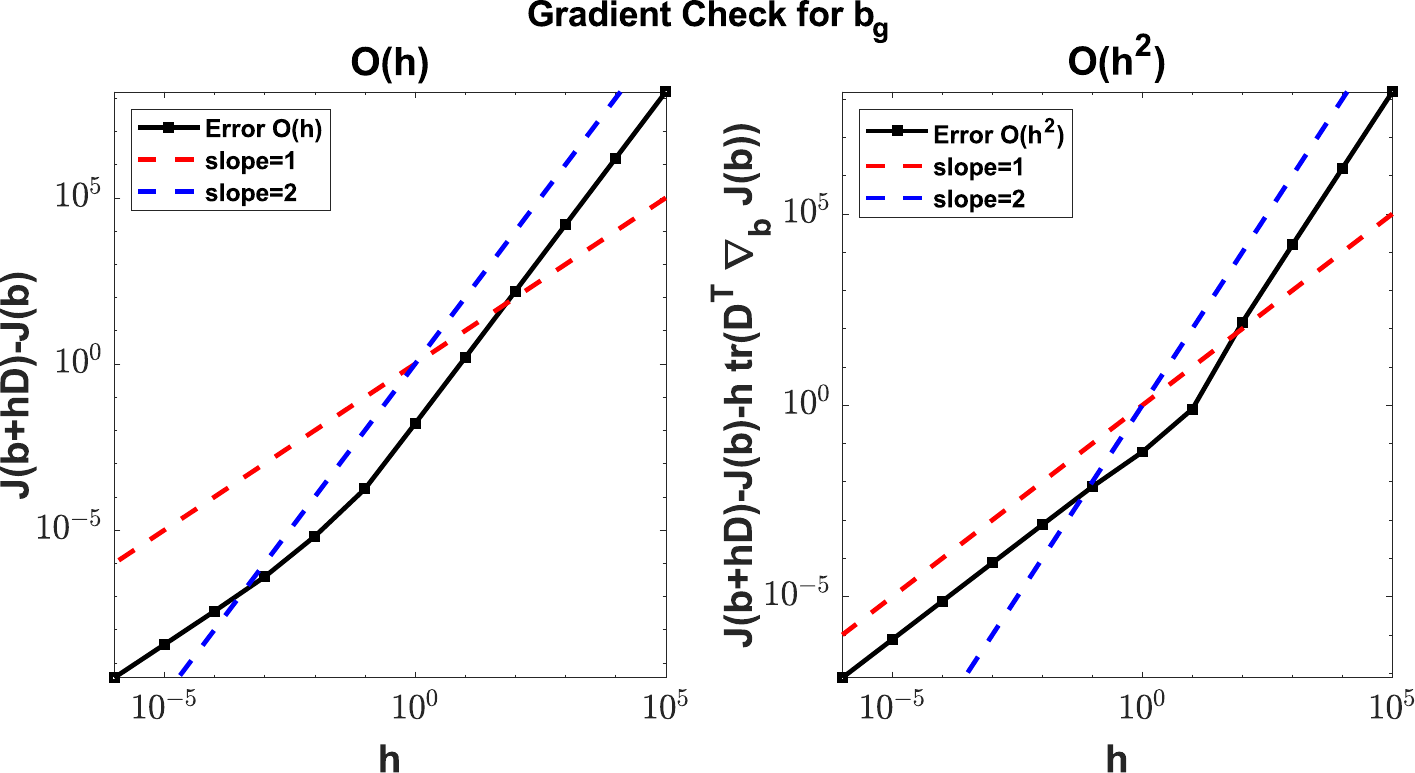}
\caption{\label{f:deriv_test}Comparison of discretized gradients with finite difference approximations. The black line shows the error, and blue and red are reference lines with slopes $1$ and $2$, respectively. The expected rate of convergence is obtained.}
\end{center}
\end{figure}

To validate the gradients in \cref{DesEq_disc} and $\nabla_g J(\mathfrak{x},\mathfrak{g}_N)$ from \cref{Jgradg}, we conduct a gradient test using synthetic image data (random vertical streaks), comparing discretized gradients against finite difference approximations from \cref{DesEq} and \cref{Jgradg_mse}. As shown in \cref{f:deriv_test}, the results align and exhibit the expected convergence order for all design variables and the objective gradient.

\subsubsection{Tensor Format and Toolbox}
We do not restrict the choice of tensor representation or computational toolbox; any format and library may be used. In our experiments, we adopt the Tensor-Train (TT) format \cite{Oseledets_2011} and employ the \texttt{TT-toolbox}\footnote{\texttt{https://github.com/oseledets/TT-Toolbox}} in MATLAB for tensor operations.

\subsubsection{Tensor Truncation with \texttt{TT-toolbox}} \label{Error_merge}

As discussed in \cref{rank_select}, any tensor representation and its associated truncation (or rounding) operation introduces an inherent approximation error. For the Tensor-Train (TT) format used in our computations, this truncation error is governed by the tolerance parameter $\varepsilon$. Specifically, \cite[Corollary 2.4]{Oseledets_2011} provides the following bound for a tensor $\mathfrak{a}$:
\begin{equation}\label{Osledet_bound}
\|\mathfrak{a} - \mathfrak{T}(\mathfrak{a})\|_2 \leq \varepsilon \|\mathfrak{a}\|_2.
\end{equation}

In our setting, the inequalities in \cref{eq:s_bound} and \cref{eq:r_bound} (and their counterparts \cref{eq:s_bound_rev}, \cref{eq:r_bound_rev}) bound the truncation error in the explicit Euler step after applying ${\mathfrak{T}}_{s}$ and ${\mathfrak{T}}_{r}$. Matching these bounds with \cref{Osledet_bound} yields the truncation tolerances:
\begin{equation}\label{merged_bounds}
\begin{aligned}
\varepsilon = \frac{M_s \tau}{\|\mathcal{N}(Y)\|_2}, \quad &\text{(resp. } \varepsilon = \frac{M_s \tau}{\|\mathcal{N}(Z)\|_2}\text{)} \\
\varepsilon = \frac{M_r \tau^2}{\|Y + \tau\,\mathfrak{T}_s(\mathcal{N}(Y))\|_2}, \quad &\text{(resp. } \varepsilon = \frac{M_r \tau^2}{\|Z + \tau\,\mathfrak{T}_s(\mathcal{N}(Z))\|_2}\text{)}.
\end{aligned}
\end{equation}

Thus, given $\tau$, $M_s$, and $M_r$, \cref{merged_bounds} defines the truncation bounds passed to the \texttt{TT-toolbox} \texttt{round()} function for rank reduction.


\subsubsection{Optimization and Xavier Initialization} 

We employ the BFGS algorithm with Armijo line search \cite{Kelly_bfgs}, terminating when the gradient norm falls below $10^{-5}$ or after $m_2$ iterations—typically the latter. All design variables are initialized using Xavier initialization \cite{Bengio2010}.

\subsubsection{Data and Data Batches}

We use the MNIST dataset \cite{MNIST} for our experiments, selecting a single digit and dividing its images into training, validation, and testing sets. Data is scaled to $[0,1]$ and stacked so that each frontal slice of the resulting tensor (used in \texttt{TT-toolbox}) corresponds to one image.

During training (\cref{alg:OCTANE_train}), $50\%$ of the training data is randomly sampled as a mini-batch in each of the $m_1$ iterations. For testing (\cref{alg:OCTANE_test}), the test data is divided into batches (of size $20$), and reconstruction errors are averaged across batches to compute the final testing error.

\subsubsection{Experimental Configuration} In \cref{t:exp_config} we provide some common configurations for all the experiments discussed in \cref{Sec:Denoising} and \cref{Sec:Deblurring}.

\begin{normalsize}
\begin{table}[h!]
\caption{\label{t:exp_config}Experimental configuration for OCTANE autoencoder used for image denoising and image deblurring tasks.}
\begin{tabular}{|c|c|c|c|c|c|c|c|c|c|c|c|c|}
\hline
\textbf{Exp. Type}  & \textbf{$n_{train}$} & \textbf{$n_{valid}$} & \textbf{$n_{test}$} & \textbf{$T$} & \textbf{$M_s$} & \textbf{$M_r$} & \textbf{$\lambda_1$} & \textbf{$\lambda_2$} & \textbf{$\lambda_3$} & \textbf{$\lambda_4$} & \textbf{$m_1$} & \textbf{$m_2$} \\ \hline
\textbf{Denoising}  & $20$                 & $20$                 & $1000$              & $10$         & $\tau^{-1}$    & $\tau^{-2}$    & $1e-5$               & $1e-5$               & $1$                    & $1$                    & $3$              & $30$              \\ \hline
\textbf{Deblurring} & $20$                 & $20$                 & $1000$              & $10$          & $\tau^{-1}$    & $\tau^{-2}$    & $0$                    & $0$                    & $1e-1$               & $1e-1$               & $3$             & $30$             \\ \hline
\end{tabular}
\end{table}
\end{normalsize}
\subsubsection{Computational Platform}
All the computations have been carried out in MATLAB R$2024$a, on a laptop with an Intel Core i$7-12700$H processor. 
\subsection{Image Denoising}\label{Sec:Denoising} 

We apply the OCTANE autoencoder (\cref{alg:OCTANE_train,alg:OCTANE_test}) to denoise MNIST images. Clean images \( x \) (e.g., digit 2) are corrupted with $5\%$ Gaussian noise to produce synthetic noisy inputs \( \hat{x} \), serving as initial conditions in \cref{AE_OC}. The autoencoder learns a low-rank representation of \( \hat{x} \) and reconstructs denoised outputs \( g(f(\hat{x})) \).

\begin{figure}[h!]
\begin{center}
\hspace{-0.4cm}\fbox{\text{\textbf{Initial Cond.}}} \hspace{1.2cm}\fbox{$\textbf{N=6}$} \hspace{1.9cm}\fbox{$\textbf{N=12}$} \hspace{1.7cm} \fbox{$\textbf{N=20}$} \hspace{1.3cm}\fbox{\text{\textbf{Ref. Sol.}}}\\~\\

\includegraphics[width=0.14\textwidth, height=0.13\textwidth]{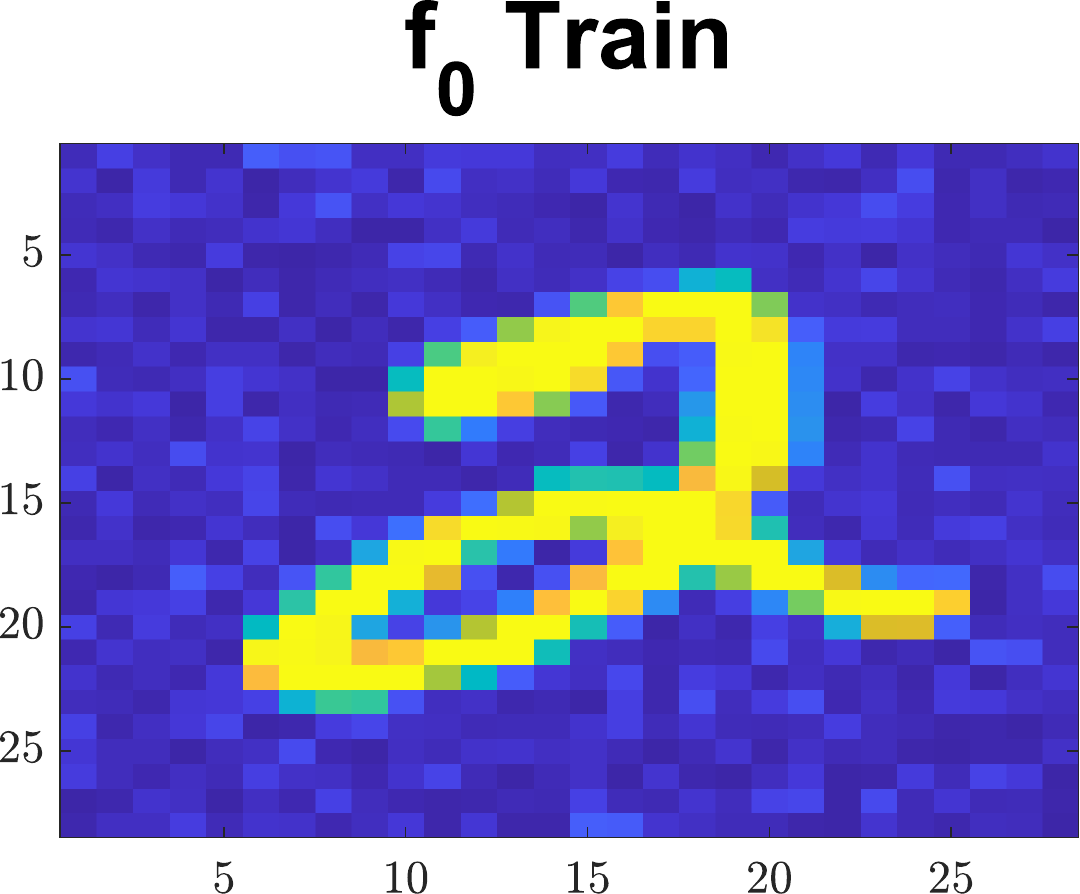} 
\hspace{0.8cm}\includegraphics[width=0.14\textwidth, height=0.14\textwidth]{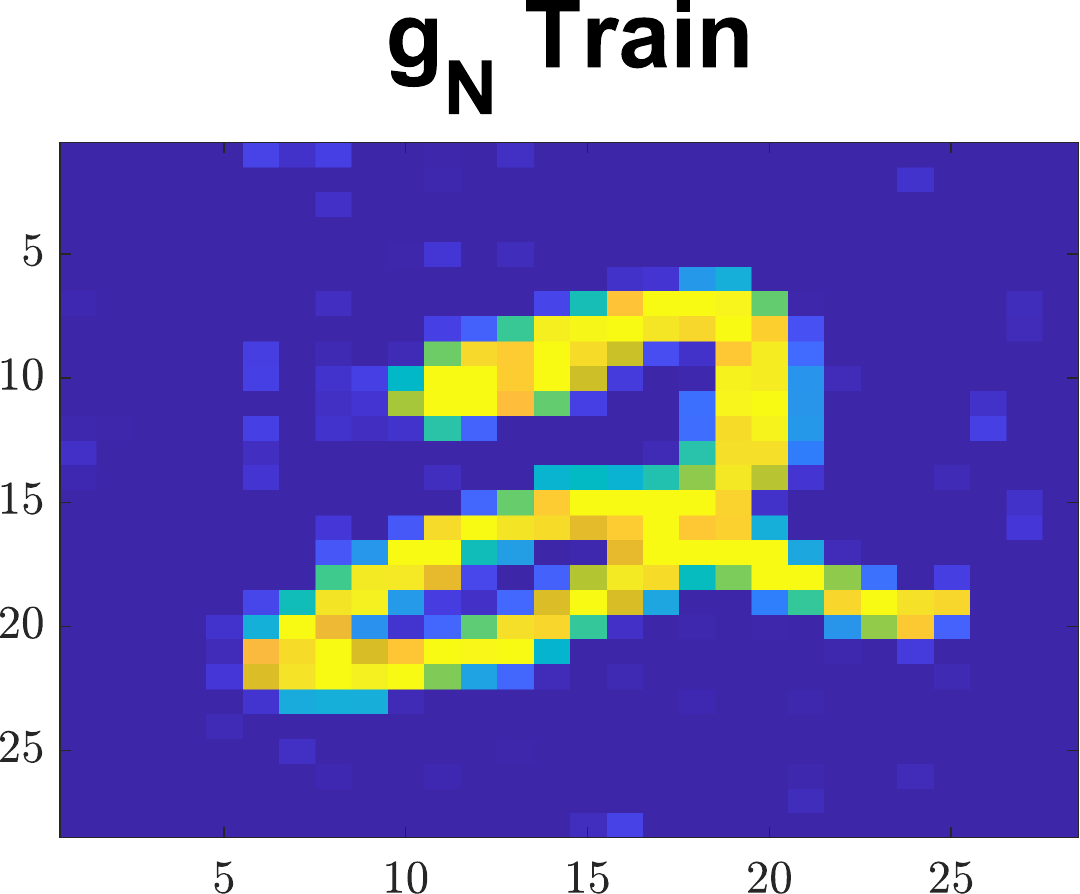}
\hspace{0.8cm}\includegraphics[width=0.14\textwidth, height=0.14\textwidth]{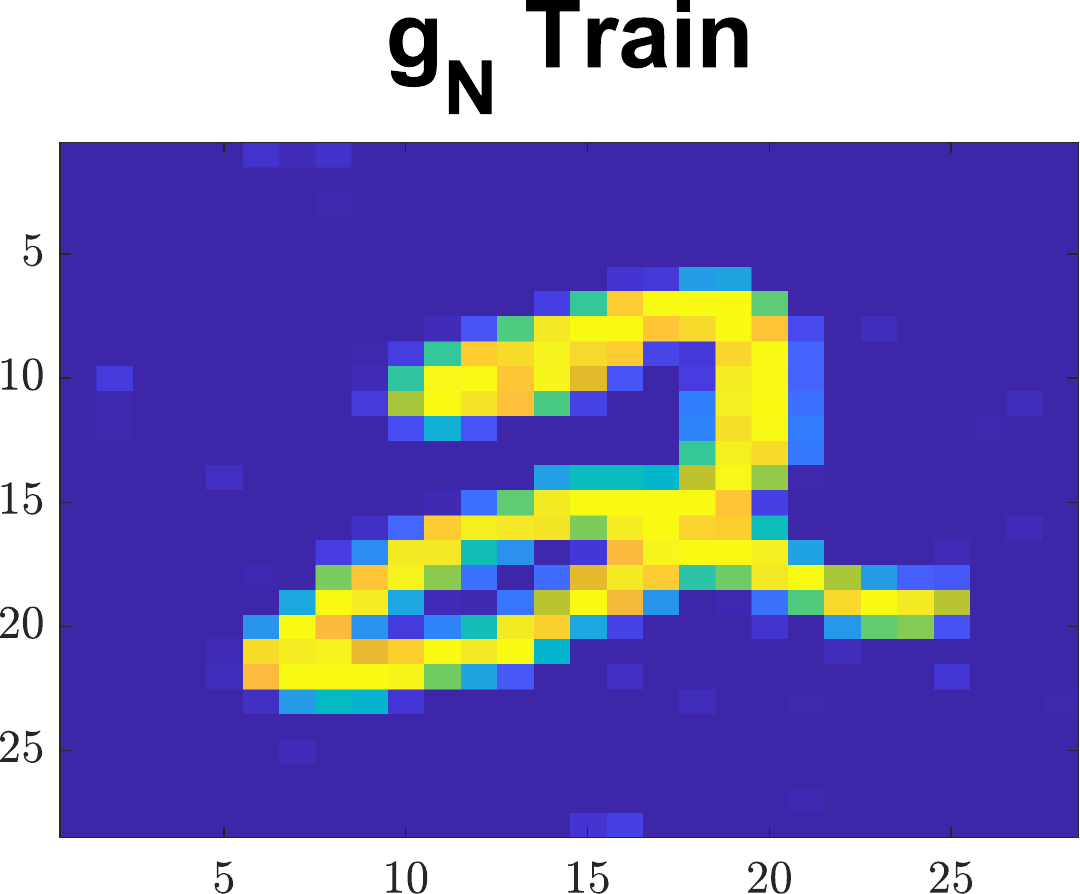}
\hspace{0.8cm}\includegraphics[width=0.14\textwidth, height=0.14\textwidth]{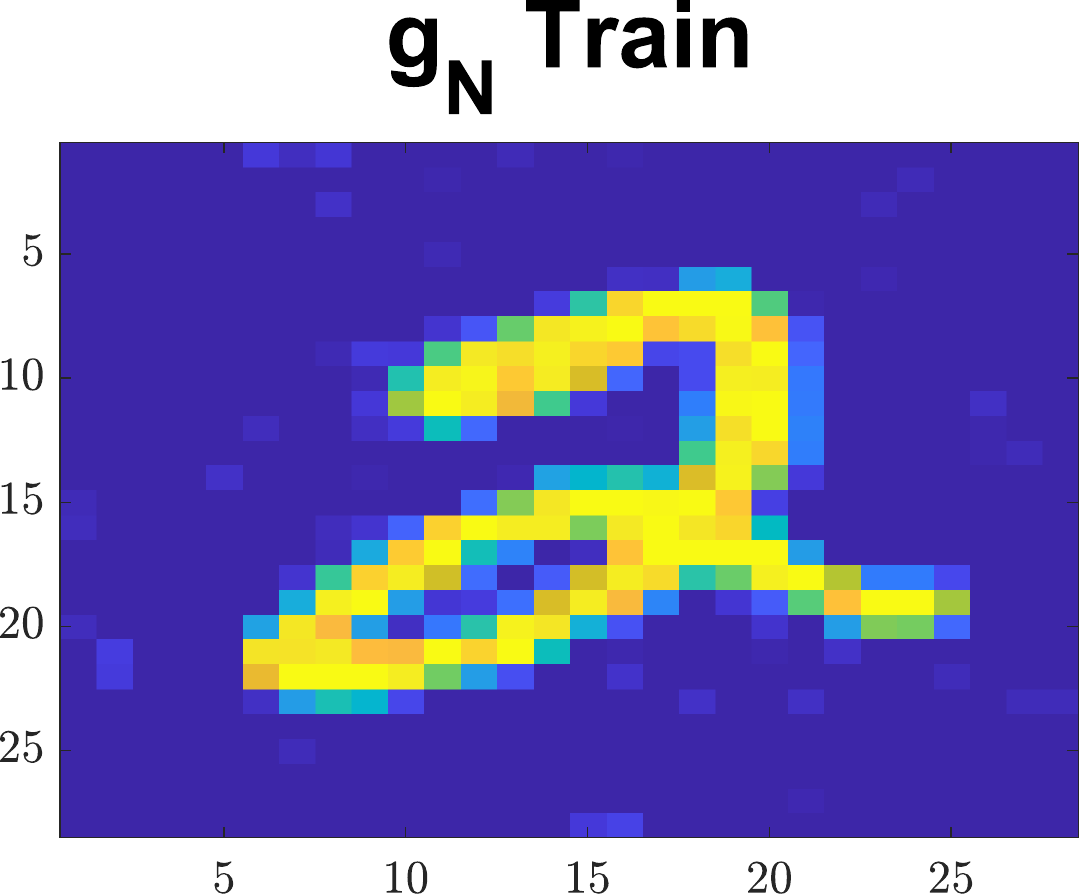}
\hspace{0.8cm}\includegraphics[width=0.15\textwidth, height=0.13\textwidth]{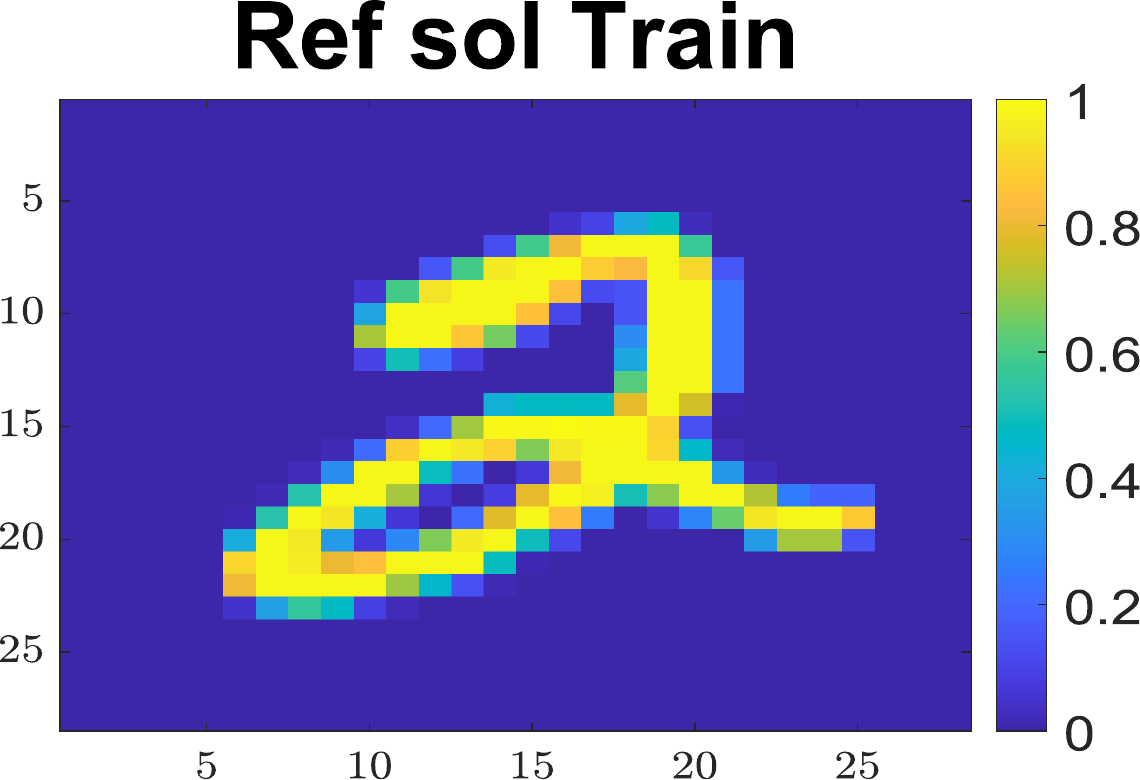}\\
\includegraphics[width=0.14\textwidth, height=0.13\textwidth]{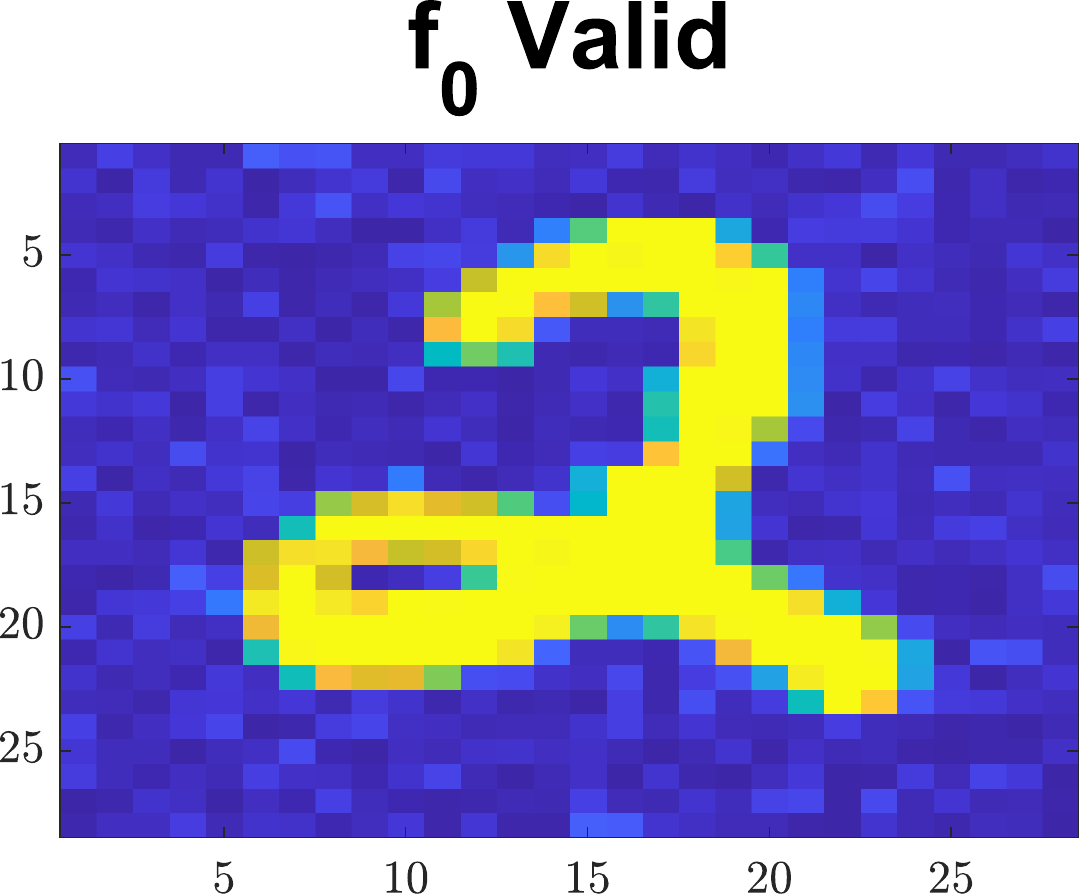} 
\hspace{0.8cm}\includegraphics[width=0.14\textwidth, height=0.14\textwidth]{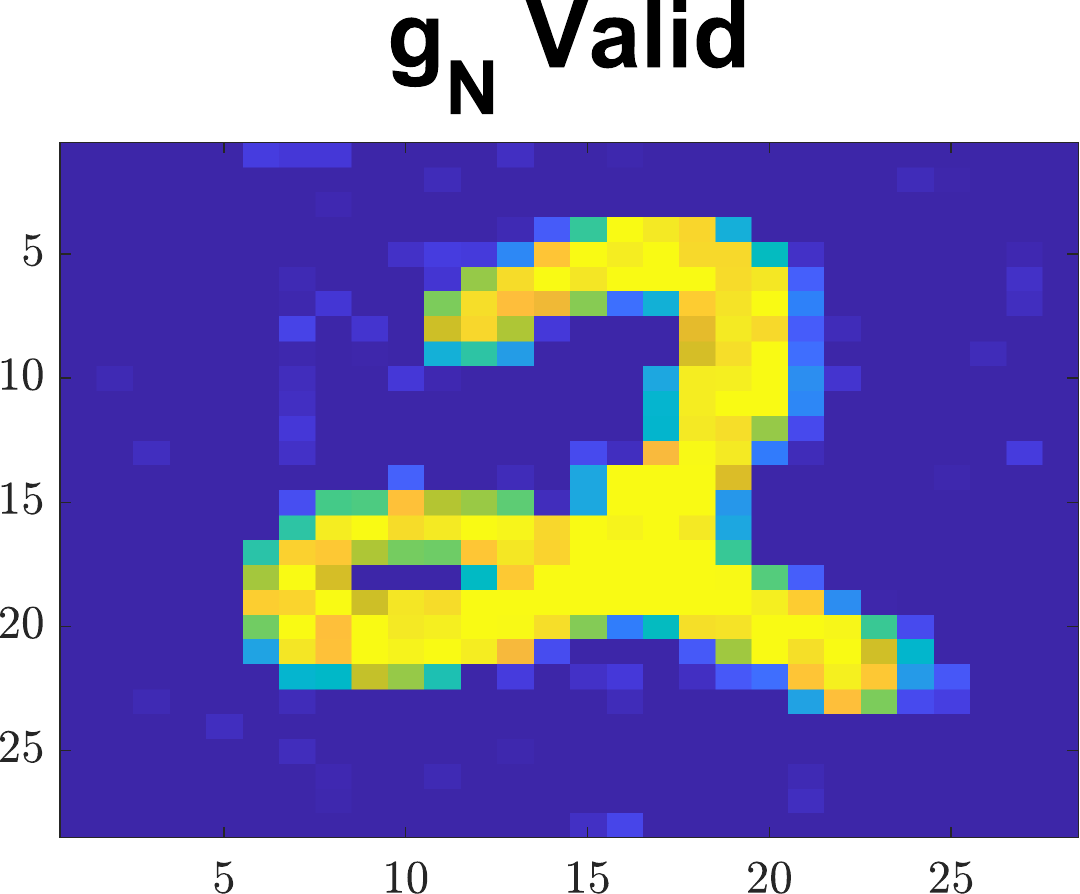} 
\hspace{0.8cm}\includegraphics[width=0.14\textwidth, height=0.14\textwidth]{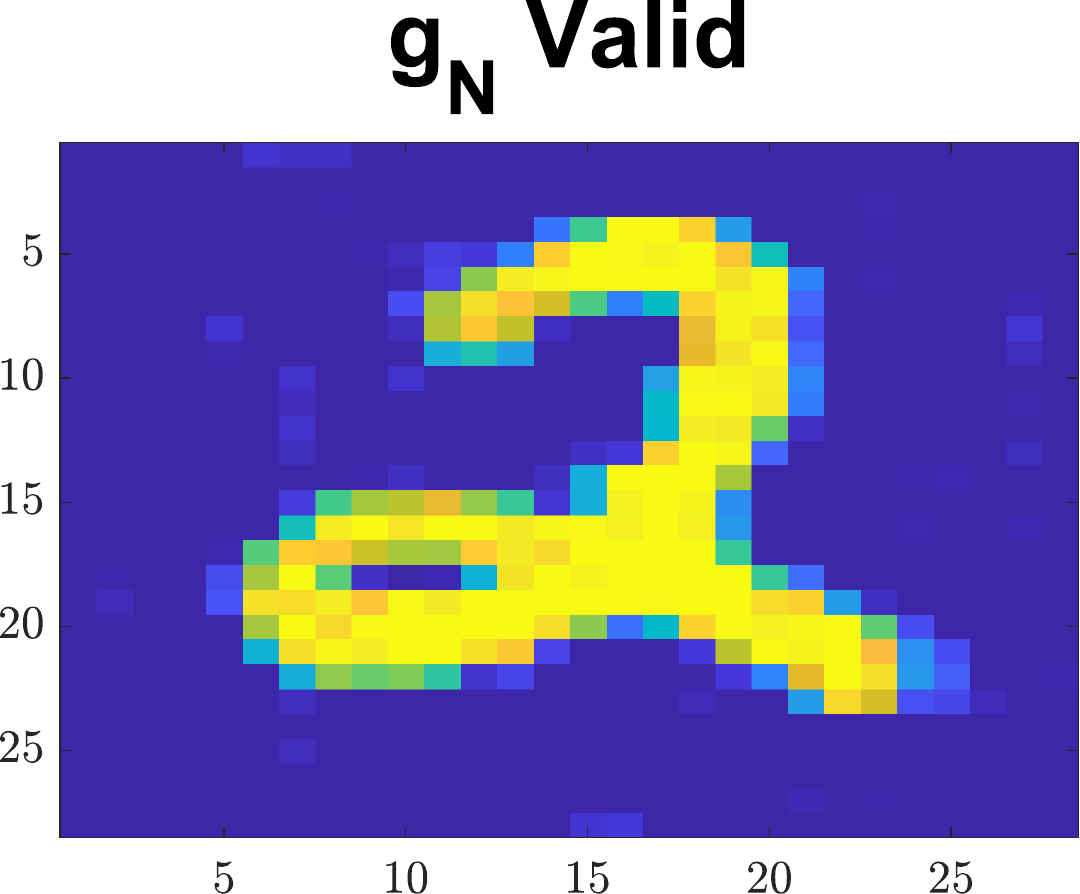}
\hspace{0.8cm}\includegraphics[width=0.14\textwidth, height=0.14\textwidth]{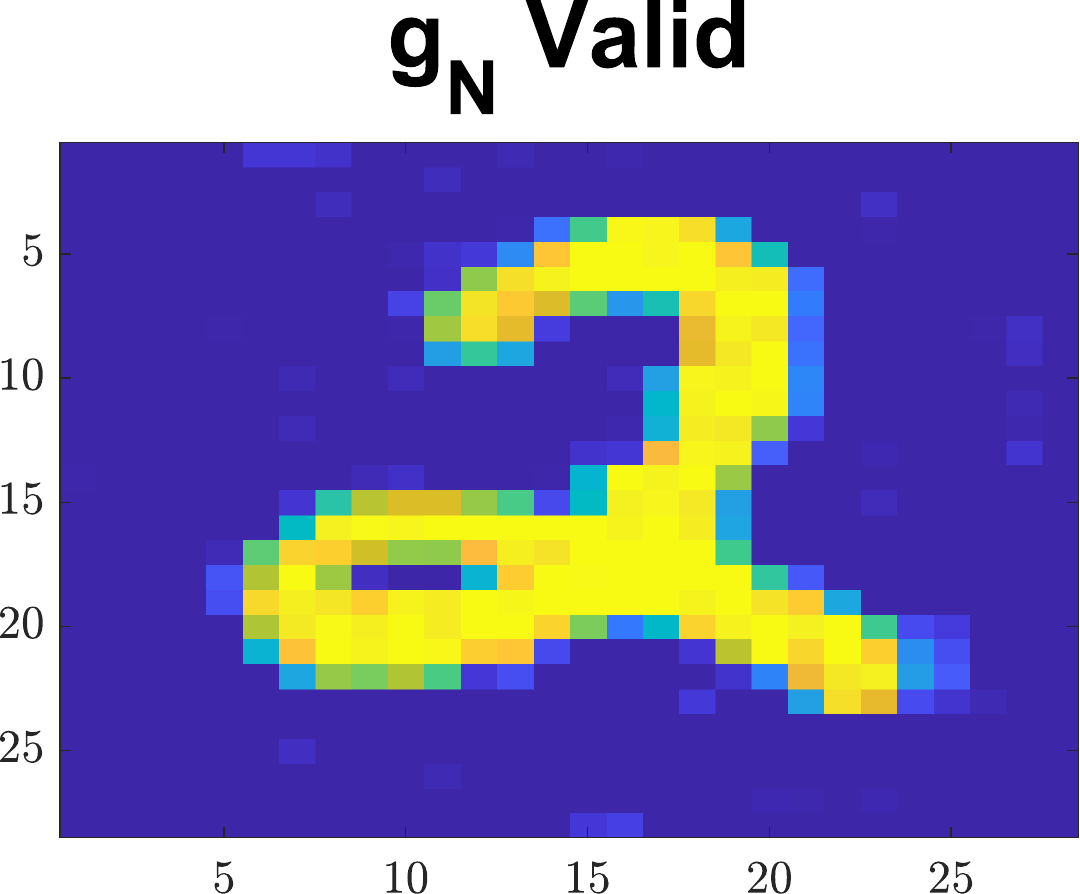}
\hspace{0.8cm}\includegraphics[width=0.15\textwidth, height=0.13\textwidth]{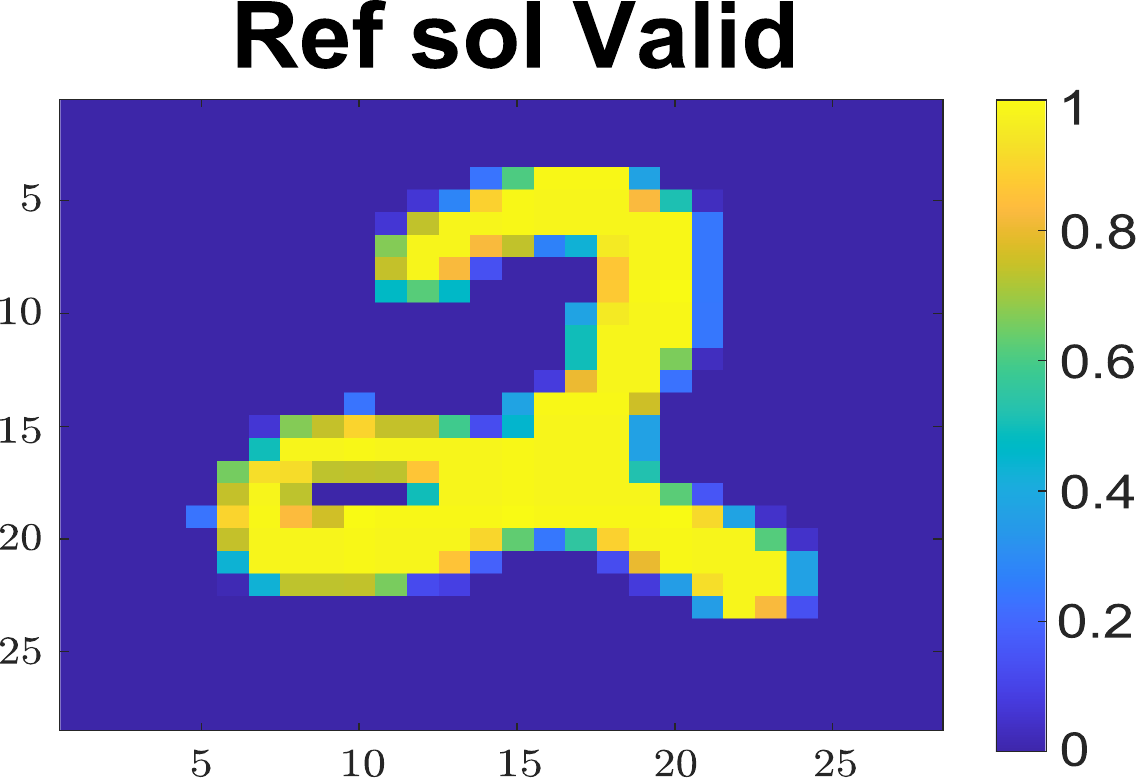}\\
\includegraphics[width=0.14\textwidth, height=0.13\textwidth]{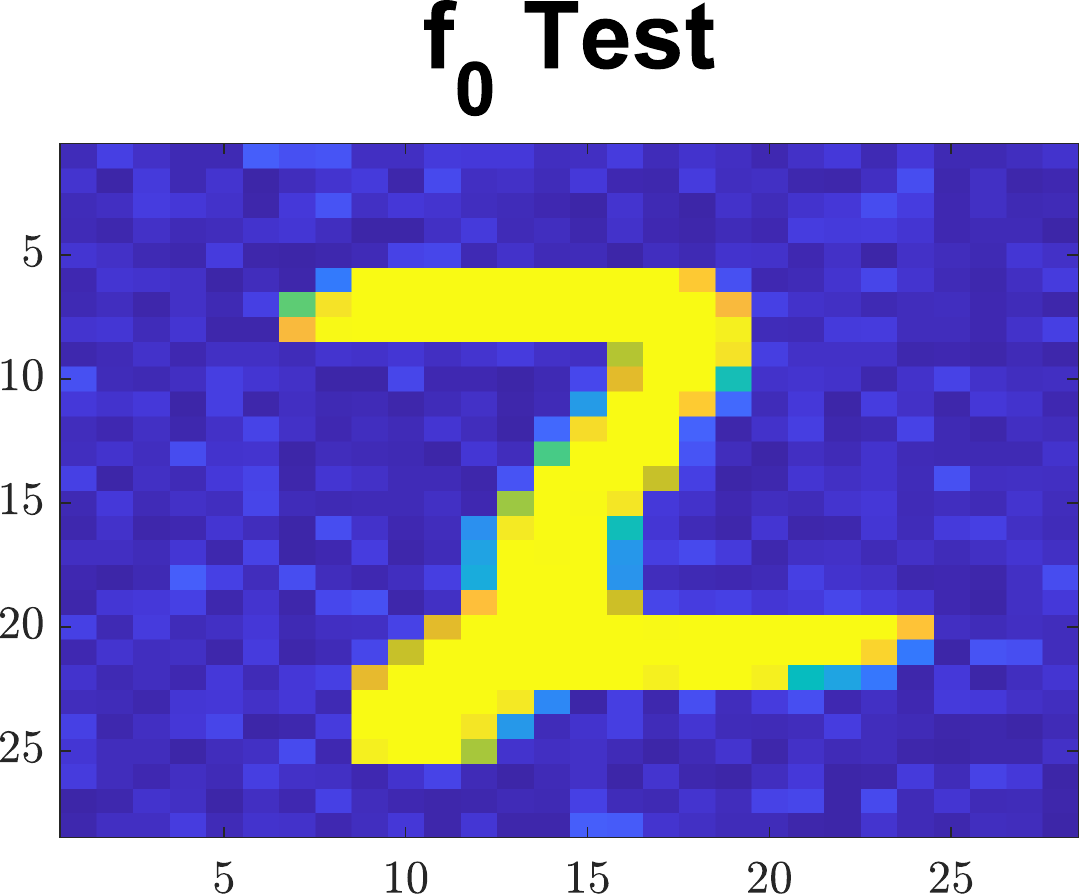}
\hspace{0.8cm}\includegraphics[width=0.14\textwidth, height=0.14\textwidth]{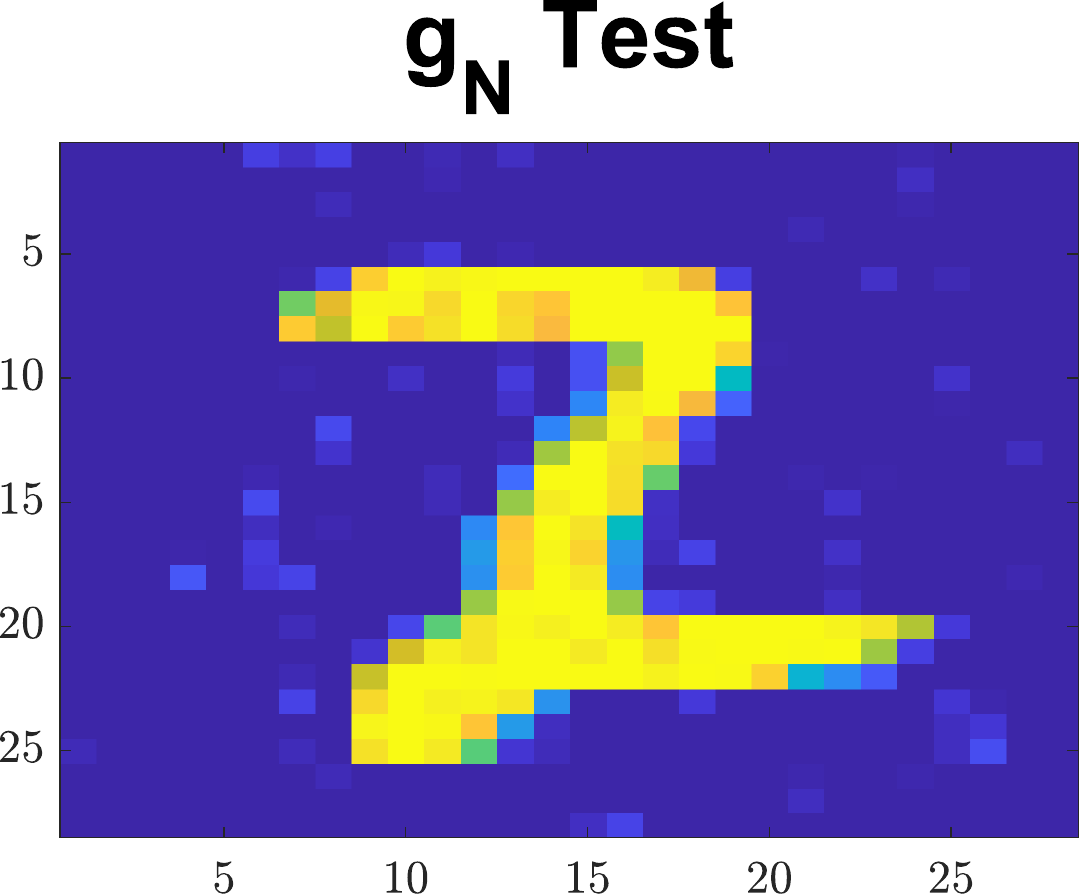}
\hspace{0.8cm}\includegraphics[width=0.14\textwidth, height=0.14\textwidth]{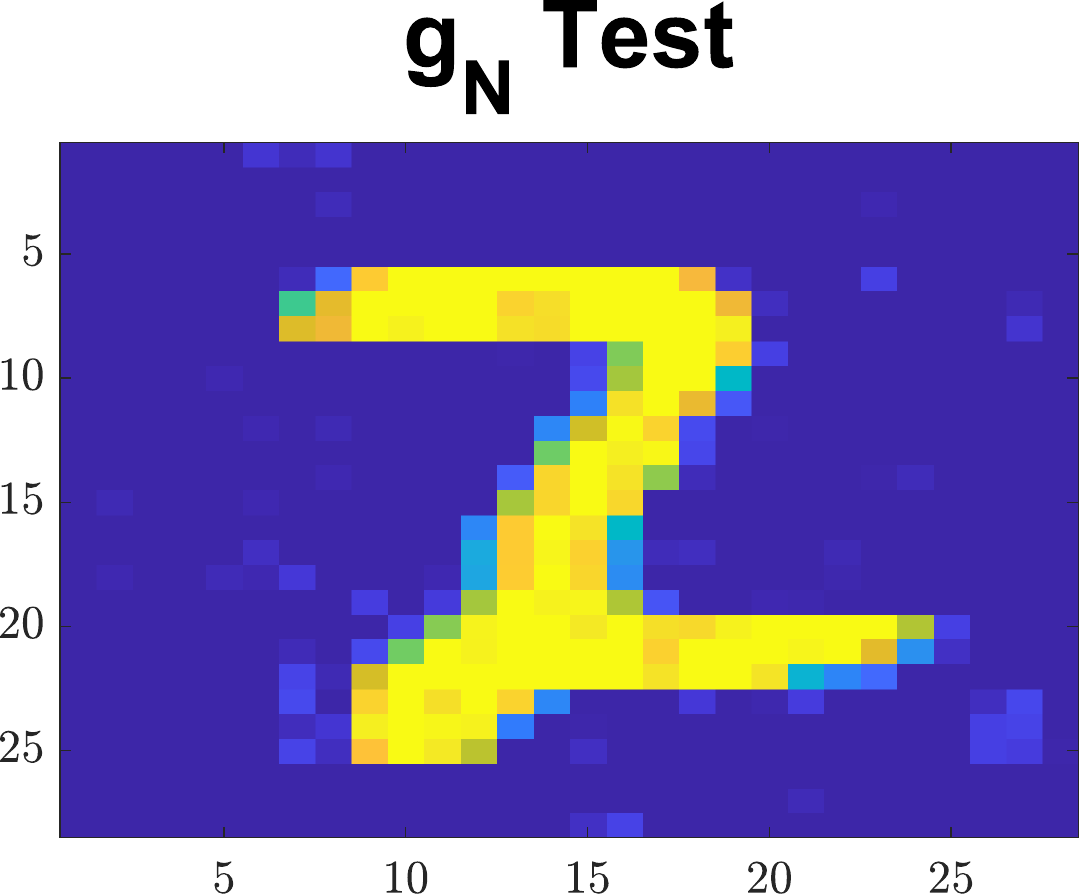}
\hspace{0.8cm}\includegraphics[width=0.14\textwidth, height=0.14\textwidth]{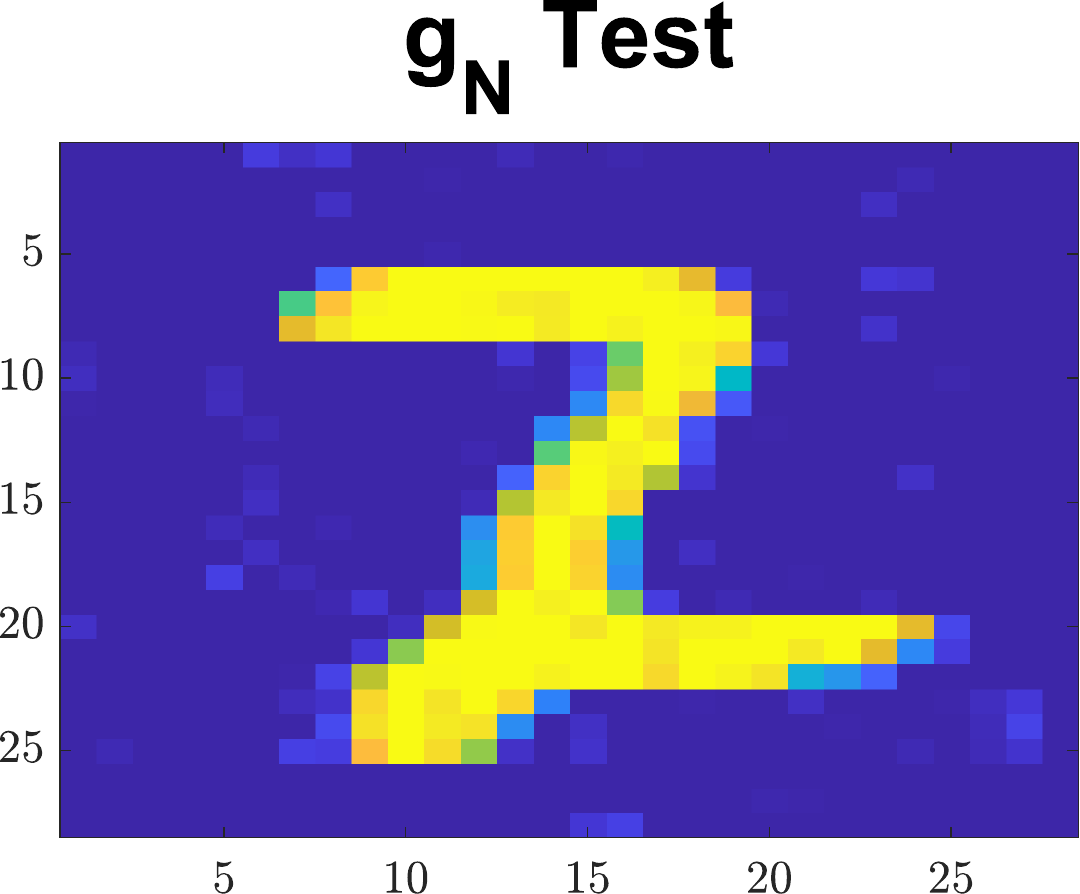}
\hspace{0.8cm}\includegraphics[width=0.15\textwidth, height=0.13\textwidth]{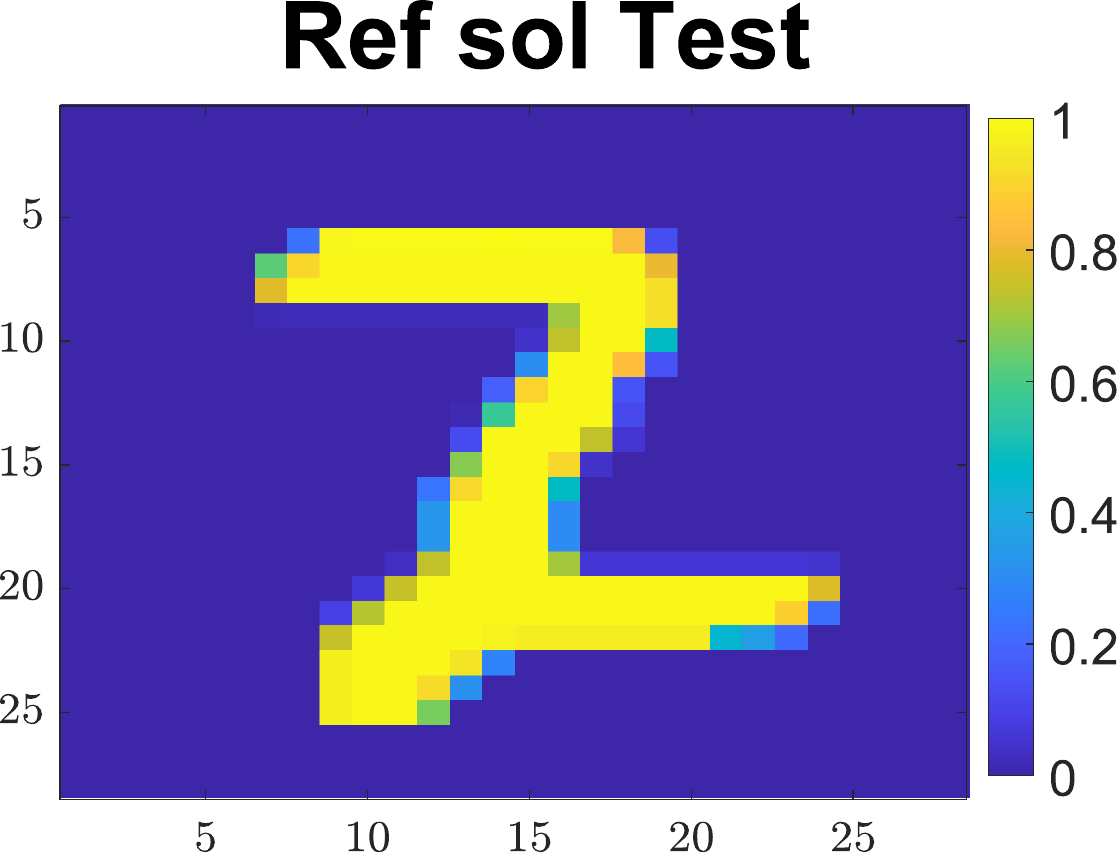}
\caption{\label{f:denoise_all}Reconstructions of representative samples from training (\textit{row} $1$), validation (\textit{row} $2$) and testing (\textit{row} $3$) data showing successful image denoising performed using OCTANE algorithm with $N=\{6,12,20\}$ layers (\textit{columns} $2$-$4$), noisy input data as the initial condition (\textit{column} $1$) and reference solution (\textit{column} $5$).}
\end{center}
\end{figure}

Experiments are conducted for \( N = \{4, 6, 10, 12, 20, 30\} \) layers, with each \( N \) corresponding to a separate training run (see \cref{t:exp_config}). \Cref{f:denoise_all} shows representative reconstructions for \( N = \{6, 12, 20\} \): each row corresponds to training, validation, or testing data; columns show initial condition \( f_0 \), reconstructions \( \hat{\sigma}(g_N) \), and reference images \( x \). Reconstructions are converted to MATLAB arrays for visualization.

Reported testing errors are: \( \alpha_{\text{test}} = \{1.51\text{e}{-3}, 9.32\text{e}{-4}, 6.03\text{e}{-4}\} \), \texttt{PSNR} = \{27.04, 28.75, 30.52\}, and \texttt{SSIM} = \{0.91, 0.94, 0.94\}.

\begin{figure}[h!]
\begin{center}
\includegraphics[width=0.3\textwidth, height=0.3\textwidth]{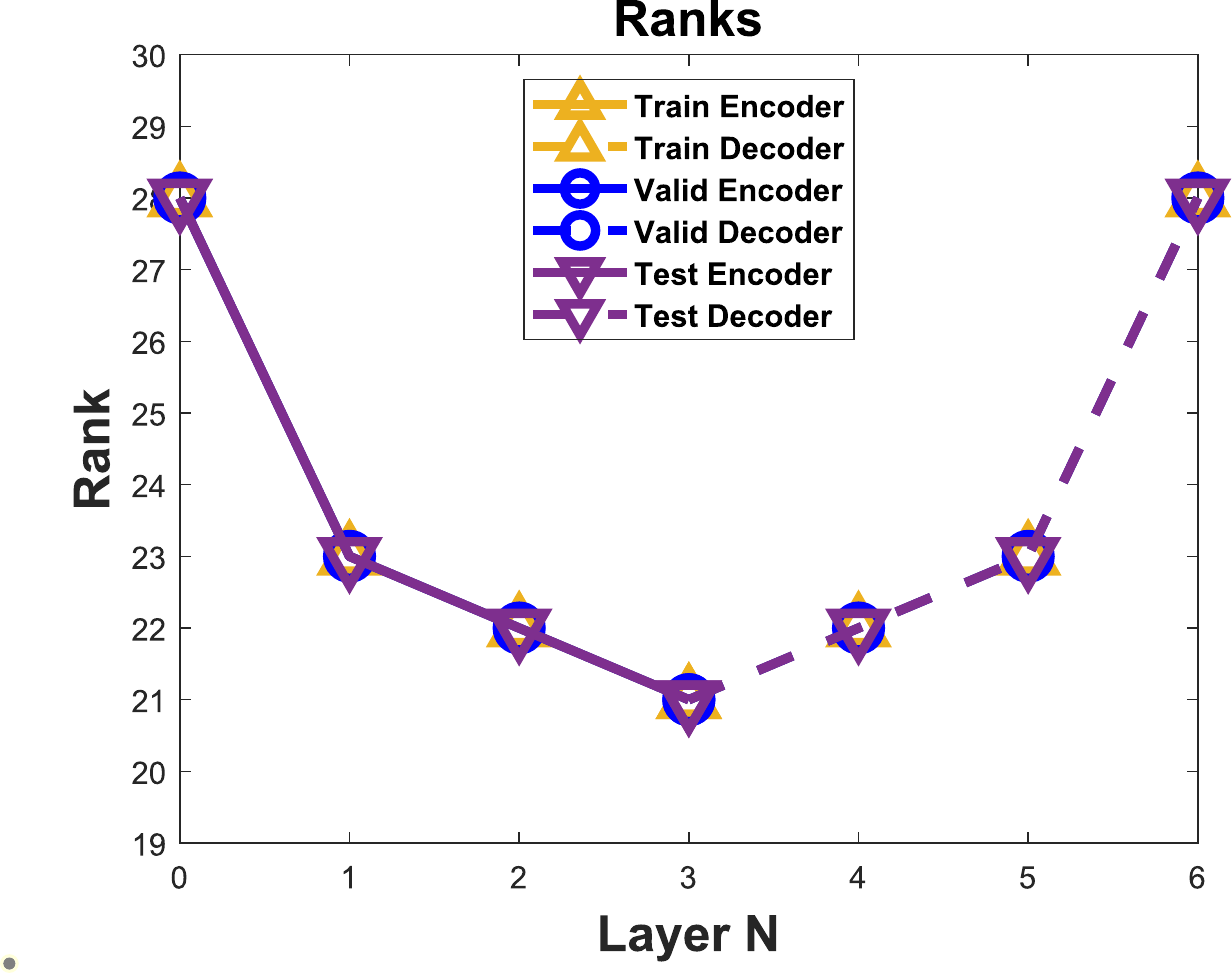} 
\hspace{0.5cm} \includegraphics[width=0.3\textwidth, height=0.3\textwidth]{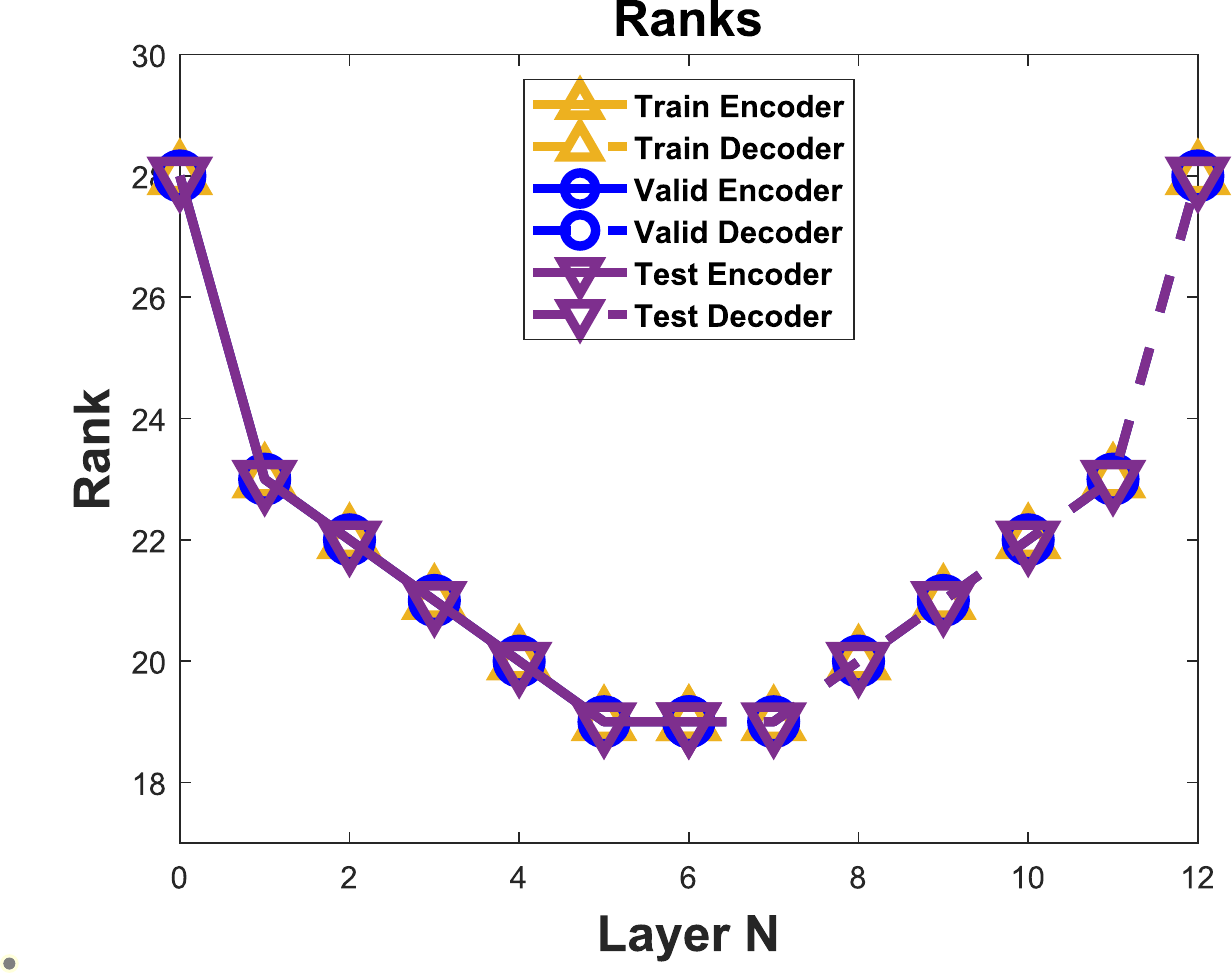}
\hspace{0.5cm} \includegraphics[width=0.3\textwidth, height=0.3\textwidth]{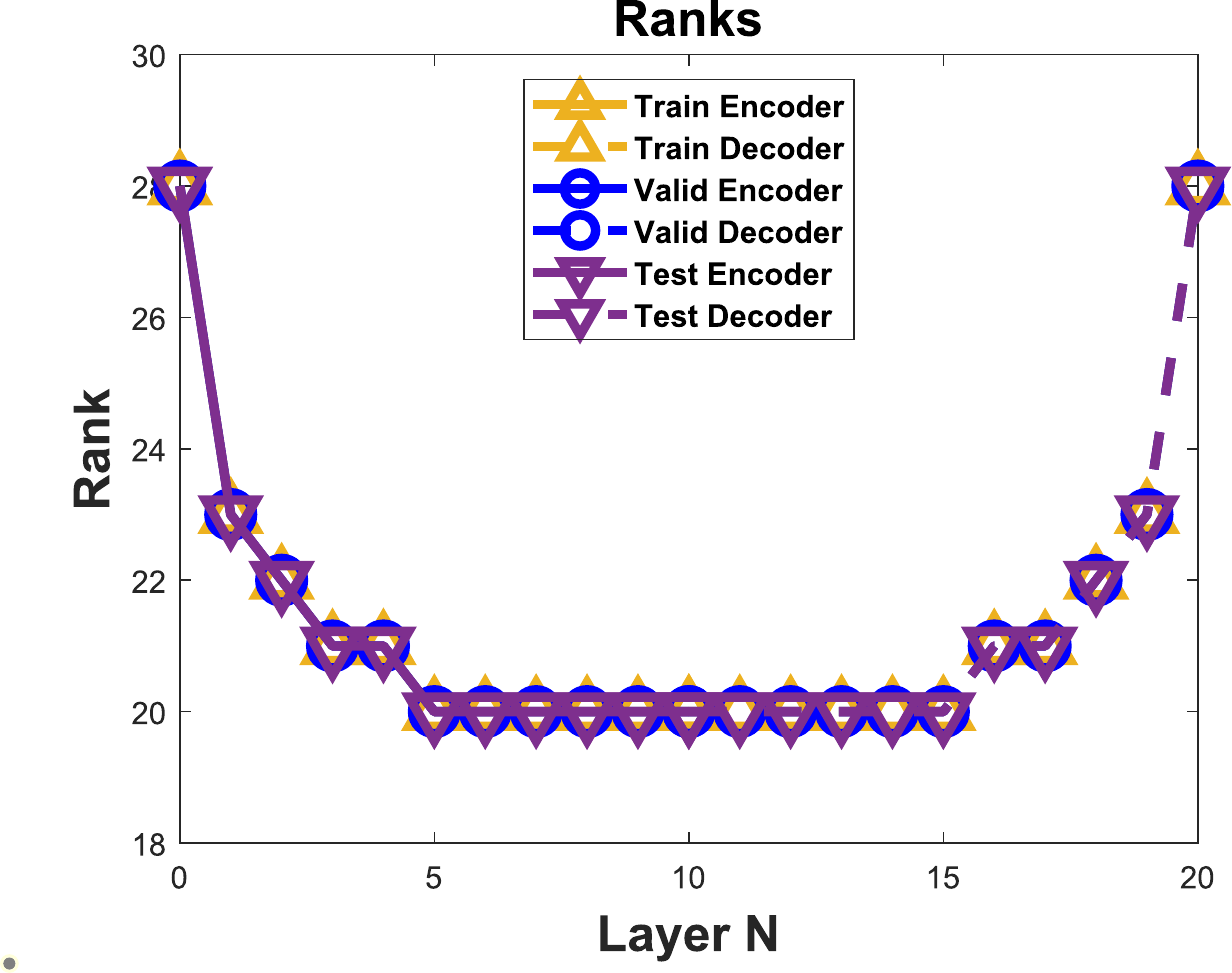} 
\caption{\label{f:denoise_ranks}Rank distributions obtained via OCTANE algorithm for various number of layers $N=6$, $N=12$, and $N=20$ in the image denoising task. Solid lines are encoder ranks and dotted lines are decoder ranks. Note that compression in the network is in the context of rank reduction.}
\end{center}
\end{figure}

A key feature of the OCTANE algorithm is its ability to adapt ranks across layers for a fixed \( N \). \Cref{f:denoise_ranks} shows the rank distributions during forward propagation for the denoising experiments in \cref{f:denoise_all}. The rank at layer 0 corresponds to the input \( f_0 \), at \( N/2 \) to the encoder-decoder interface \( f_{N_e}=g_0 \), and at layer \( N \) to the final output \( \hat{\sigma}(g_N) \). Solid lines represent encoder ranks (state \( f \)), and dotted lines denote decoder ranks (state \( g \)). Colors indicate training (orange), validation (blue), and testing (purple) data; testing ranks are shown for the last mini-batch of 20 images.

These profiles reveal the emergent optimal rank architecture learned by the model, with ranks effectively serving as layer widths (cf. \cref{f:OCTANEnet}).

\begin{figure}[h!]
\begin{center}
\includegraphics[width=0.3\textwidth, height=0.3\textwidth]{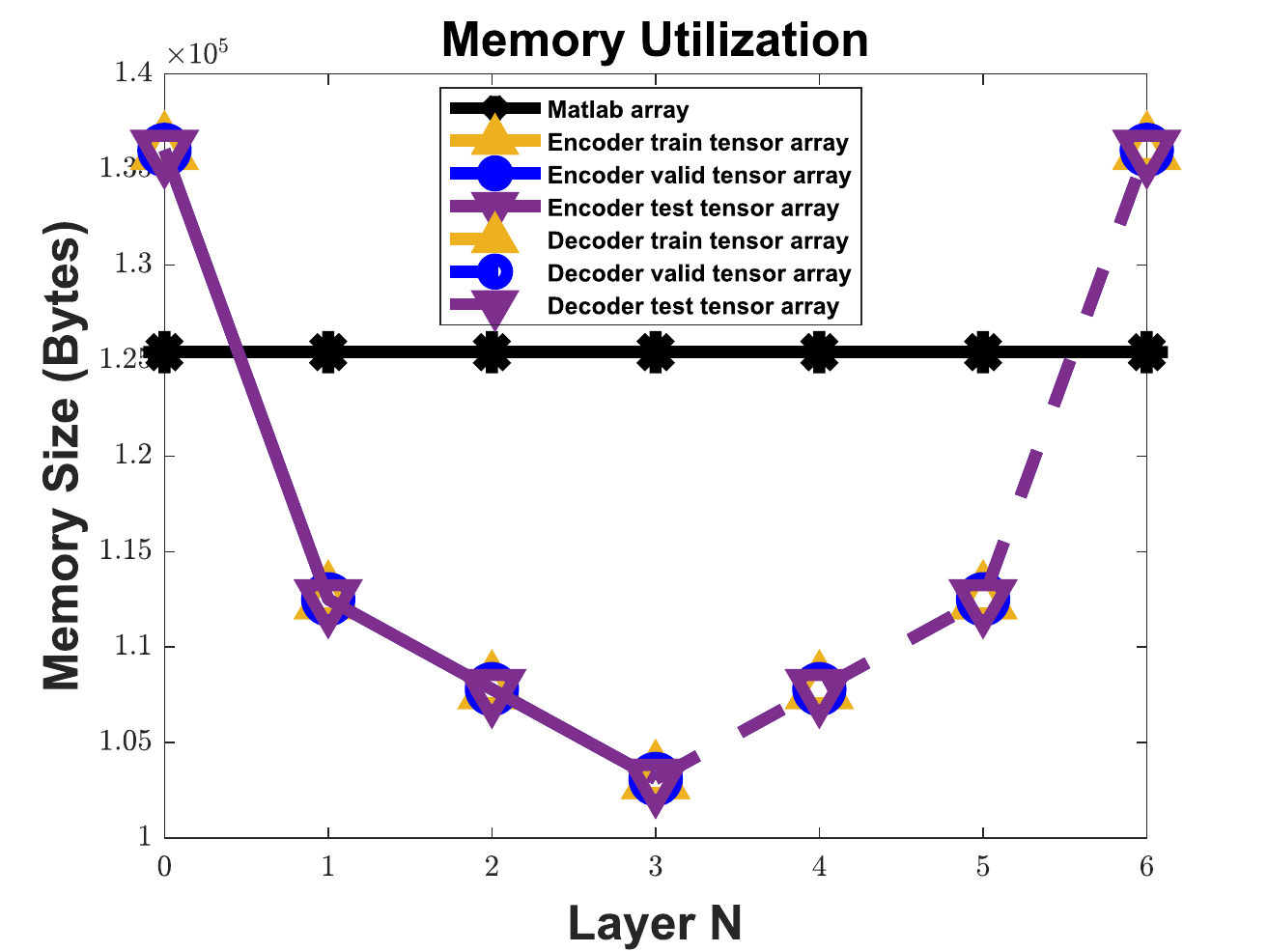} 
\hspace{0.5cm} \includegraphics[width=0.3\textwidth, height=0.3\textwidth]{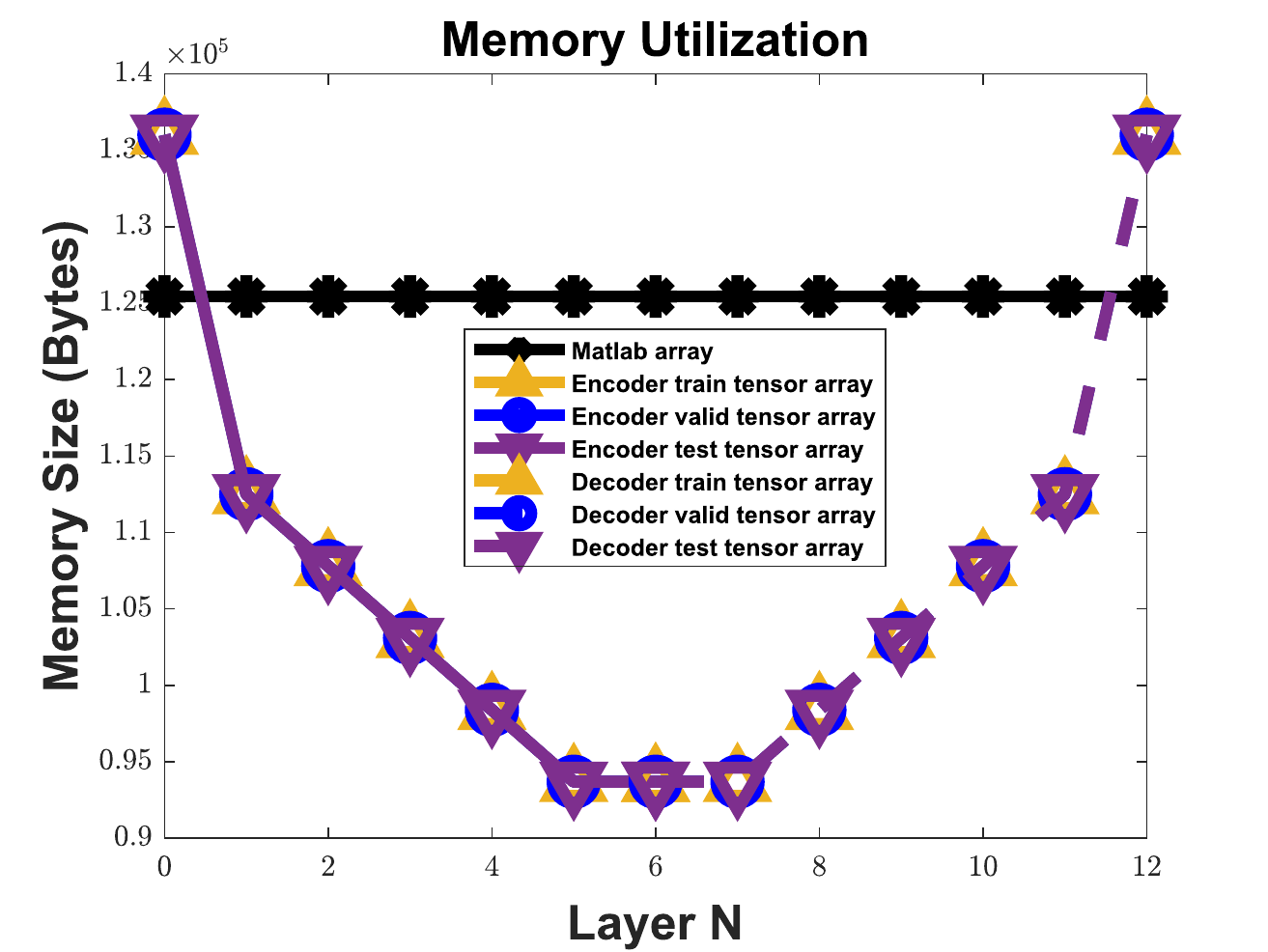}
\hspace{0.5cm} \includegraphics[width=0.3\textwidth, height=0.3\textwidth]{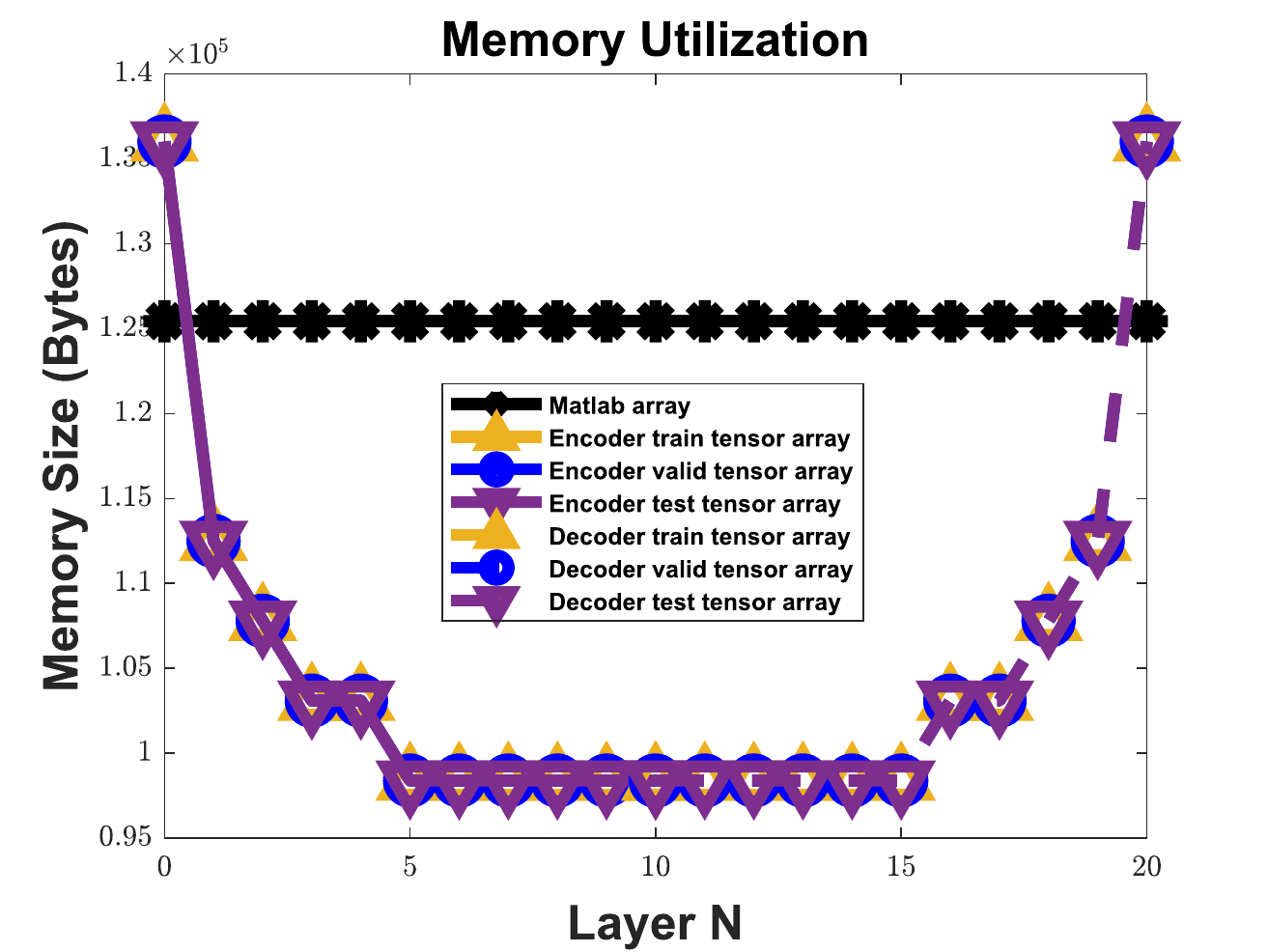} 
\caption{\label{f:denoise_memory}Comparison of memory size in bytes utilized by data in each layer stored as a tensor vs. MATLAB array (non-tensor version) for various number of layers for the image denoising task. Note that the OCTANE architecture saves an average of $7.10\%$ ($N=6$ layers), $14.35\%$ ($N=12$ layers), and $16.21\%$ ($N=20$) memory for each data type (training/validation/testing). Owing to rank-reduction, solving differential equation on the low-rank manifold saves significant memory.}
\end{center}
\end{figure}

By solving the differential equations on a low-rank tensor manifold, OCTANE significantly reduces memory usage. \Cref{f:denoise_memory} compares memory consumption of state variables (training: orange, validation: blue, testing: purple) in tensor format against standard MATLAB arrays (black), as used in conventional RNNs. Total memory is computed as the sum across all layers, and similar trends are expected for adjoint variables. The results highlight OCTANE's efficiency: average memory savings are approximately $7.10\%$ for \( N=6 \), $14.35\%$ for \( N=12 \), and $16.21\%$ for \( N=20 \), across all datasets.

\begin{figure}[h!]
\begin{center}
\includegraphics[width=0.3\textwidth, height=0.3\textwidth]{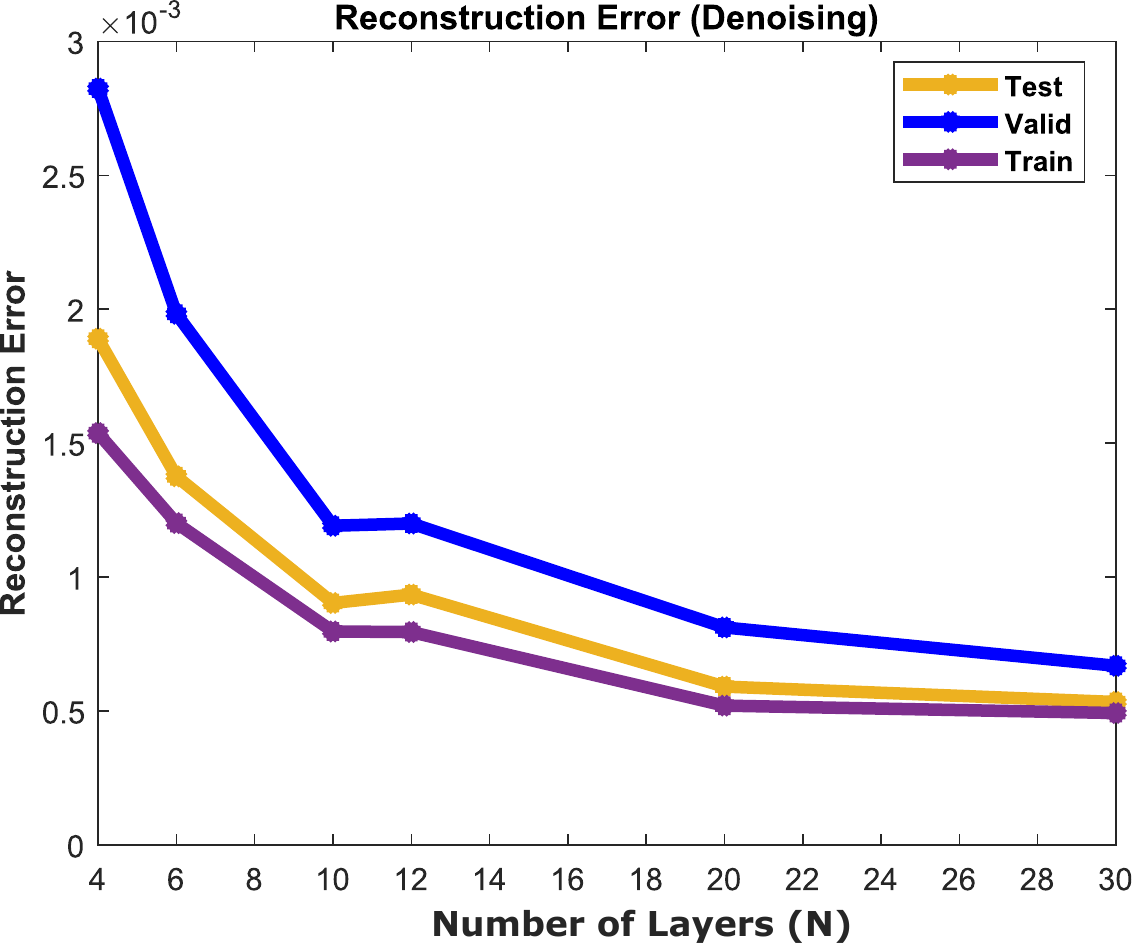} 
\hspace{0.5cm} \includegraphics[width=0.3\textwidth, height=0.3\textwidth]{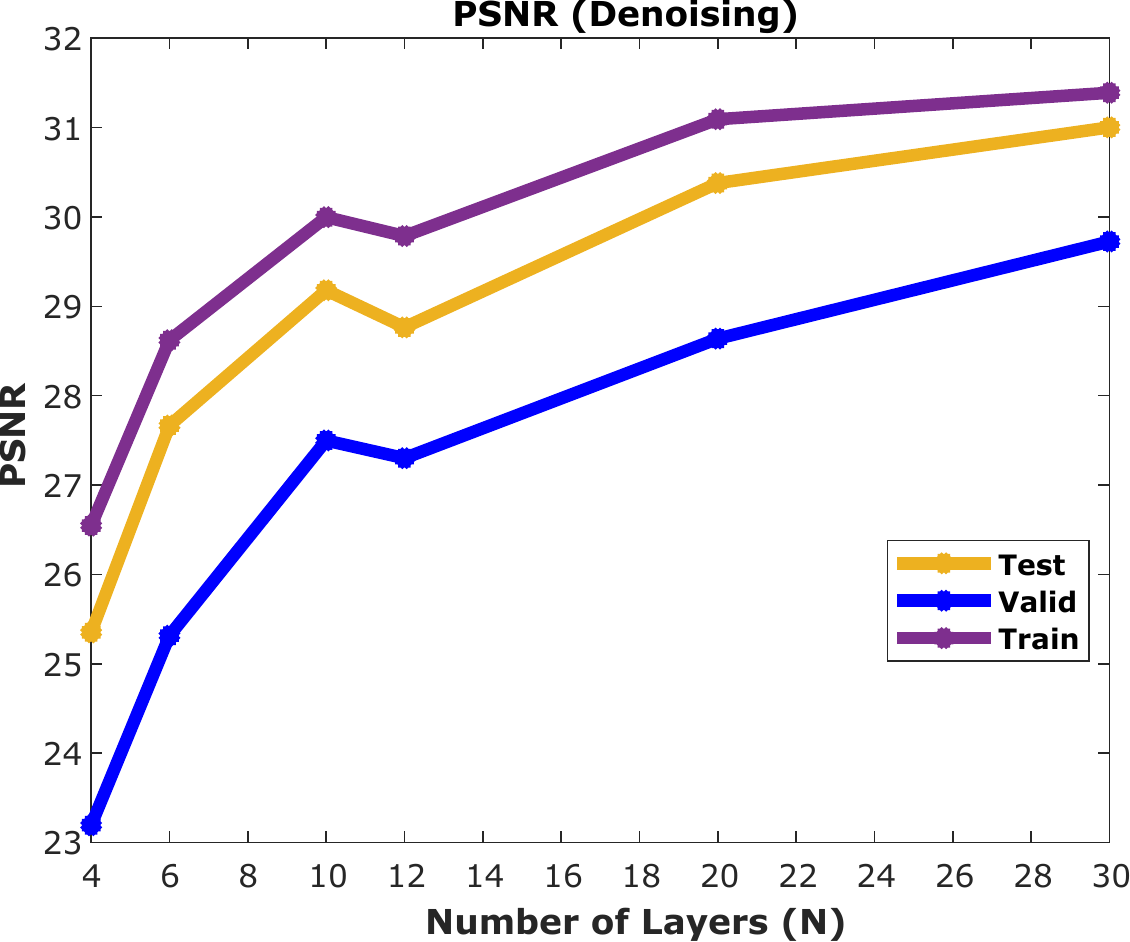}
\hspace{0.5cm} \includegraphics[width=0.31\textwidth, height=0.3\textwidth]{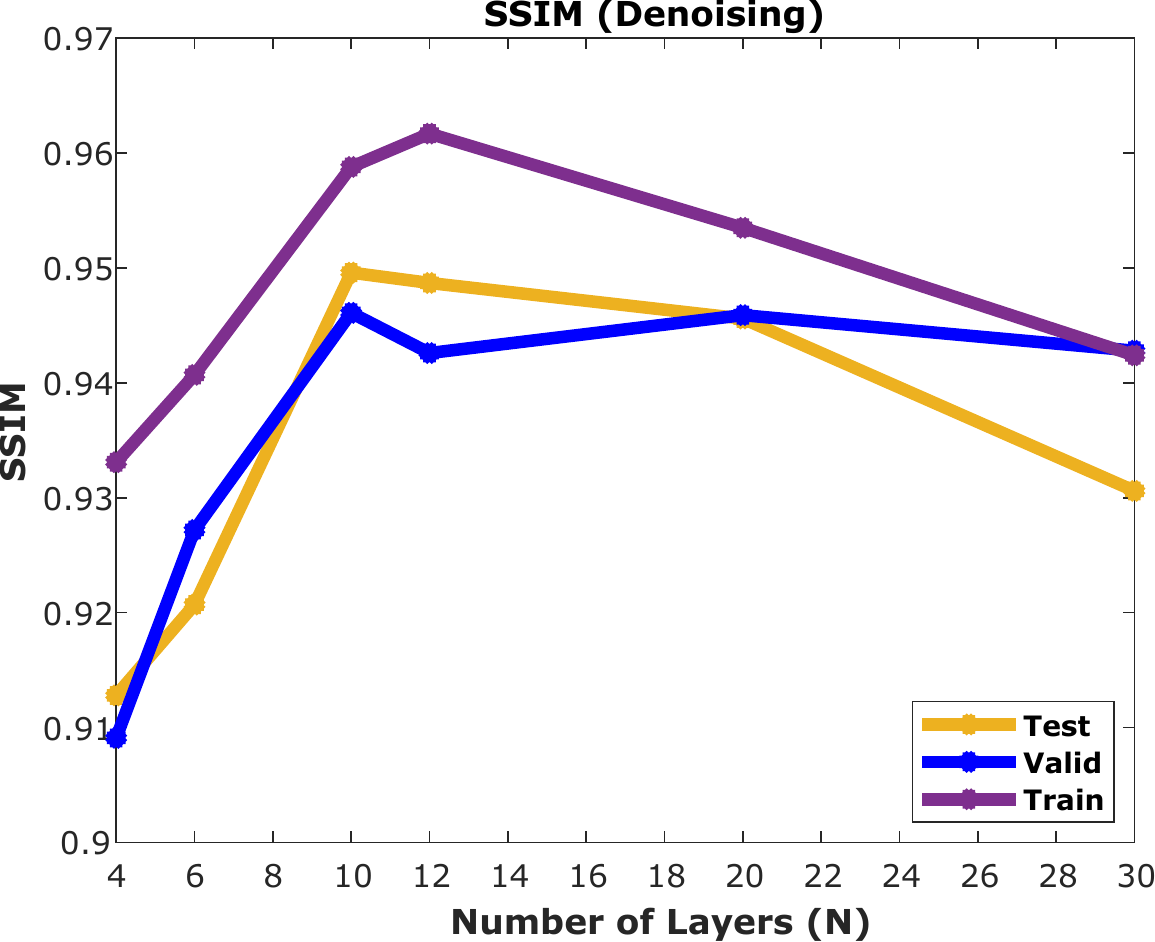}
\caption{\label{f:denoise_results_per_N}Plots of the reconstruction error, \texttt{PSNR}, and \texttt{SSIM} against the number of layers in the image denoising experiment. More layers generally corresponds to improvement in reconstruction error, but the marginal improvement comes at a relatively greater computational cost. \texttt{SSIM} asserts larger number of layers may not be necessary. \texttt{PSNR} follows the same trend as the reconstruction error.}
\end{center}
\end{figure}

A key question is how many layers are needed for effective reconstruction. In \cref{f:denoise_results_per_N}, we report the reconstruction error \cref{recon_err}, \texttt{PSNR}, and \texttt{SSIM} for $N = \{4, 6, 10, 12, 20, 30\}$ layers across training, validation, and test datasets. Each point on the x-axis represents an independent training run. Reconstruction error (lower is better) is based on MSE, while \texttt{PSNR} and \texttt{SSIM} (higher is better) are computed using MATLAB’s built-in functions. Training and validation metrics are averaged over 20 samples; testing over 1000 samples.

From the reconstructions and similarity metrics, it is evident that OCTANE effectively performs image denoising. As expected, increasing the number of layers improves reconstructions, since a larger $N$ implies a smaller Euler step-size $\tau$, yielding finer discretizations and more stable solutions.

However, while the reconstruction error decreases with $N$, the improvement from $N=12$ to $N=30$ is marginal (on the order of $1e\text{-}4$), despite significantly increased computational cost. Interestingly, \texttt{SSIM} peaks at $N=12$ and slightly deteriorates for larger $N$, with a drop of about $1e\text{-}2$ between $N=12$ and $N=30$ in testing. \texttt{PSNR} trends mirror reconstruction error. These observations suggest that while deeper networks reduce error, high-quality reconstructions can be obtained with fewer layers at lower computational cost. This motivates a deeper investigation into optimal choices of $N$, $\tau$, and $T$ in \cref{Sec:TauN}. All experiments complete within 40 minutes.

\subsection{Image Deblurring} \label{Sec:Deblurring}

In this experiment, we use the OCTANE algorithm to deblur images from the MNIST dataset (e.g., digit $2$). Synthetic blurred images $\hat{x}$ are generated by applying a Gaussian filter (\texttt{imgaussfilt($x$)} with mean $0$ and standard deviation $1$) to clean images $x$. These serve as initial conditions for \cref{AE_OC}. The autoencoder then learns a low-rank representation of $\hat{x}$ and reconstructs a deblurred image via $g(f(\hat{x}))$, where $g \circ f$ acts as a learned deblurring operator. Experiments follow the configuration in \cref{t:exp_config}.

\begin{figure}[h!]
\begin{center}
\hspace{-0.4cm}\fbox{\text{\textbf{Initial Cond.}}} \hspace{1.2cm}\fbox{$\textbf{N=6}$} \hspace{1.9cm}\fbox{$\textbf{N=12}$} \hspace{1.7cm} \fbox{$\textbf{N=20}$} \hspace{1.3cm}\fbox{\text{\textbf{Ref. Sol.}}}\\~\\

\includegraphics[width=0.14\textwidth, height=0.14\textwidth]{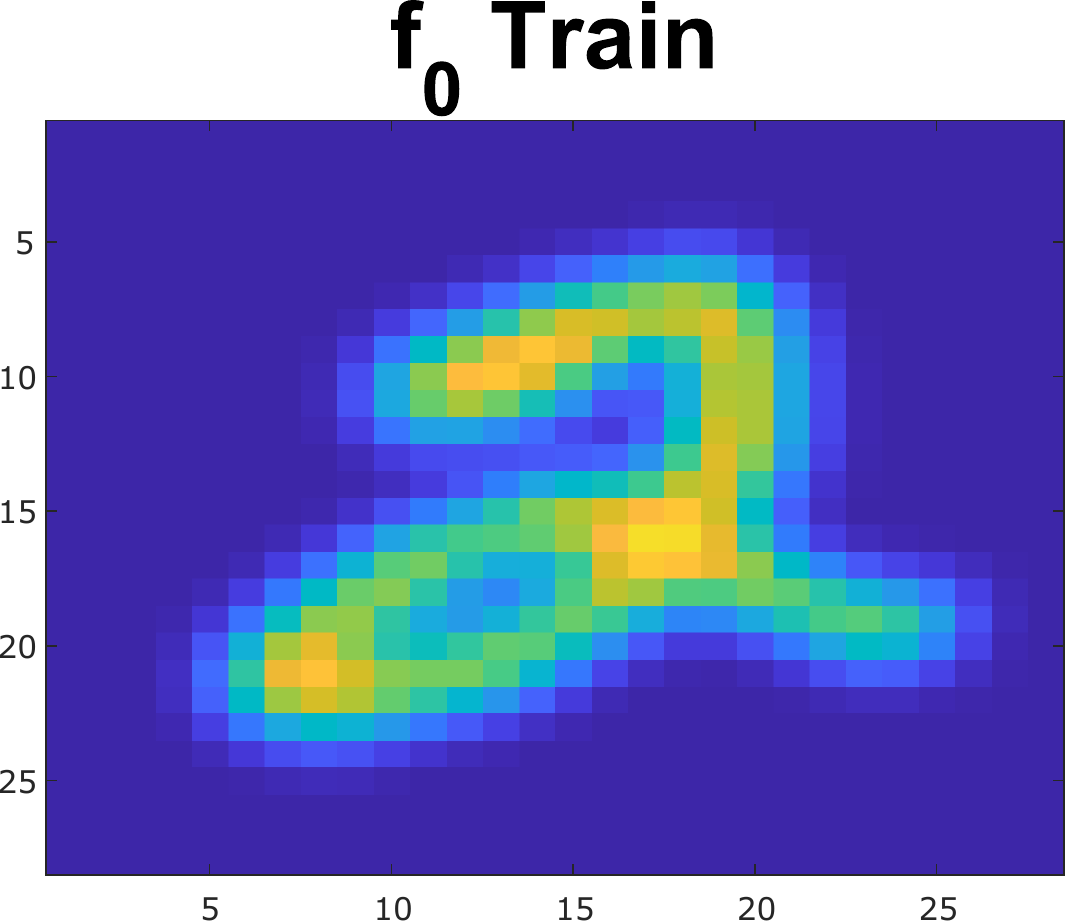} 
\hspace{0.8cm}\includegraphics[width=0.14\textwidth, height=0.14\textwidth]{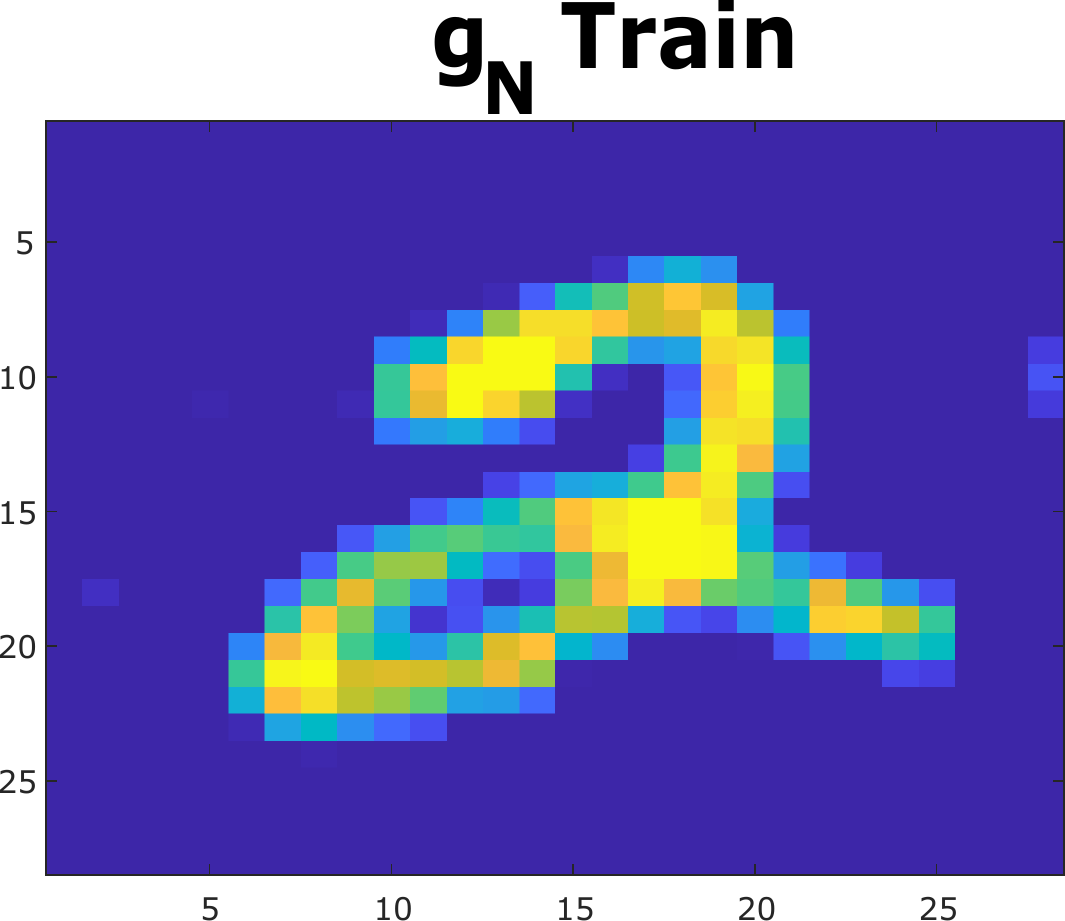}
\hspace{0.8cm}\includegraphics[width=0.14\textwidth, height=0.14\textwidth]{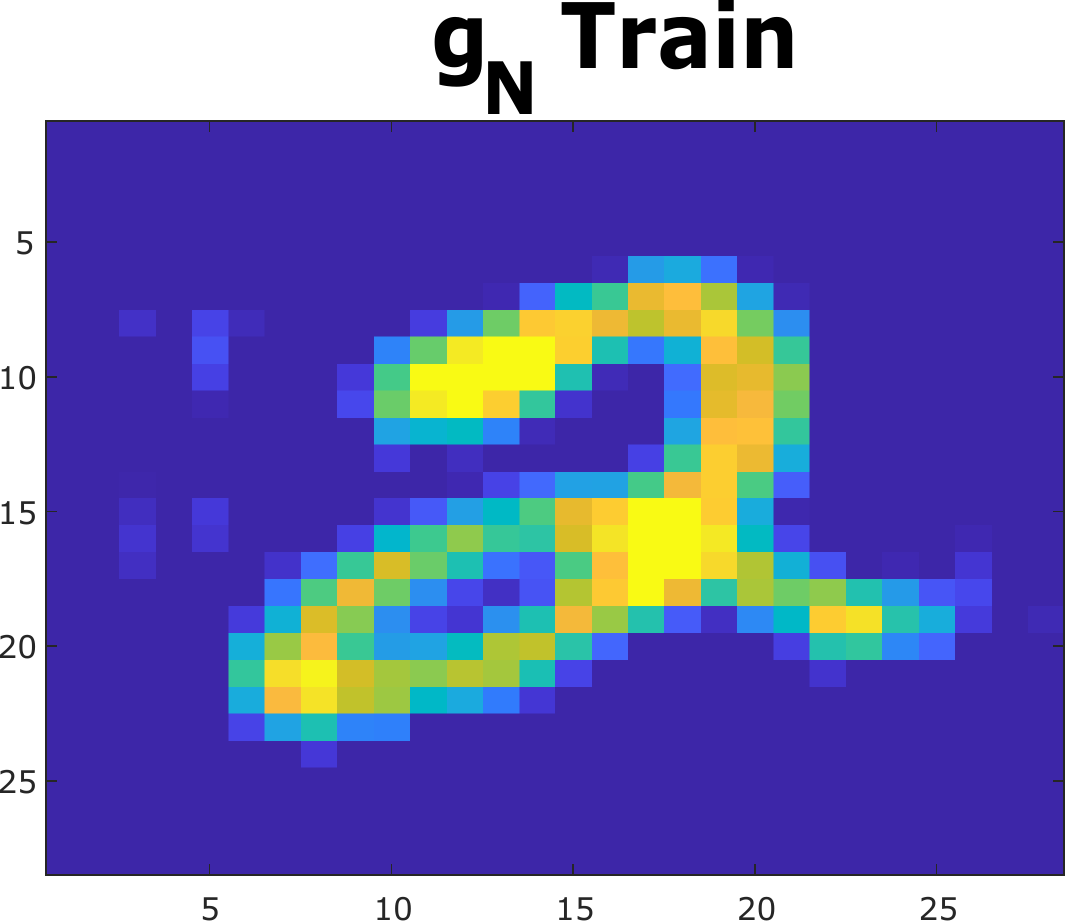}
\hspace{0.8cm}\includegraphics[width=0.14\textwidth, height=0.14\textwidth]{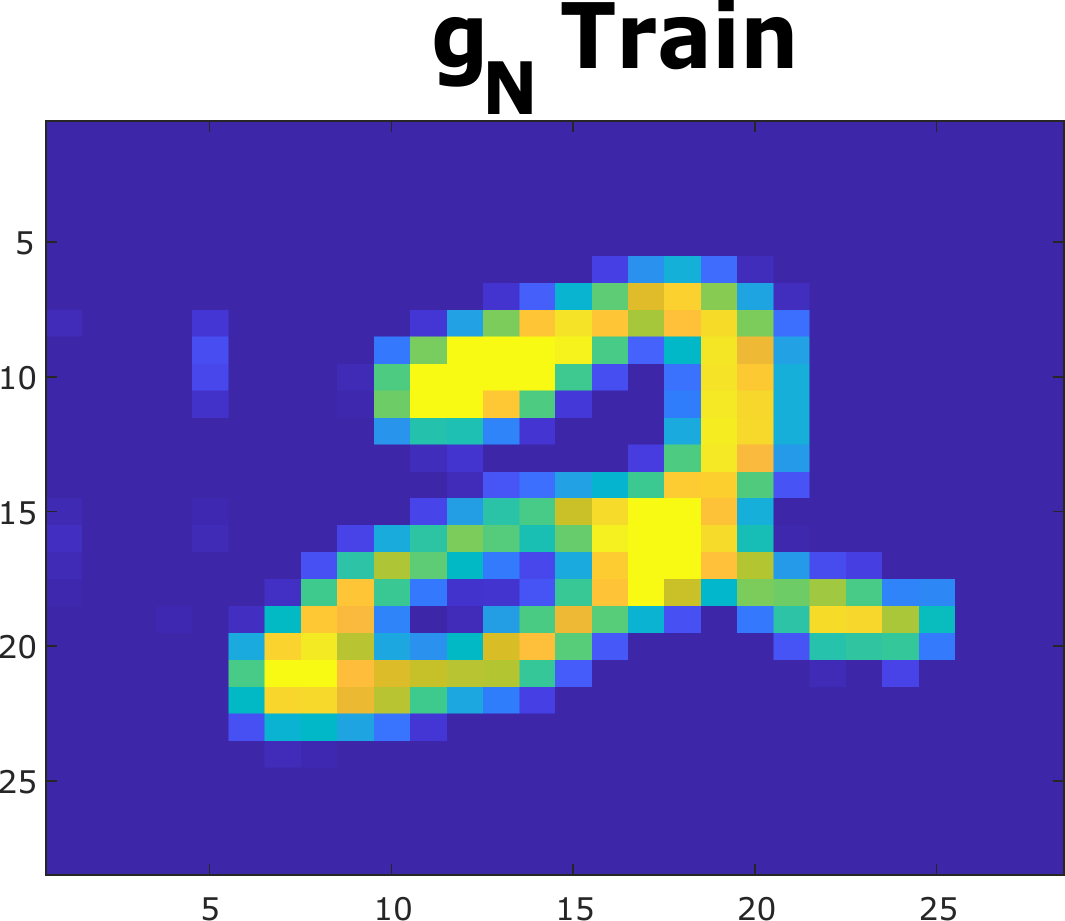}
\hspace{0.8cm}\includegraphics[width=0.15\textwidth, height=0.135\textwidth]{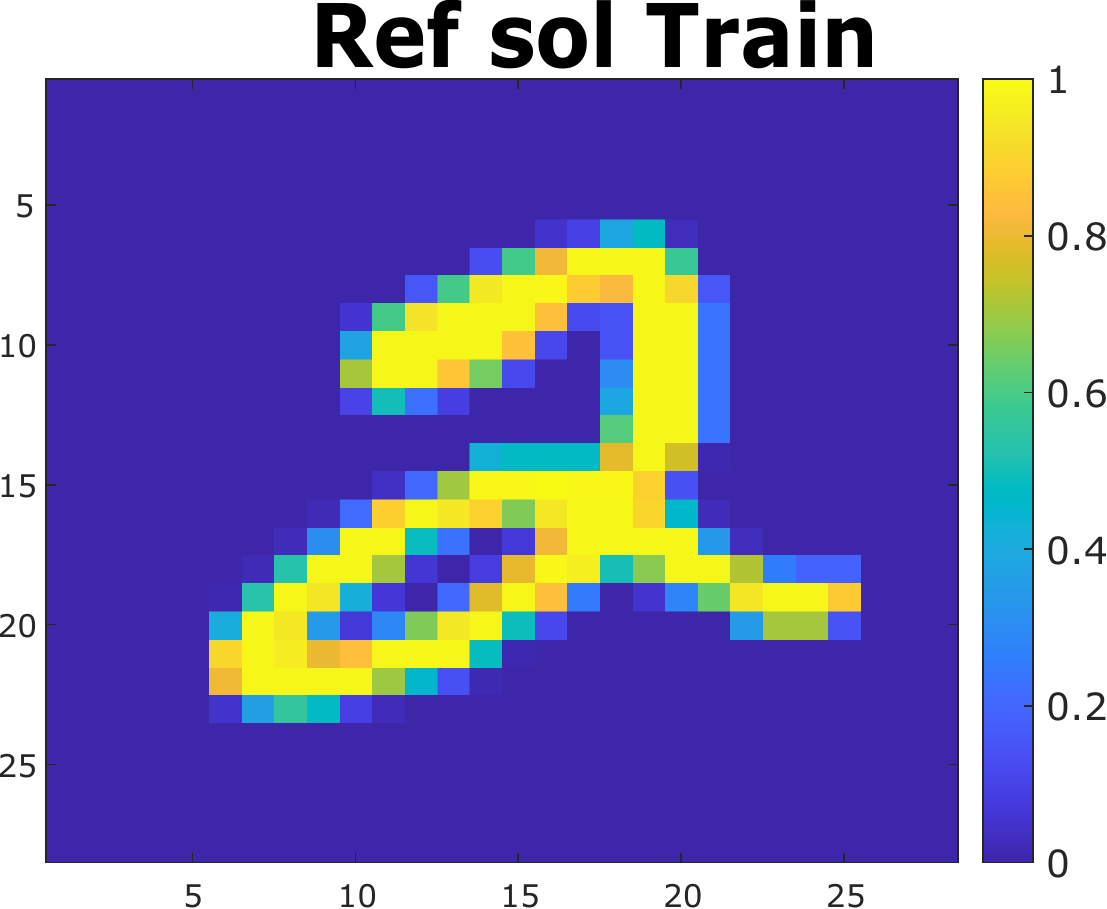}\\
\includegraphics[width=0.14\textwidth, height=0.14\textwidth]{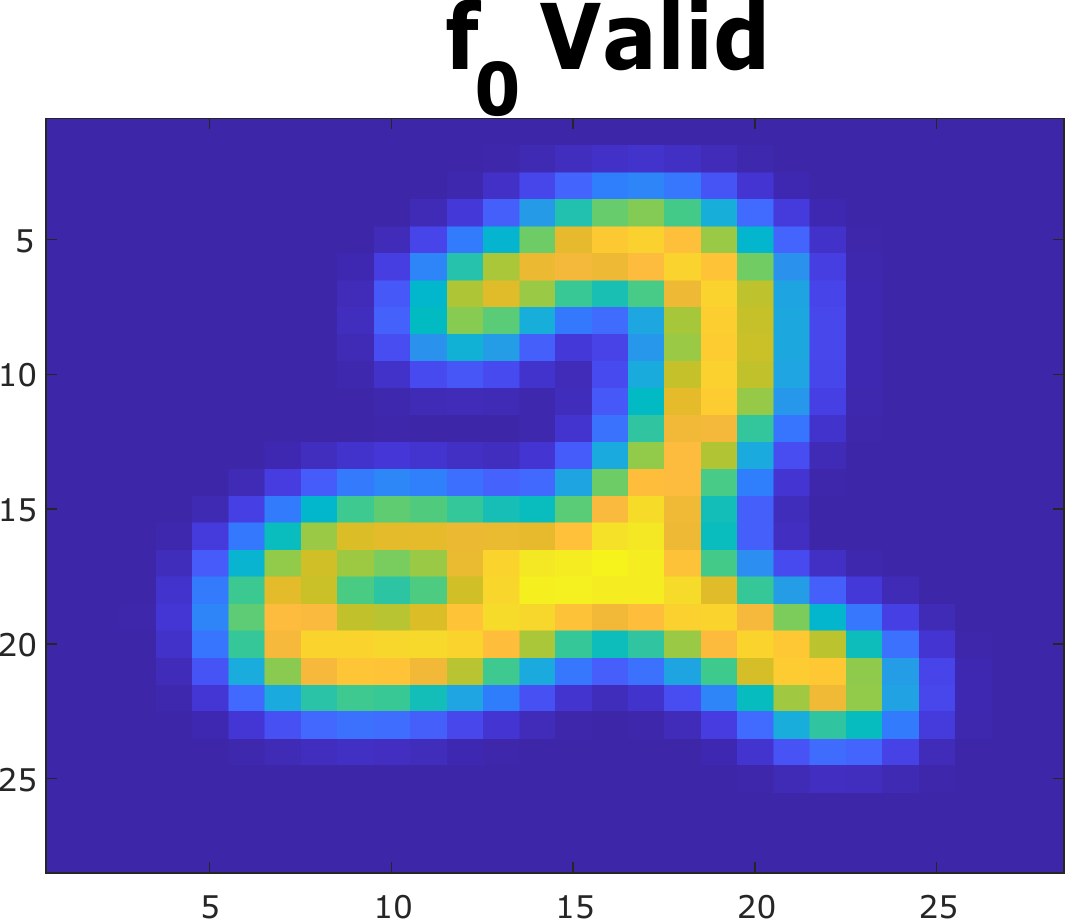} 
\hspace{0.8cm}\includegraphics[width=0.14\textwidth, height=0.14\textwidth]{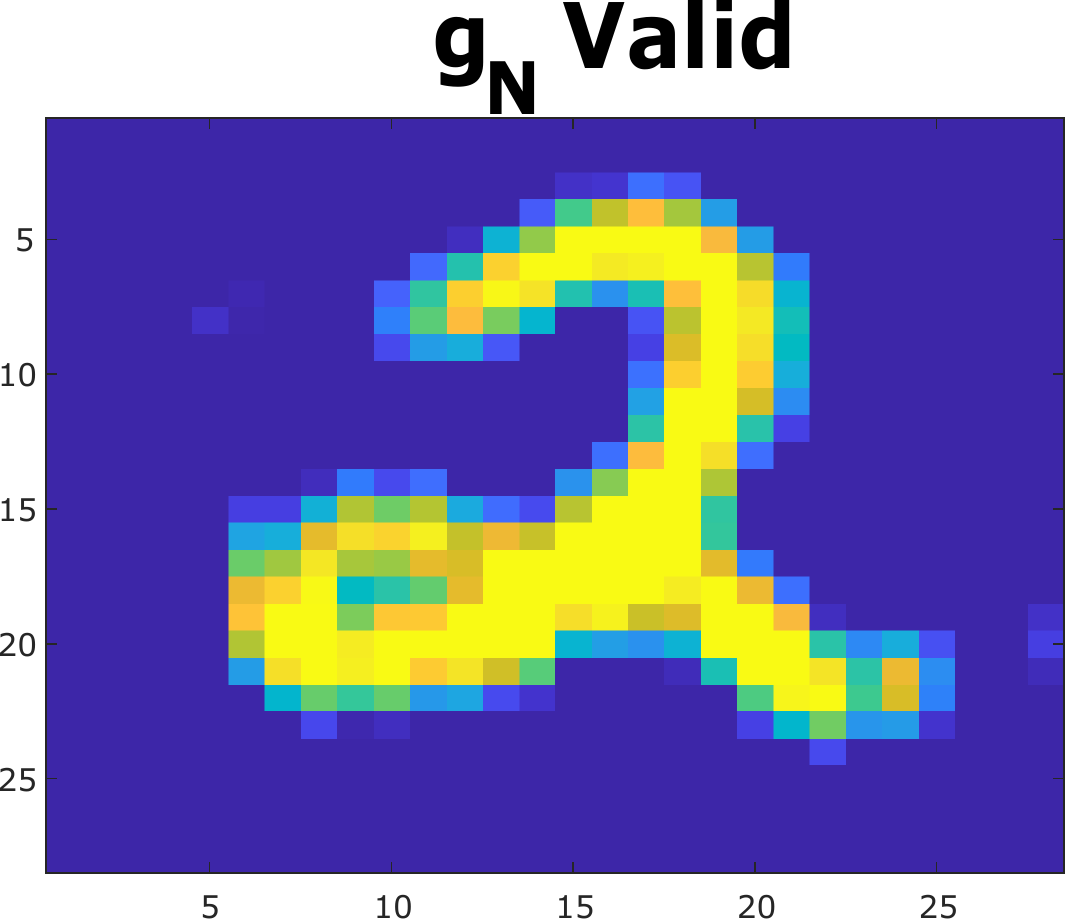}
\hspace{0.8cm}\includegraphics[width=0.14\textwidth, height=0.14\textwidth]{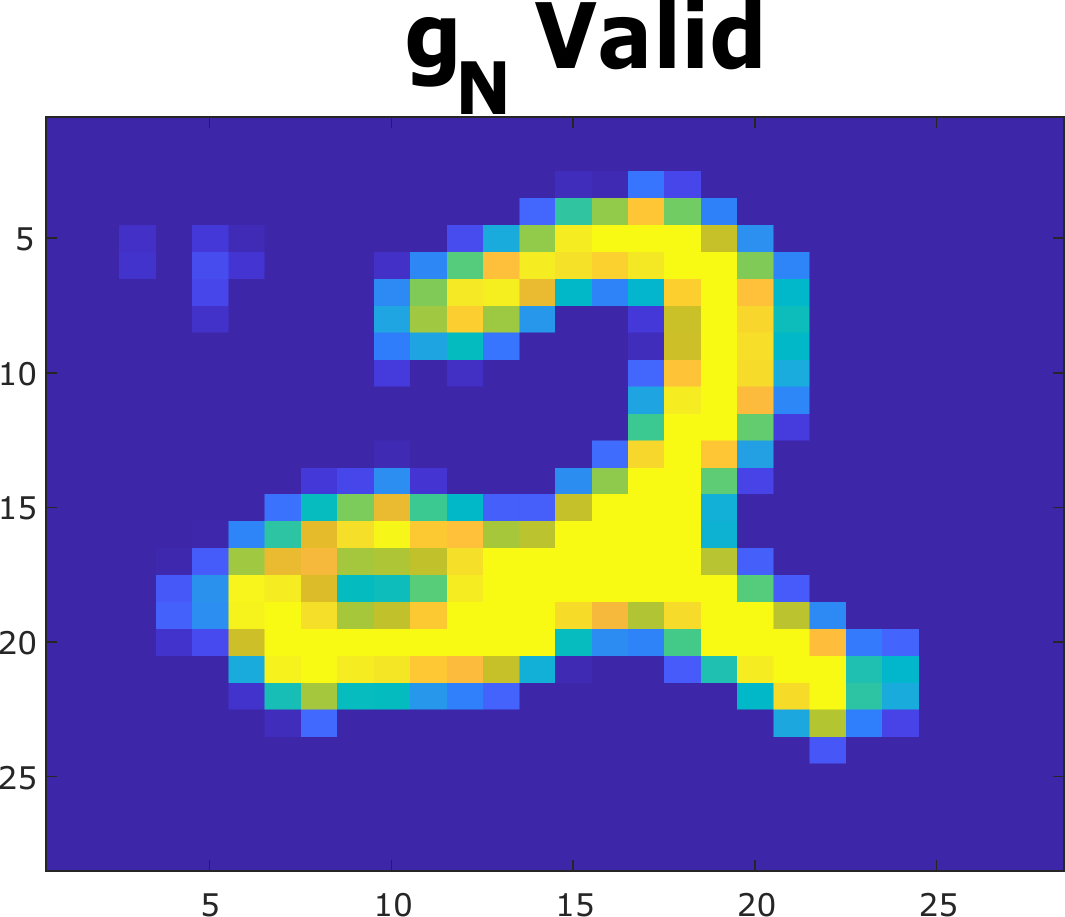}
\hspace{0.8cm}\includegraphics[width=0.14\textwidth, height=0.14\textwidth]{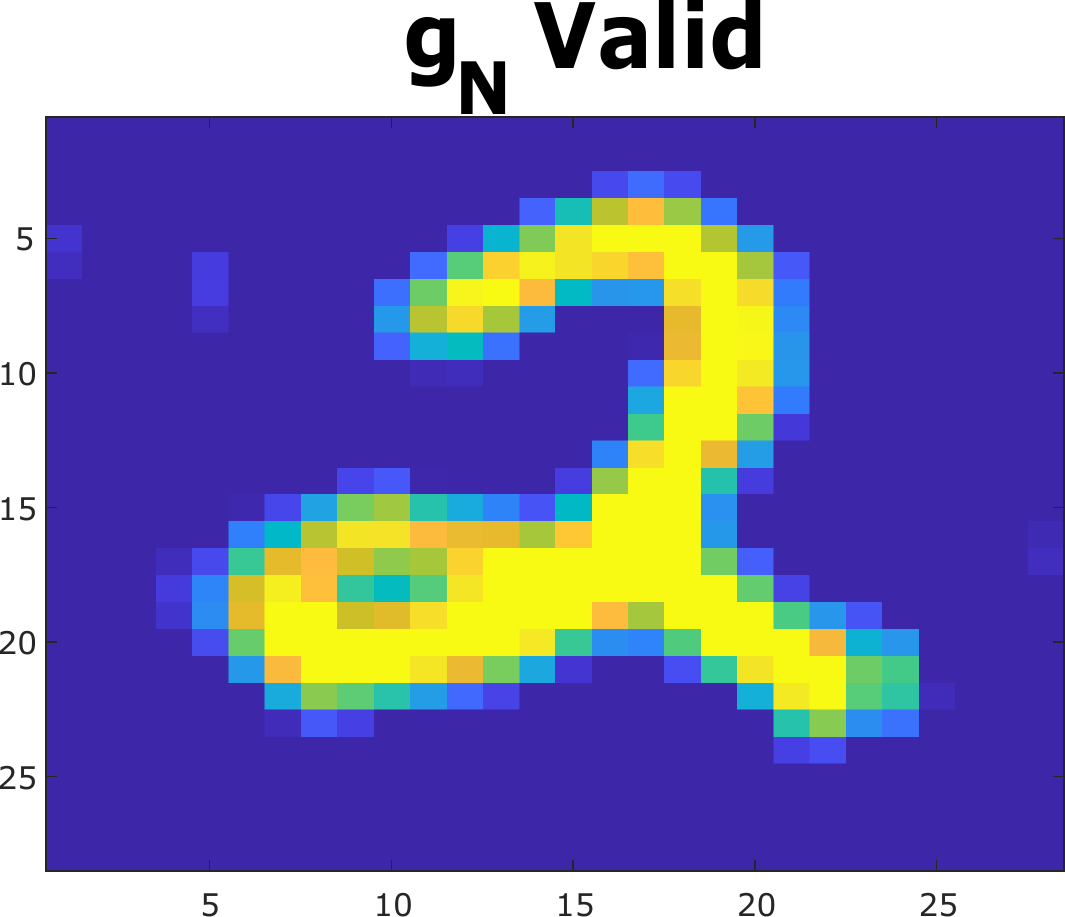}
\hspace{0.8cm}\includegraphics[width=0.15\textwidth, height=0.135\textwidth]{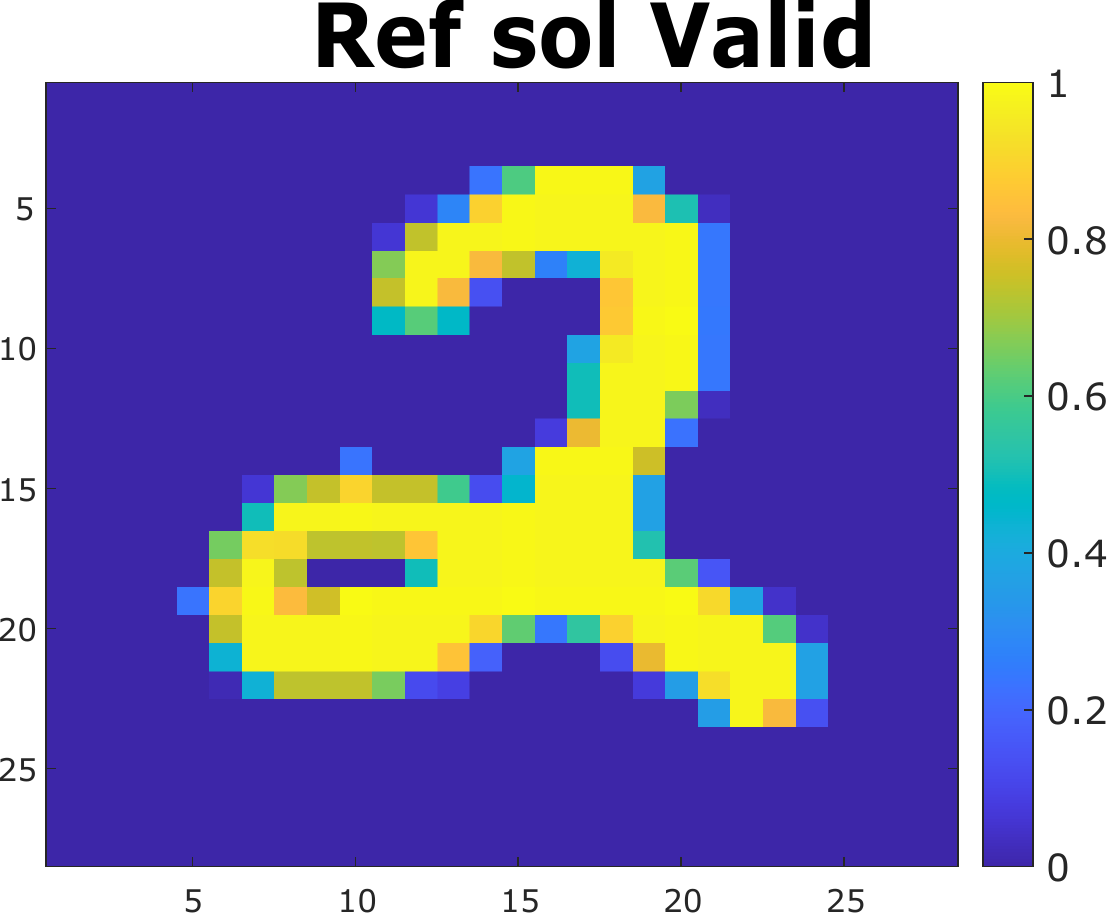}\\
\includegraphics[width=0.14\textwidth, height=0.14\textwidth]{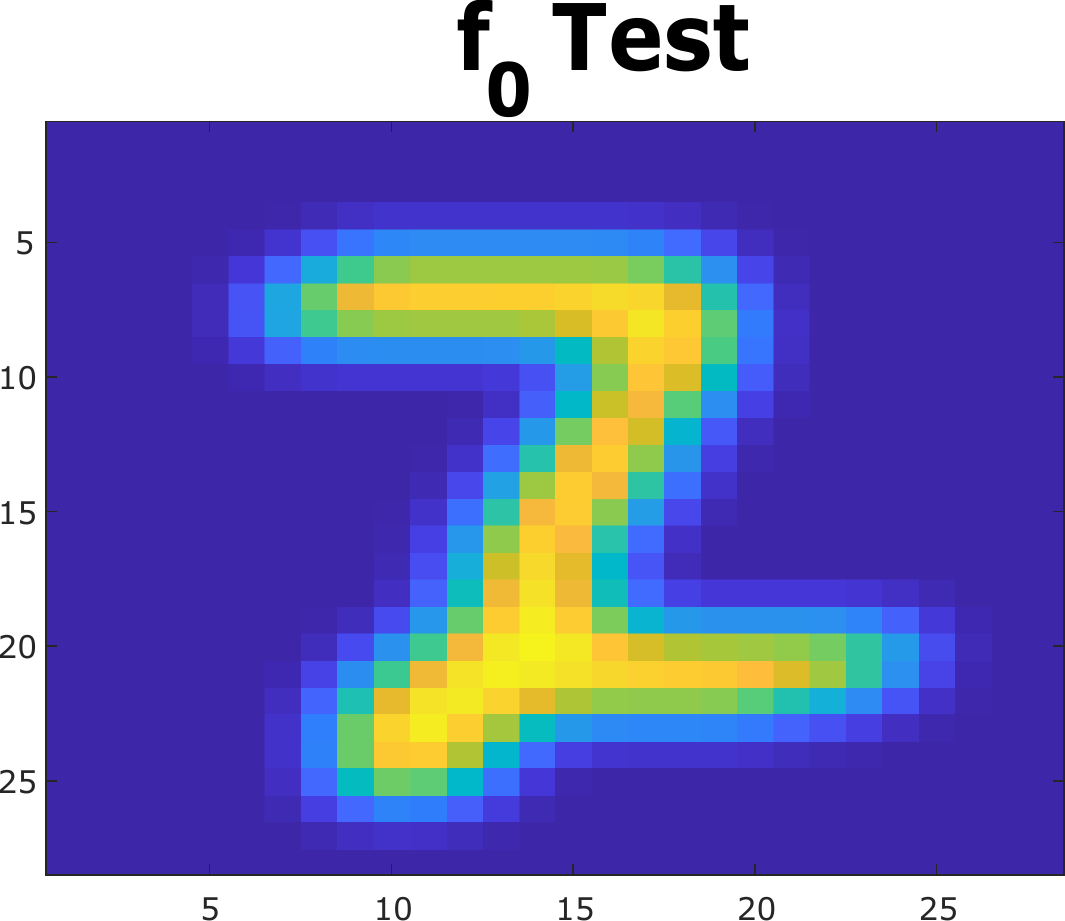} 
\hspace{0.8cm}\includegraphics[width=0.14\textwidth, height=0.14\textwidth]{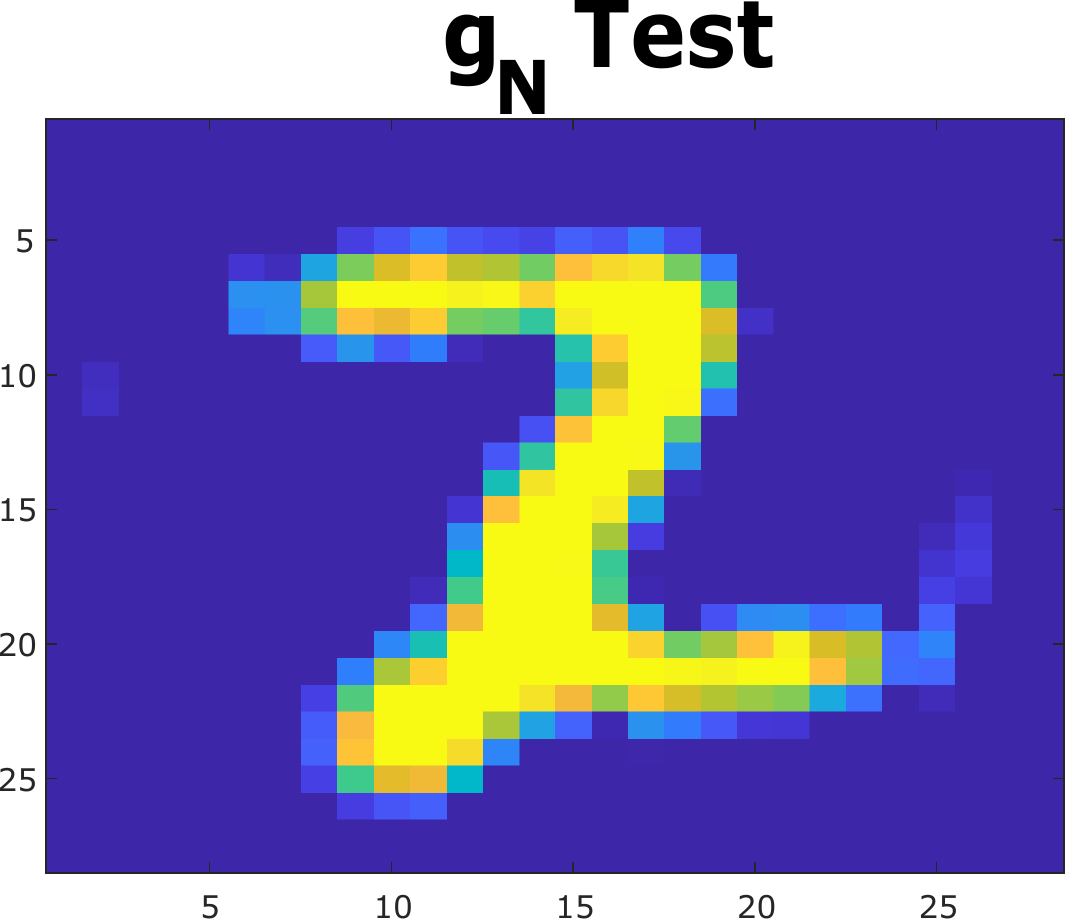}
\hspace{0.8cm}\includegraphics[width=0.14\textwidth, height=0.14\textwidth]{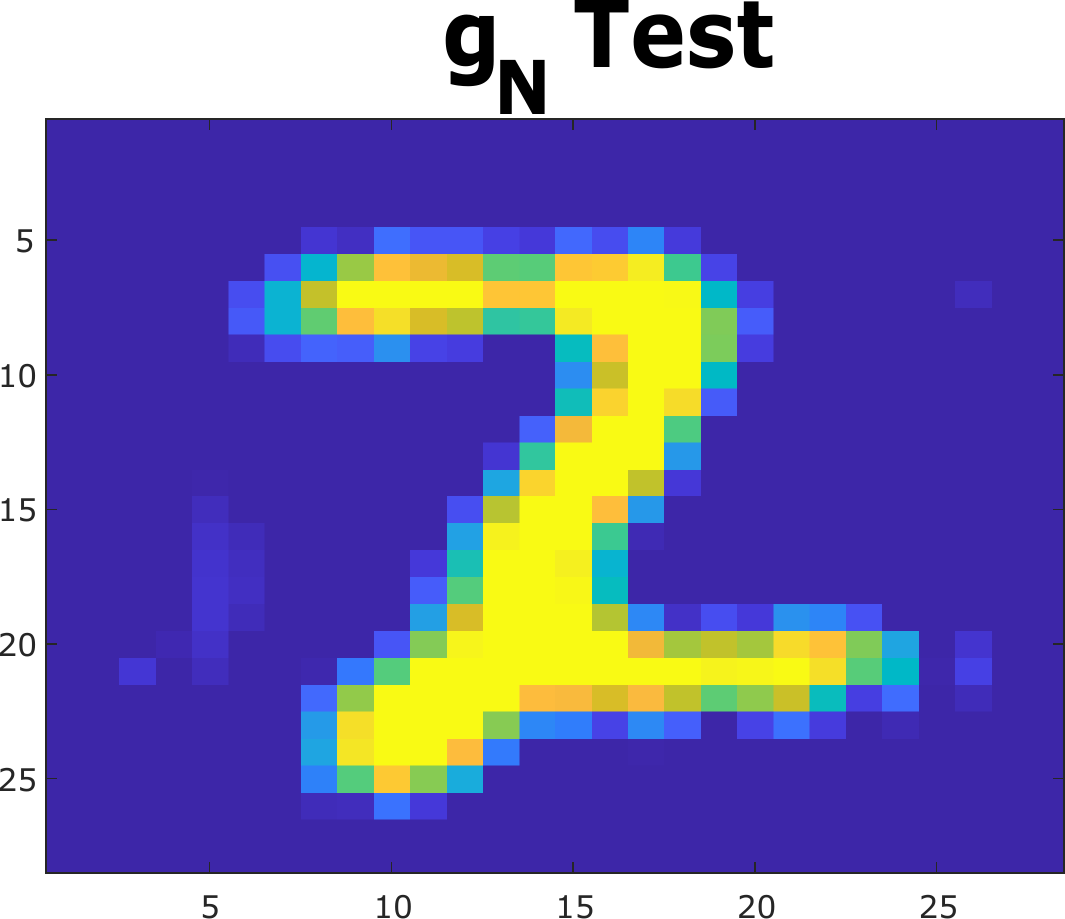}
\hspace{0.8cm}\includegraphics[width=0.14\textwidth, height=0.14\textwidth]{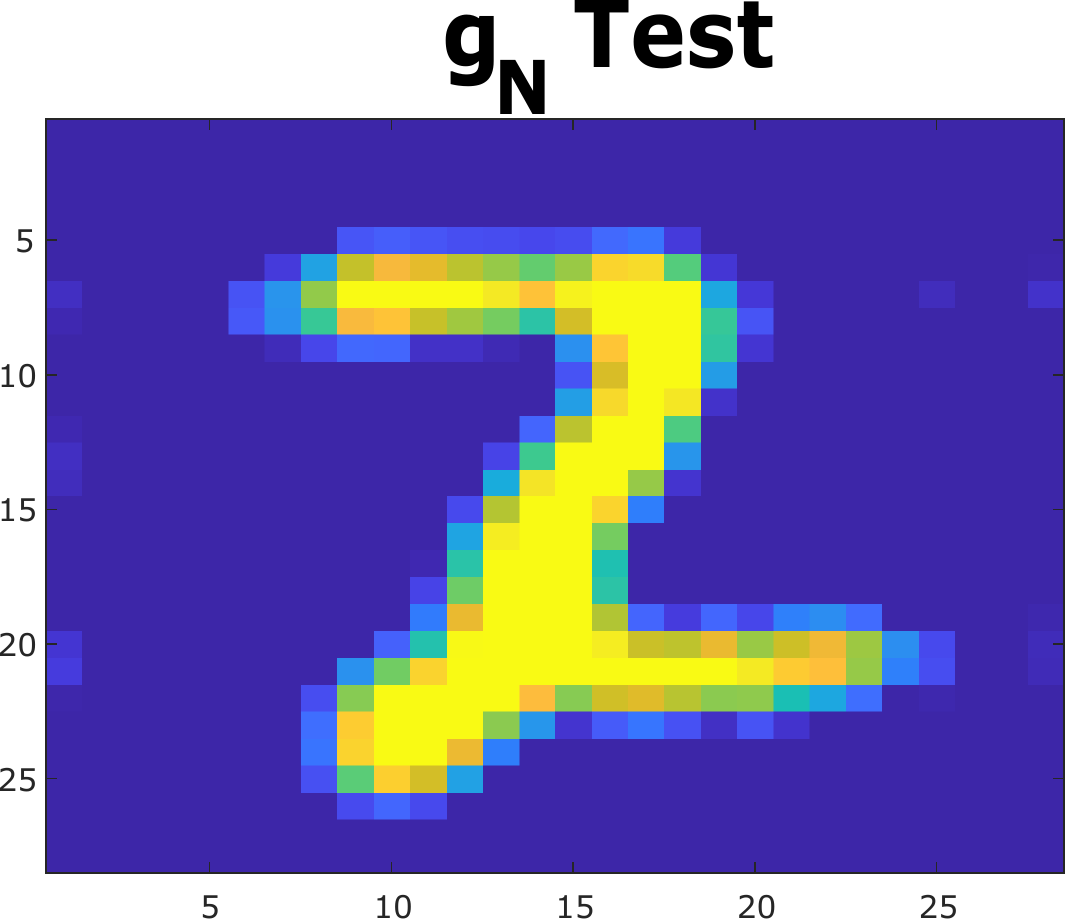}
\hspace{0.8cm}\includegraphics[width=0.15\textwidth, height=0.135\textwidth]{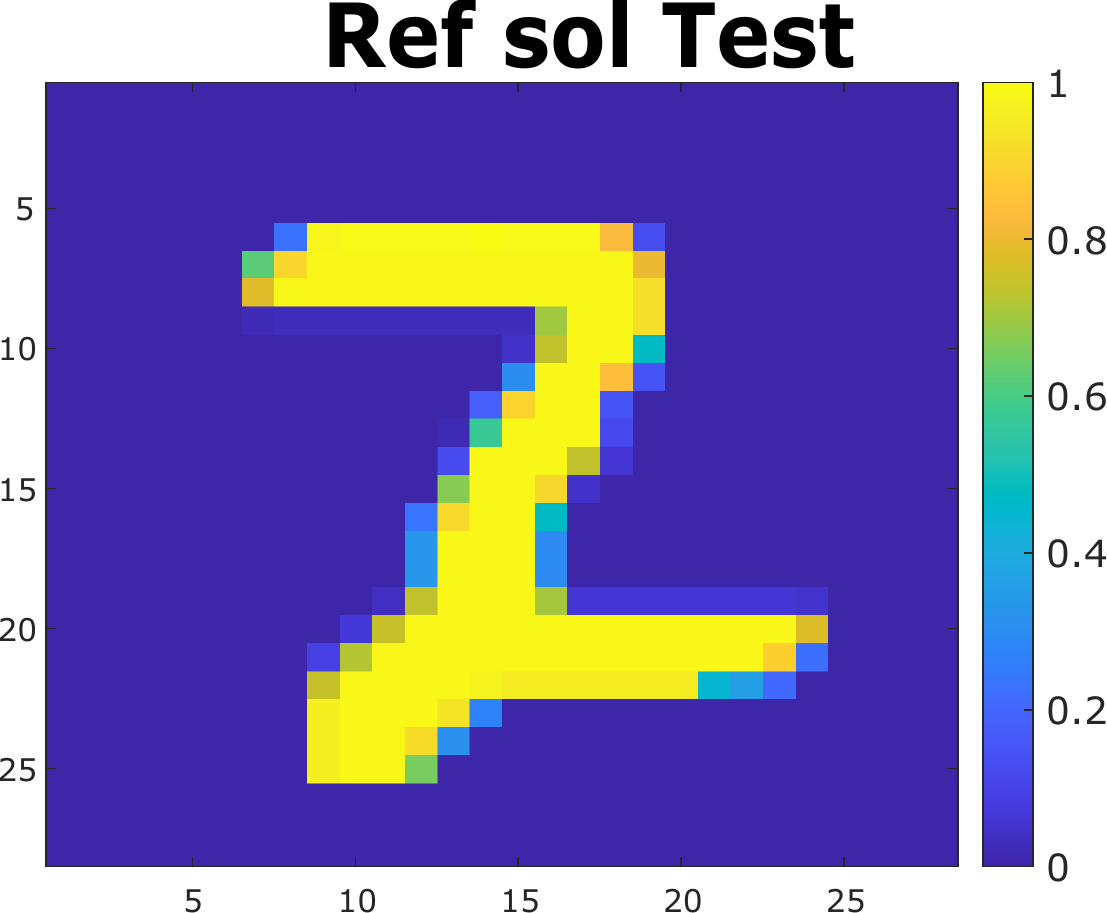}
\caption{\label{f:deblurr_all}Reconstructions of representative samples from training (\textit{row} $1$), validation (\textit{row} $2$) and testing (\textit{row} $3$) data showing successful image delurring performed using OCTANE algorithm with $N=\{6,12,20\}$ layers (\textit{Columns} $2$-$4$), blurred input data as the initial condition (\textit{column} $1$) and reference solution (\textit{column} $5$).}
\end{center}
\end{figure}

We perform runs for $N = \{4, 6, 10, 12, 20, 30\}$ layers. \Cref{f:deblurr_all} presents representative reconstructions for $N = \{6, 12, 20\}$: input $\hat{x}$ (\textit{column} 1), deblurred outputs $\hat{\sigma}(g_N)$ (\textit{columns} 2–4), and reference $x$ (\textit{column} 5). The testing reconstruction errors are $\alpha_{\text{test}} = \{5.9\text{e-}3,\,5.6\text{e-}3,\,5.3\text{e-}3\}$, with corresponding \texttt{PSNR} values $\{19.3,\,19.6,\,19.8\}$ and \texttt{SSIM} values $\{0.87,\,0.893,\,0.891\}$.

\begin{figure}[h!]
\begin{center}
\includegraphics[width=0.3\textwidth, height=0.3\textwidth]{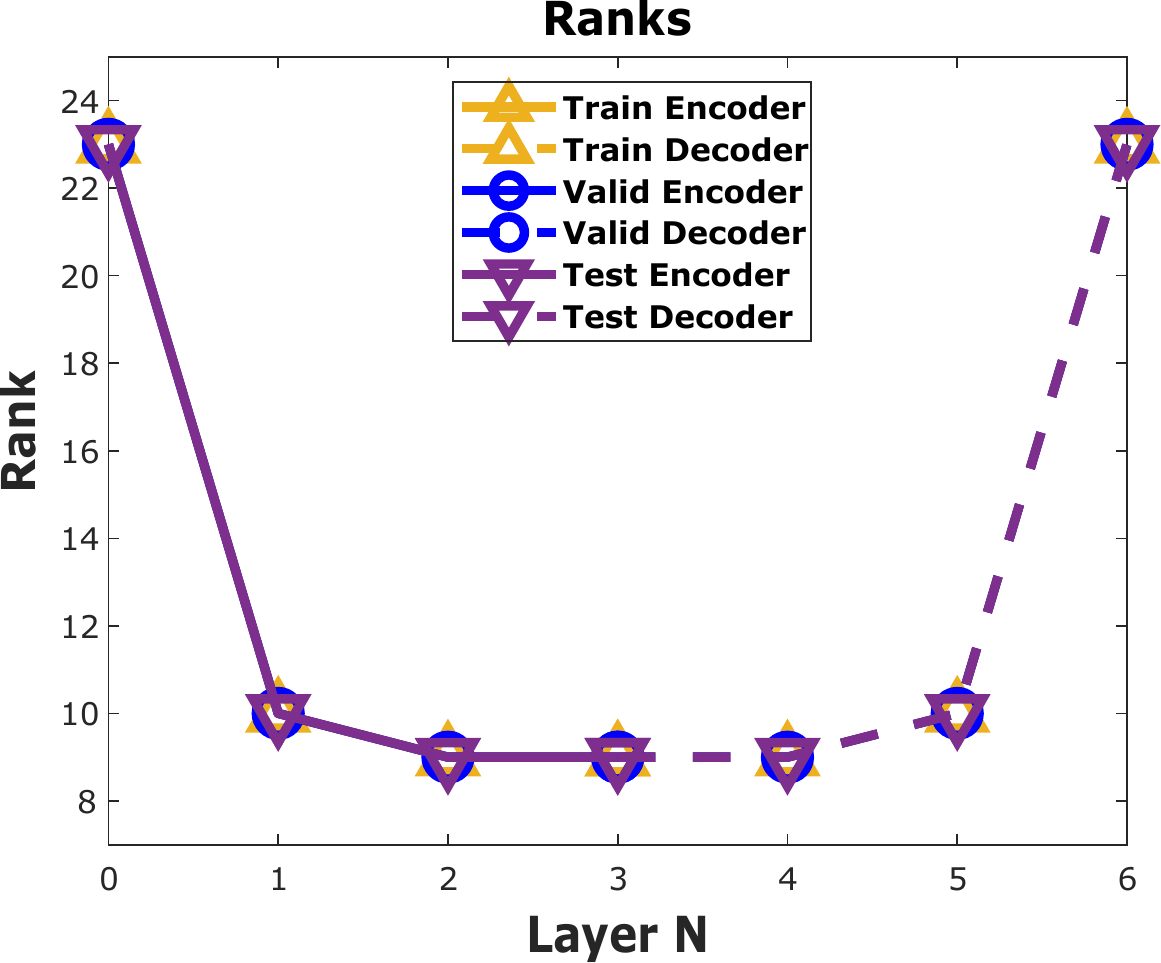} 
\hspace{0.5cm} \includegraphics[width=0.3\textwidth, height=0.3\textwidth]{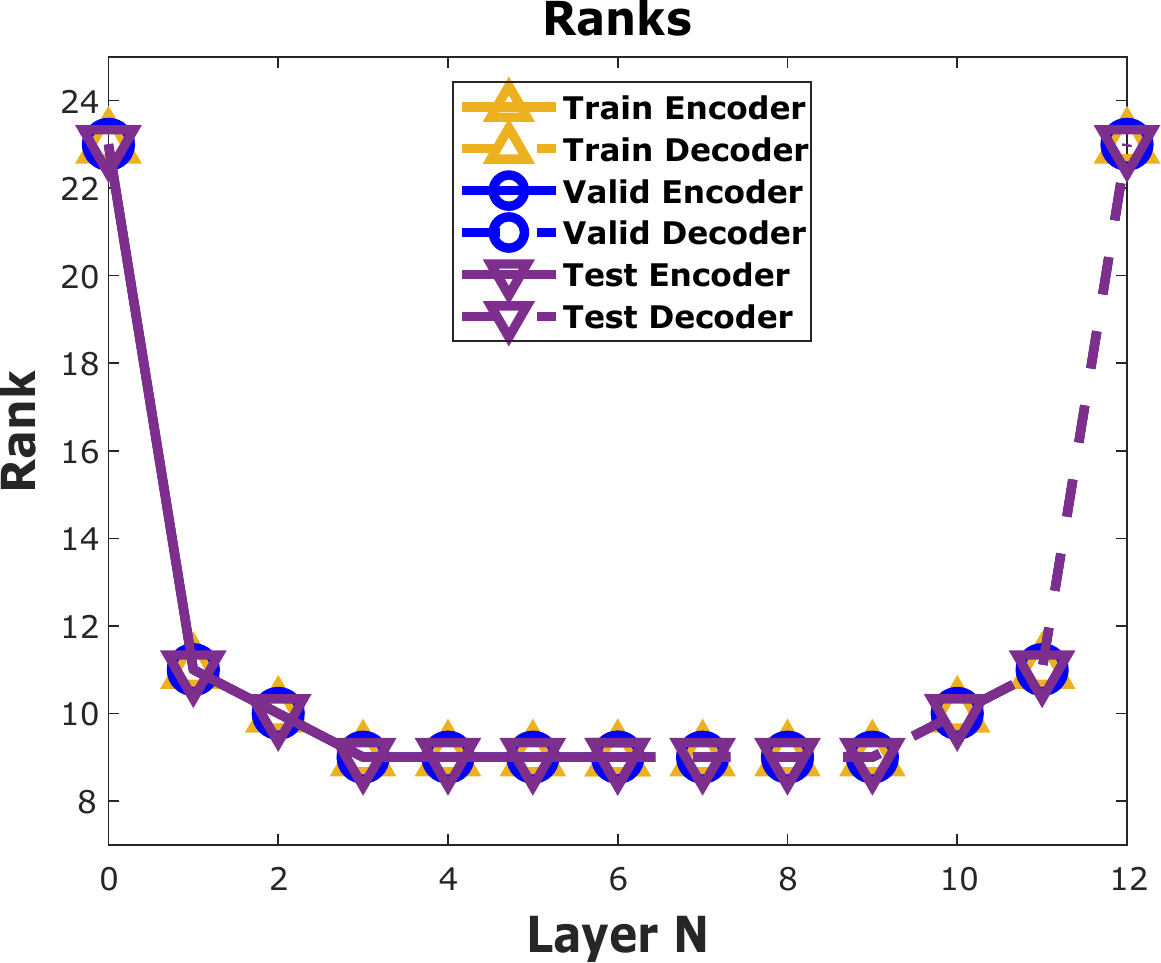}
\hspace{0.5cm} \includegraphics[width=0.3\textwidth, height=0.3\textwidth]{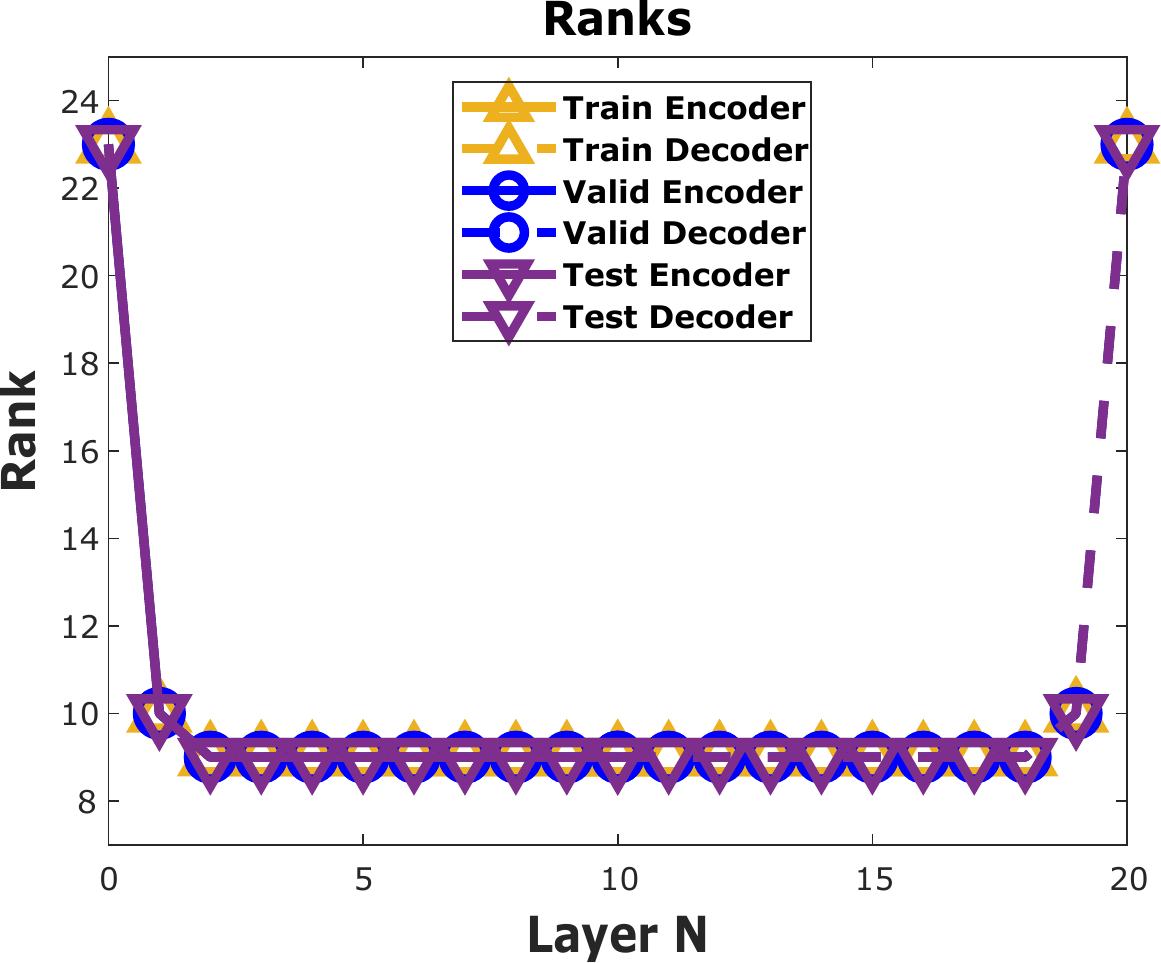} 
\caption{\label{f:deblurr_ranks}Rank distributions obtained via OCTANE algorithm for layers $N=6$, $N=12$, and $N=20$ in the image deblurring task. Solid lines are encoder ranks and dotted lines are decoder ranks. Note that compression in the network is in the context of rank reduction.}
\end{center}
\end{figure}

\begin{figure}[h!]
\begin{center}
\includegraphics[width=0.3\textwidth, height=0.3\textwidth]{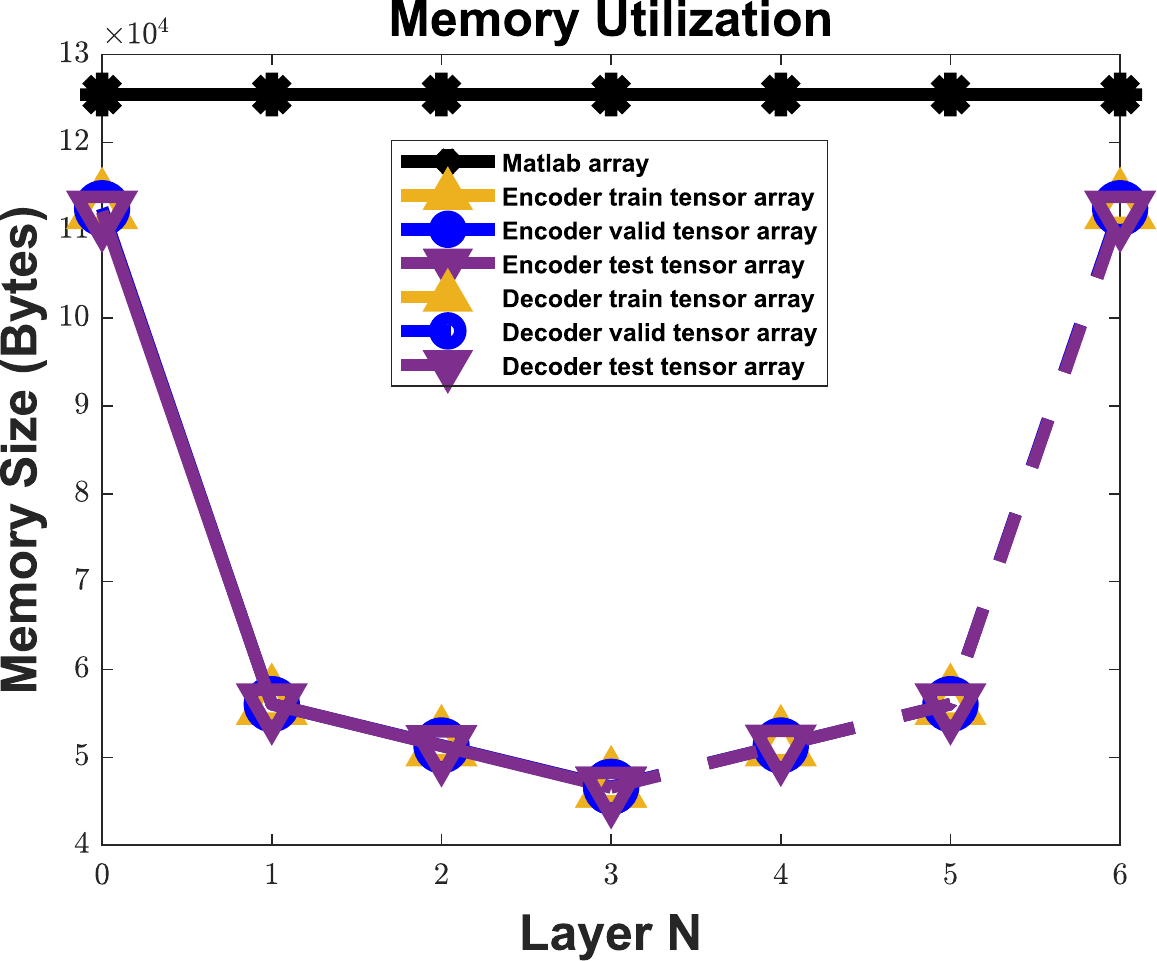} 
\hspace{0.5cm} \includegraphics[width=0.3\textwidth, height=0.3\textwidth]{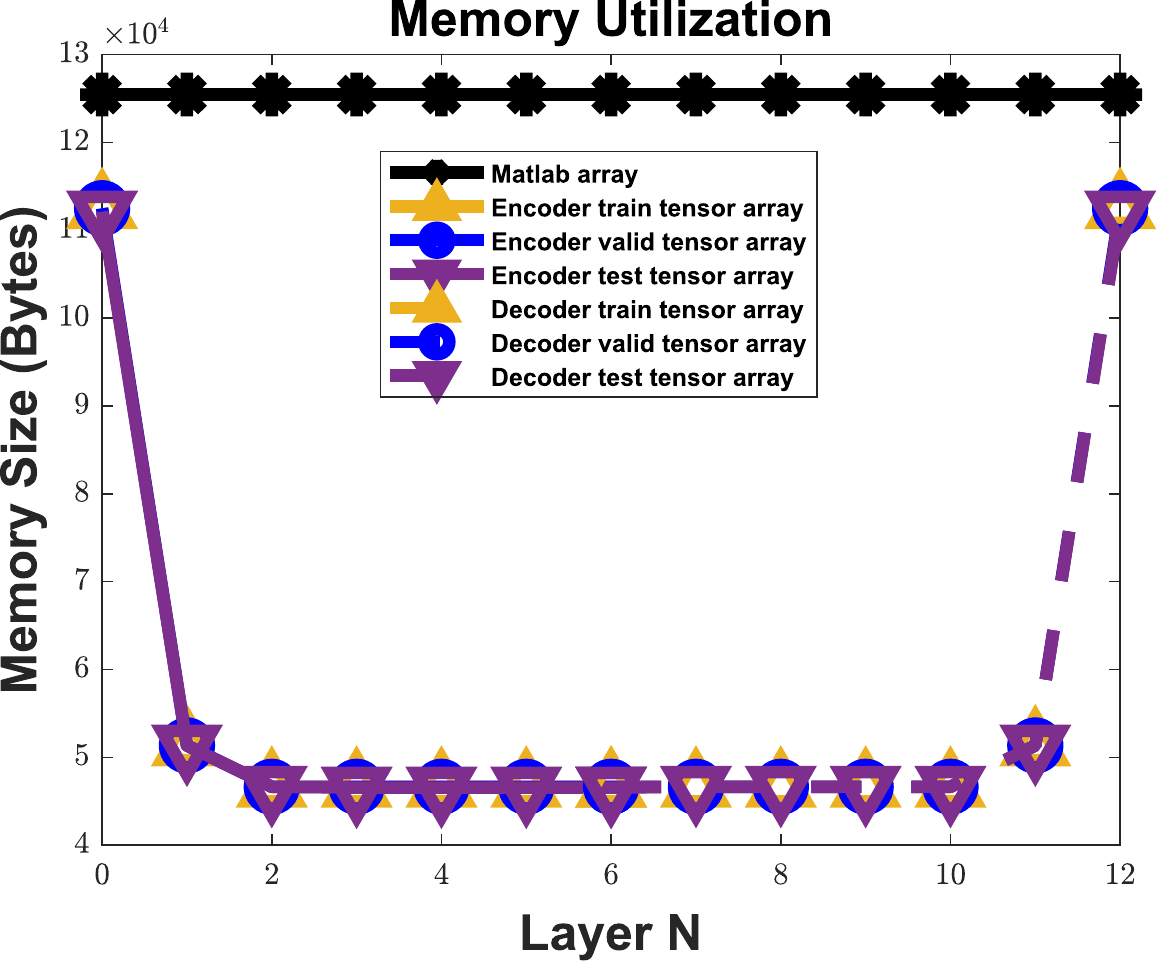}
\hspace{0.5cm} \includegraphics[width=0.3\textwidth, height=0.3\textwidth]{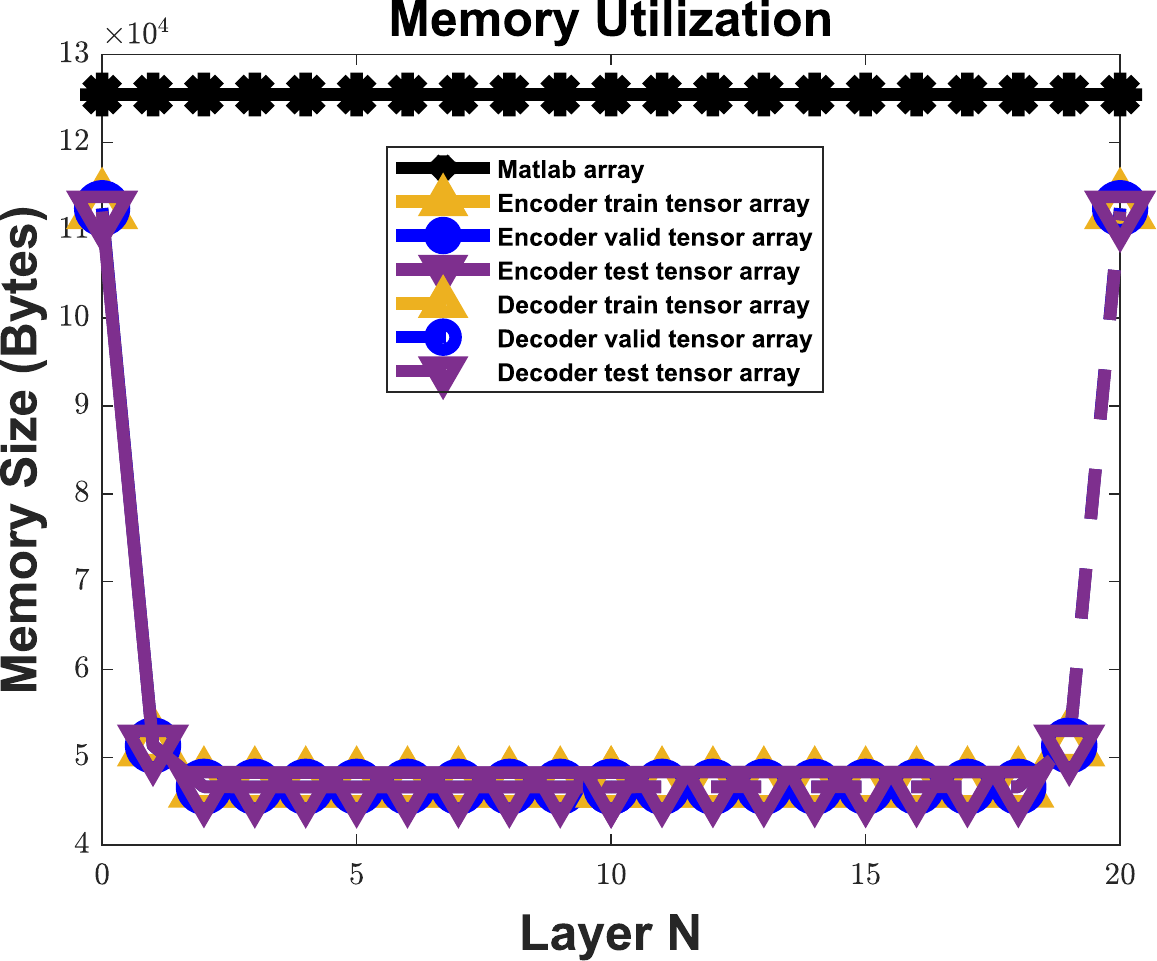} 
\caption{\label{f:deblurr_mem}Rank distributions obtained via OCTANE algorithm for $N=6$, $N=12$, and $N=20$ in the image deblurring task. Solid lines are encoder ranks and dotted lines are decoder ranks. Owing to rank-reduction, solving differential equation on the low-rank manifold saves significant memory.}
\end{center}
\end{figure}

\Cref{f:deblurr_ranks} displays the rank distributions during forward propagation for the deblurring task with $N = \{6, 12, 20\}$. Compared to denoising, the ranks drop even further (as low as $9$), indicating that OCTANE adapts its architecture based on the task. \Cref{f:deblurr_mem} compares the memory used by state variables (orange: training, blue: validation, purple: testing) in tensor format against a standard MATLAB array (black). Due to the lower-rank latent structure, memory savings are substantial—averaging $46.74\%$, $53\%$, and $57.46\%$ for $N = 6$, $12$, and $20$ layers, respectively—significantly outperforming the denoising case.

\begin{figure}[h!]
\begin{center}
\includegraphics[width=0.3\textwidth, height=0.3\textwidth]{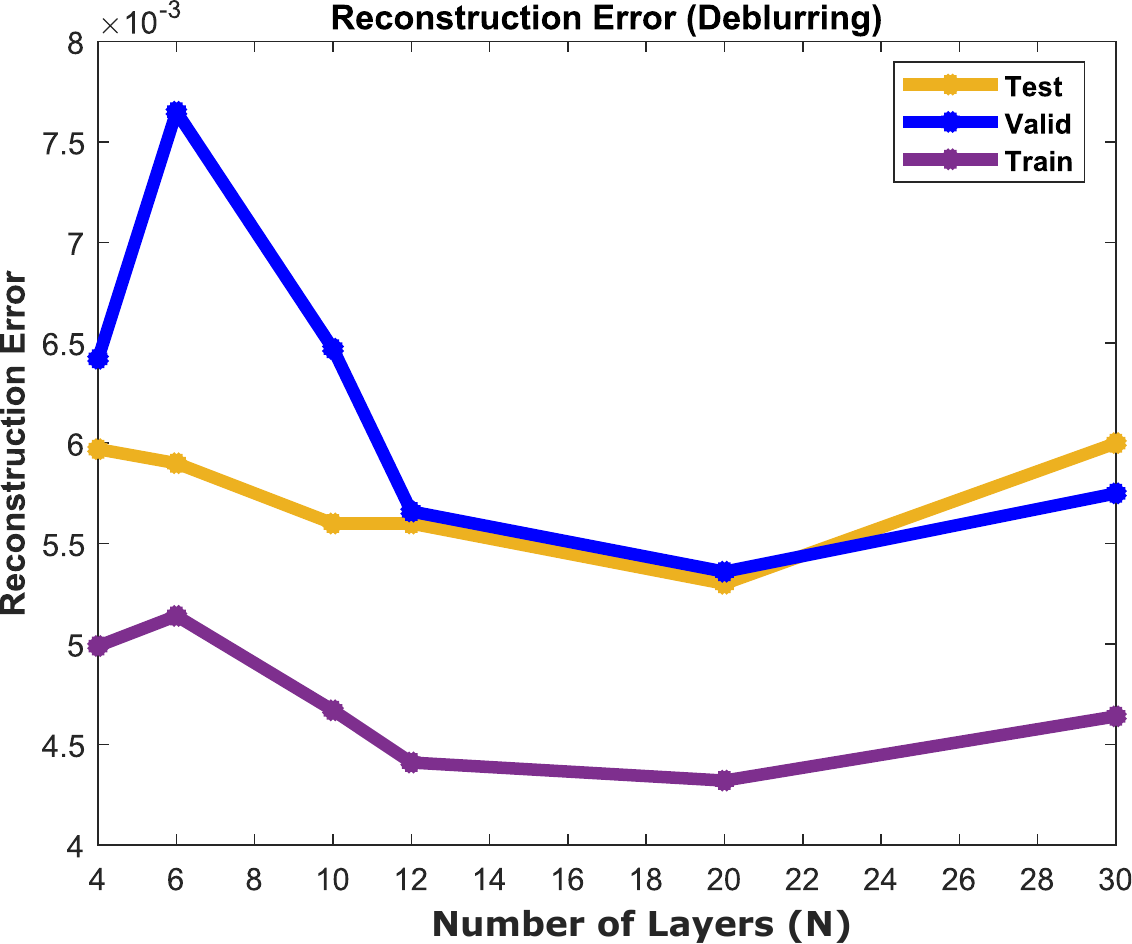} 
\includegraphics[width=0.3\textwidth, height=0.3\textwidth]{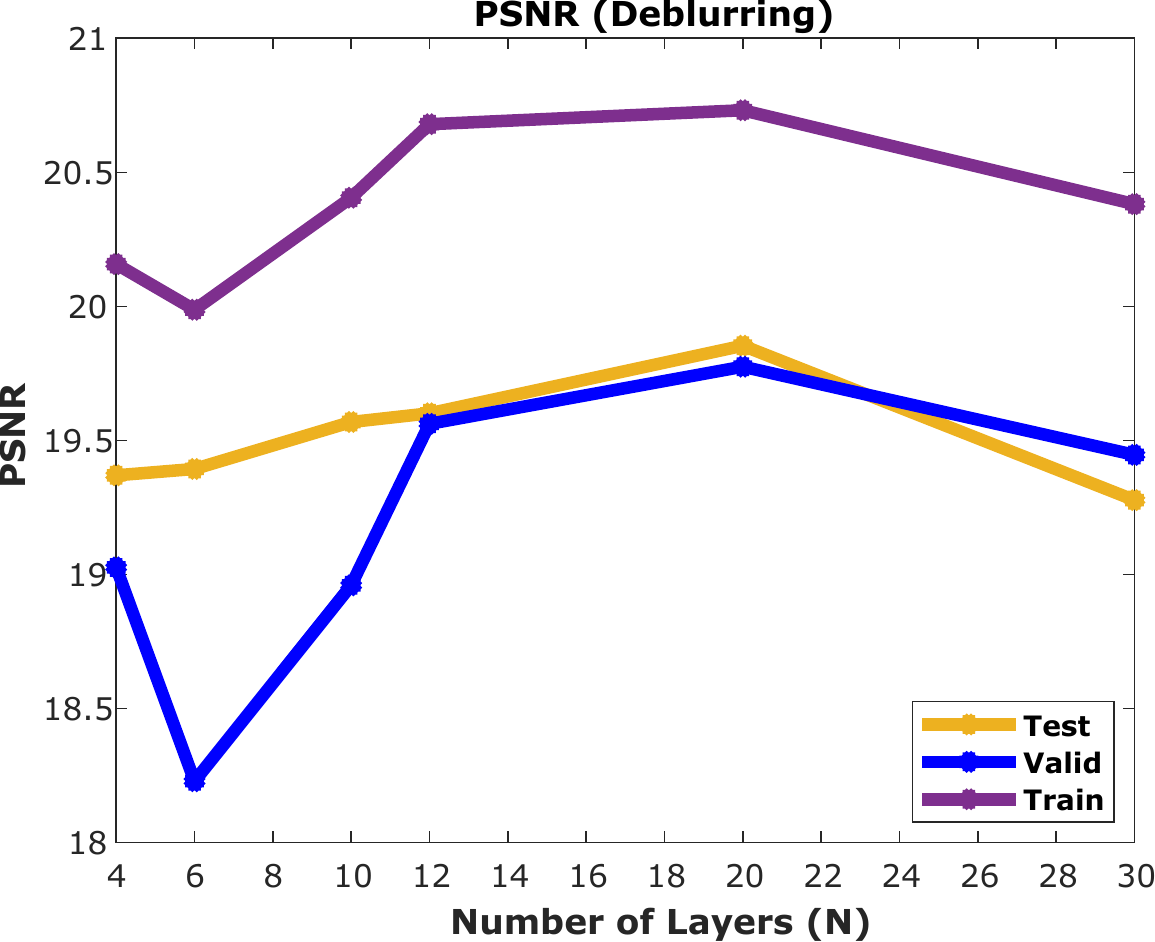}
\includegraphics[width=0.3\textwidth, height=0.3\textwidth]{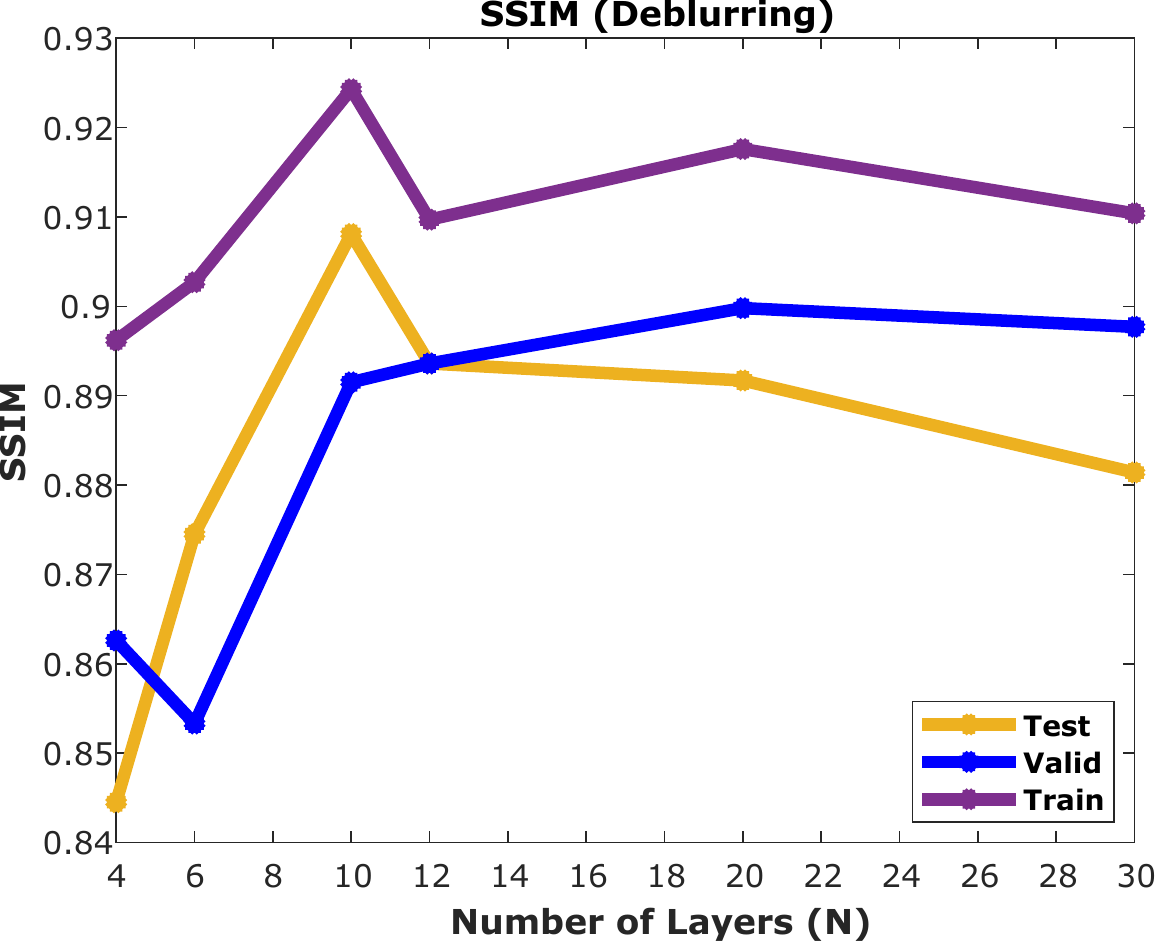}
\caption{\label{f:deblurr_results_per_N}Plots of the Reconstruction Errors, \texttt{PSNR}, and \texttt{SSIM} for various number of layers in the deblurring experiments. The best results are achieved at $N=20$ layers. The testing data generally follows the same trend as training data. Overall, there is minimal improvement in the metrics around $N=10$ to $N=20$ (relatively expensive), therefore lesser layers are sufficient due to memory saving benefit.}
\end{center}
\end{figure}

To assess the impact of layer count in deblurring, \cref{f:deblurr_results_per_N} presents reconstruction error \cref{recon_err}, \texttt{PSNR}, and \texttt{SSIM} for $N = \{4, 6, 10, 12, 20, 30\}$. The results confirm that OCTANE effectively deblurs images. While $N=20$ yields the lowest training error, improvements beyond $N=10$ are marginal across all metrics, suggesting that smaller architectures suffice. Testing curves mirror training trends, and validation results stabilize as $N$ increases. See \cref{Sec:TauN} for further discussion. Each experiment completes within $40$ minutes.

\subsection{The interplay between hyper-parameters step-size $\tau$ and layer-count $N$.} \label{Sec:TauN}

In \cref{Sec:Denoising,Sec:Deblurring}, we demonstrated that OCTANE can effectively denoise and deblur images. However, its performance and computational cost are sensitive to the initial network selection—specifically, in terms of the number of layers $N$, terminal time $T$, and step size $\tau = T/N$ in the rank-adaptive Euler scheme. While a finer $\tau$ generally improves accuracy, the interplay between $N$ and $T$ introduces two degrees of freedom. Notably, choosing $M_s = \tau^{-1}$ and $M_r = \tau^{-2}$ decouples the rank reduction tolerance from $\tau$, enforcing a Lipschitz-type condition. Yet, the hyperparameters $N$ and $\tau$ must still be selected beforehand. In this section, we investigate this trade-off through numerical experiments and propose a heuristic strategy for choosing $N$ and $\tau$ to ensure consistently good performance.

\begin{figure}[htb!]
\begin{center}
\includegraphics[width=0.3\textwidth, height=0.3\textwidth]{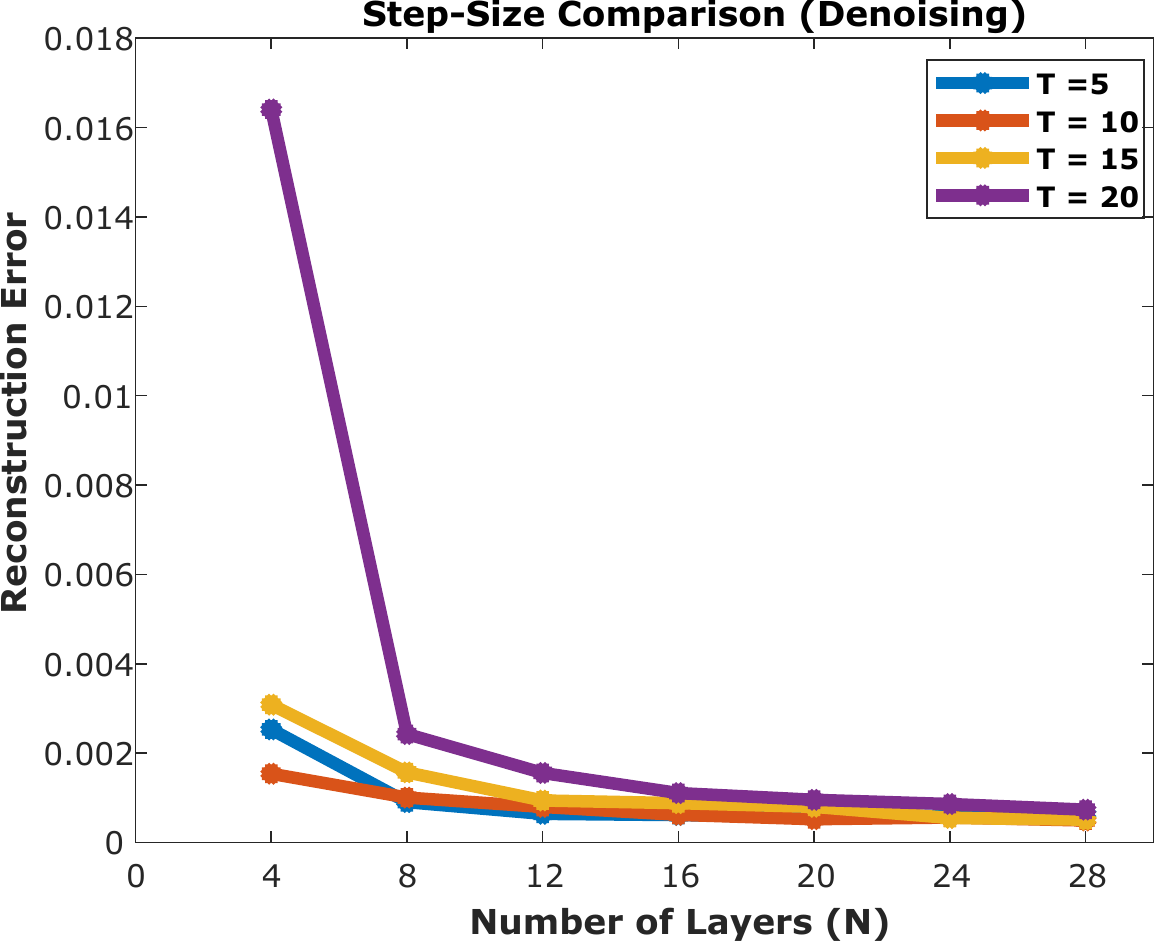} 
\includegraphics[width=0.3\textwidth, height=0.3\textwidth]{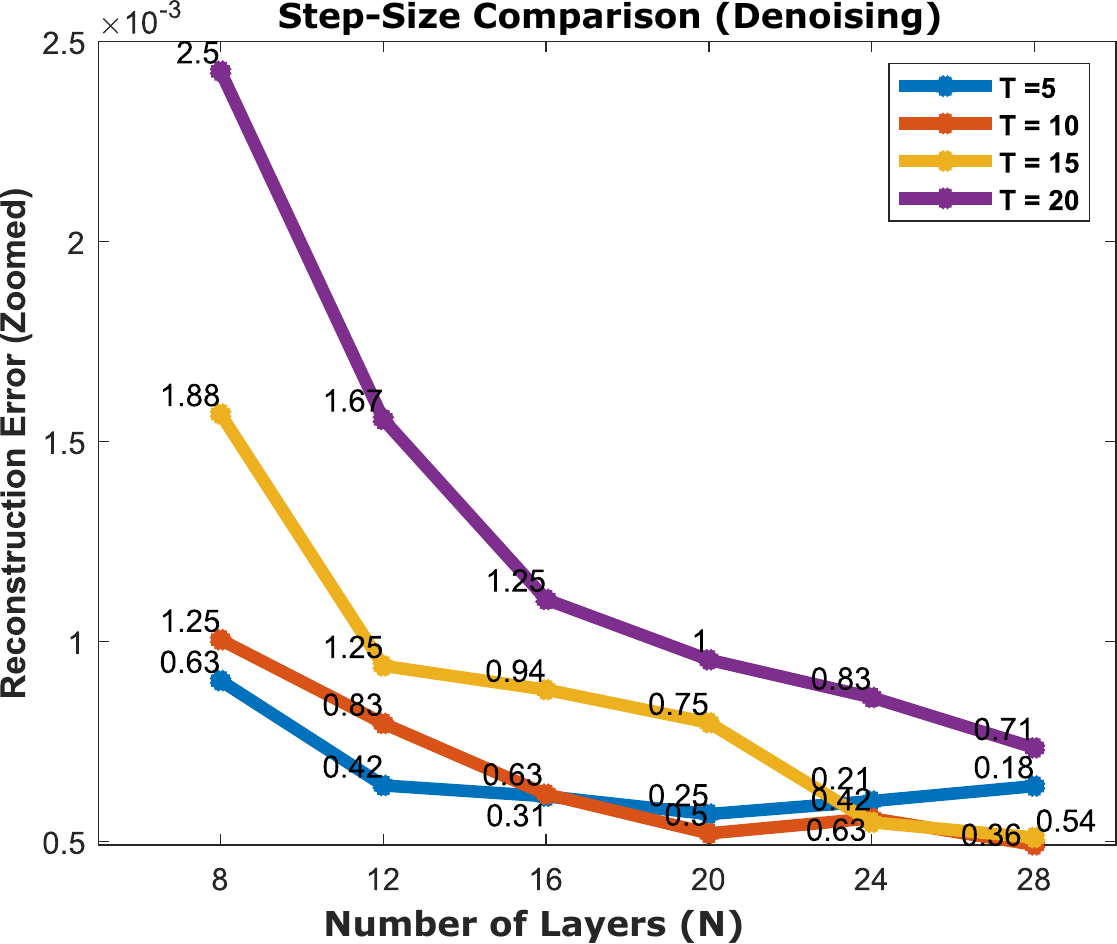}
\includegraphics[width=0.3\textwidth, height=0.3\textwidth]{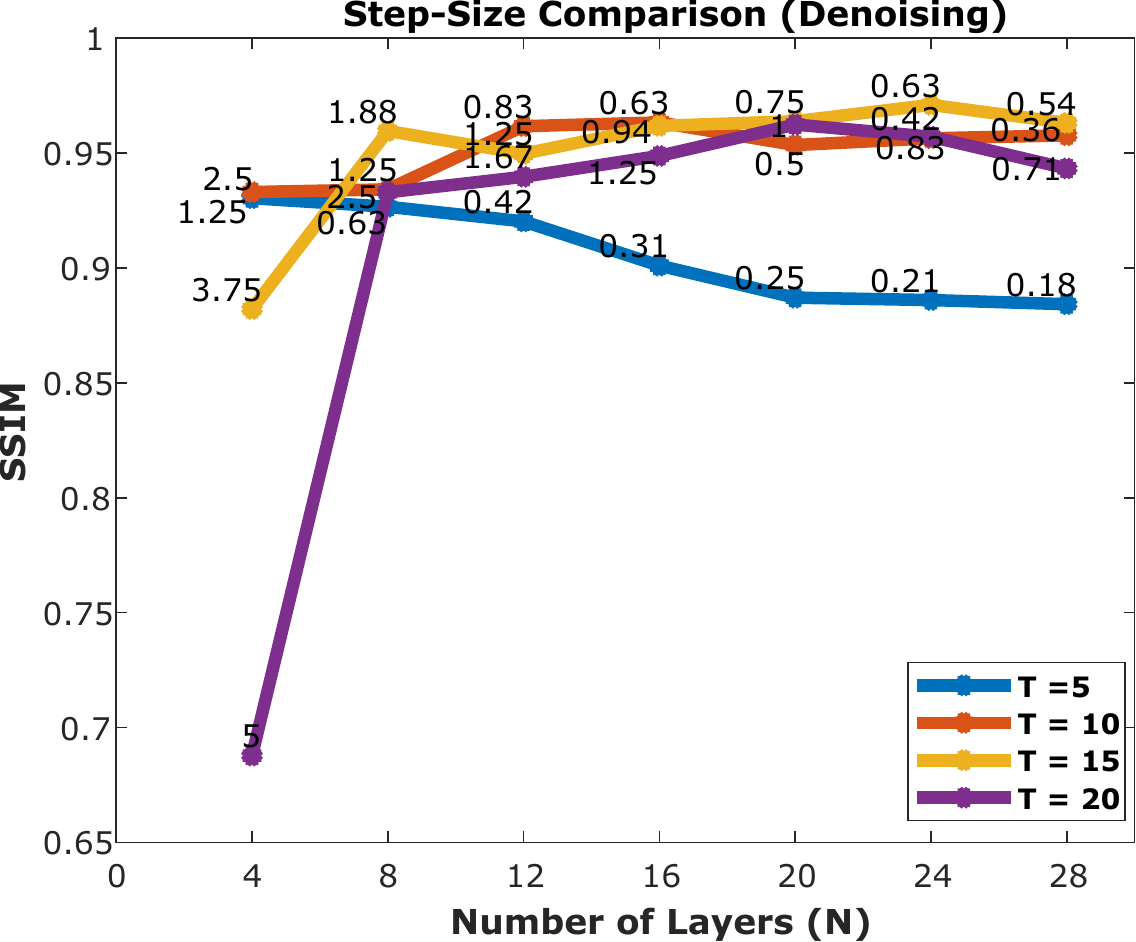}
\caption{\label{f:comparison_denoise}Comparison of reconstruction errors (zoomed-out (\textit{left}) and zoomed-in (\textit{center})) and \texttt{SSIM} (\textit{right}) against the number of layers $N$ for various values of final time $T$, with corresponding $\tau$ values mentioned on the plots for each image denoising experiment. These plots demonstrate that $0.3 \leq \tau \leq 1.3$ typically gives the best results, even for small $N$.}
\end{center}
\end{figure}

\begin{figure}[htb!]
\begin{center}
\includegraphics[width=0.3\textwidth, height=0.3\textwidth]{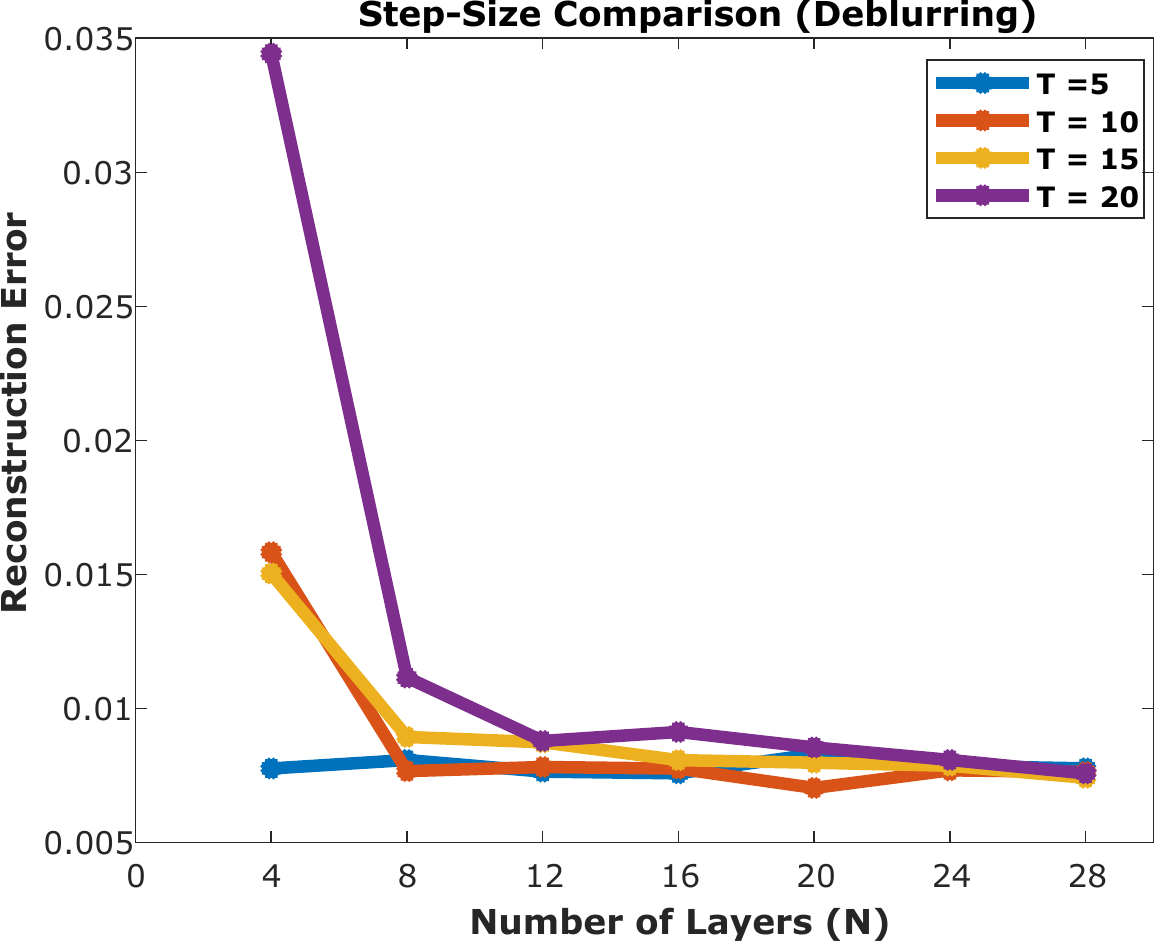} 
\includegraphics[width=0.3\textwidth, height=0.3\textwidth]{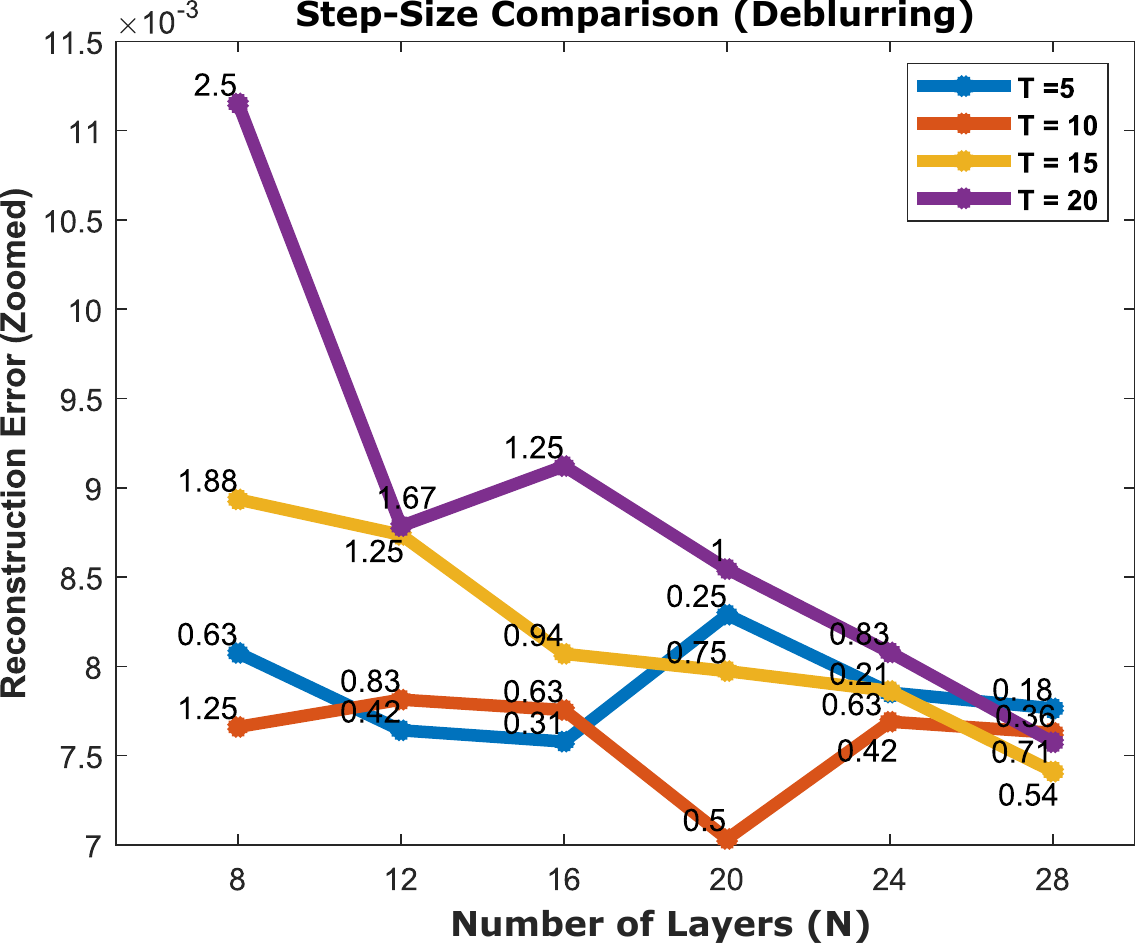}
\includegraphics[width=0.3\textwidth, height=0.3\textwidth]{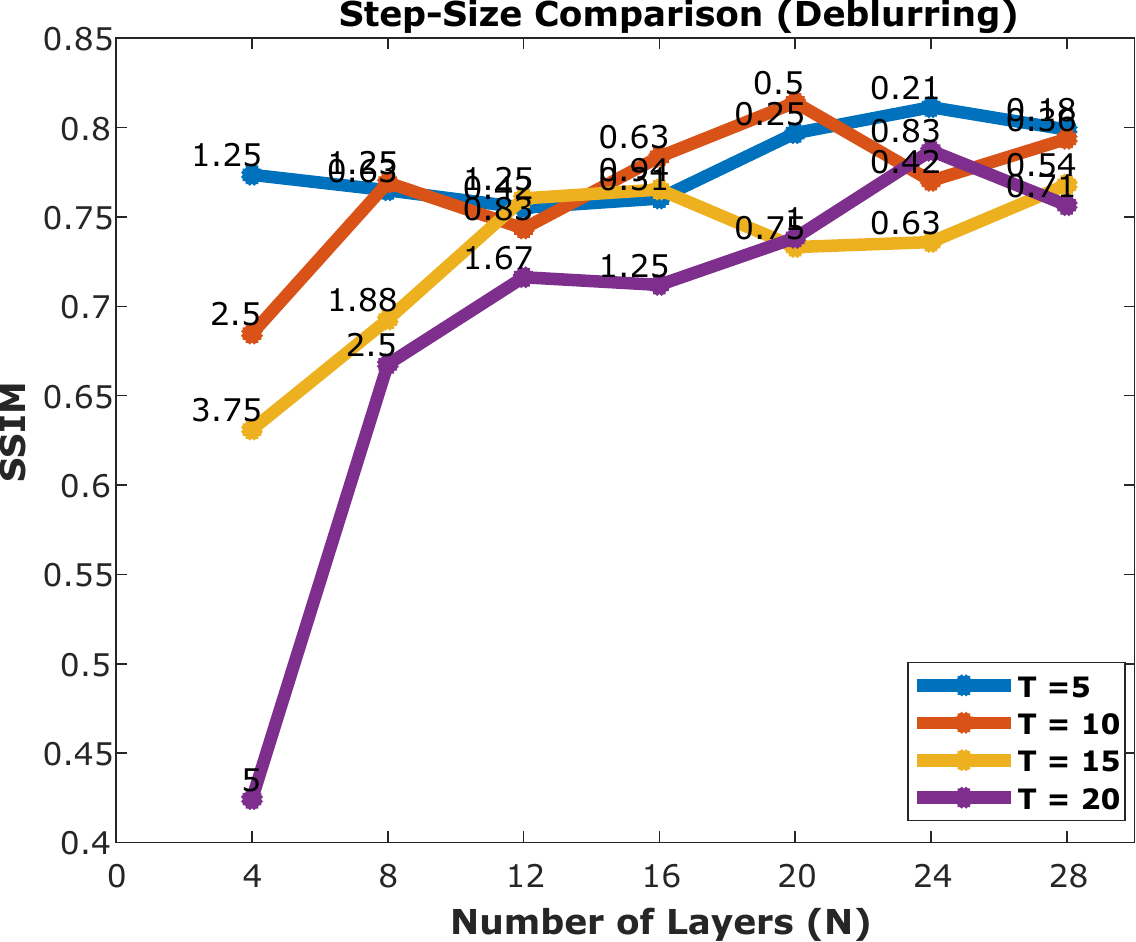}
\caption{\label{f:comparison_deblurr}Comparison of reconstruction errors (zoomed-out (\textit{left}) and zoomed-in (\textit{center})) and \texttt{SSIM} (\textit{right}) against the number of layers $N$ for various values of final time $T$, with corresponding $\tau$ values mentioned on the plots for each image deblurring experiment. These plots demonstrate that $0.3 \leq \tau \leq 1.3$ typically gives the best results, even for small $N$.}
\end{center}
\end{figure}

We now study the influence of layer count $N$ and terminal time $T$ on network performance. Fixing $T \in \{5, 10, 15, 20\}$, we vary $N \in \{4, 8, 12, \dots, 28\}$, setting $\tau = T/N$ accordingly. Each $(T, N)$ pair defines an independent run of the imaging task (denoising or deblurring) on MNIST digit $2$, using the same setup and hyperparameters as in \cref{Sec:Denoising,Sec:Deblurring,t:exp_config}. We collect training reconstruction errors and \texttt{SSIM} indices for each experiment and plot them in \cref{f:comparison_denoise} (denoising) and \cref{f:comparison_deblurr} (deblurring), with $T$ color-coded and $\tau$ annotated. The left panels show the full error trends, center panels zoom into key regions, and right panels display \texttt{SSIM}. These plots help identify practical $(\tau, T, N)$ combinations that yield robust performance across tasks.

The reconstruction error plots (\textit{left/center} in \cref{f:comparison_denoise,f:comparison_deblurr}) generally show that increasing $N$ reduces error, as expected from ODE discretization theory, though not strictly monotonically.

Focusing on denoising (\cref{f:comparison_denoise}), most good reconstructions fall in the bottom-right of the center plot. Discarding $T=20$ due to poor performance and restricting to $\alpha < 1\text{e--}3$, we obtain the feasible set,
\[
\tau_a := \{\tau \mid \alpha_{T = \{5,10,15\}} \leq 1\text{e--}3\} = (0.1,1.3).
\]
Meanwhile, for \texttt{SSIM} $\geq 0.96$, we discard $T=5$ (due to degradation at larger $N$), yielding,
\[
\tau_b := \{\tau \mid \texttt{SSIM}_{T = \{10,15,20\}} \geq 0.96\} = (0.3,1.88).
\]
Intersecting both gives a good denoising range:
\[
\tau_{\text{denoising}} = \tau_a \cap \tau_b = (0.3,1.3),\quad T = \{10,15\}.
\]

For deblurring (\cref{f:comparison_deblurr}), restricting $\alpha < 8.5\text{e--}3$ and \texttt{SSIM} $\geq 0.75$ (discarding $T=20$), gives,
\[
\tau_c := (0.1,1.3),\qquad \tau_d := (0.1,1.3) \Rightarrow \tau_{\text{deblurring}} = \tau_c \cap \tau_d = (0.1,1.3),\quad T = \{5,10,15\}.
\]

Combining both tasks:
\begin{equation}\label{tau_prop}
\tau_{\text{proposed}} = \tau_{\text{denoising}} \cap \tau_{\text{deblurring}} = (0.3,1.3),\quad T = \{10,15\}.
\end{equation}

With $T = \tau N$, this corresponds to
\begin{equation}\label{N_prop}
N_{\text{proposed}} = \{2k \mid k = 6,\dots,17\},
\end{equation}

Based on performance and cost, we recommend selecting $(N, \tau)$ from this range \cref{tau_prop}-\cref{N_prop}, favoring smaller $N$ when possible.

\begin{figure}[h!]
\begin{center}
\hspace{-0.5cm}\fbox{\text{\textbf{Initial Cond.}}} \hspace{1cm}\fbox{Subjacent} \hspace{1.3cm}\fbox{Within} \hspace{1.3cm} \fbox{Superjacent} \hspace{0.9cm}\fbox{\text{\textbf{Ref. Sol.}}}\\~\\
\includegraphics[width=0.14\textwidth, height=0.14\textwidth]{figures/denoising/f0_train.pdf} 
\hspace{0.8cm}\includegraphics[width=0.14\textwidth, height=0.14\textwidth]{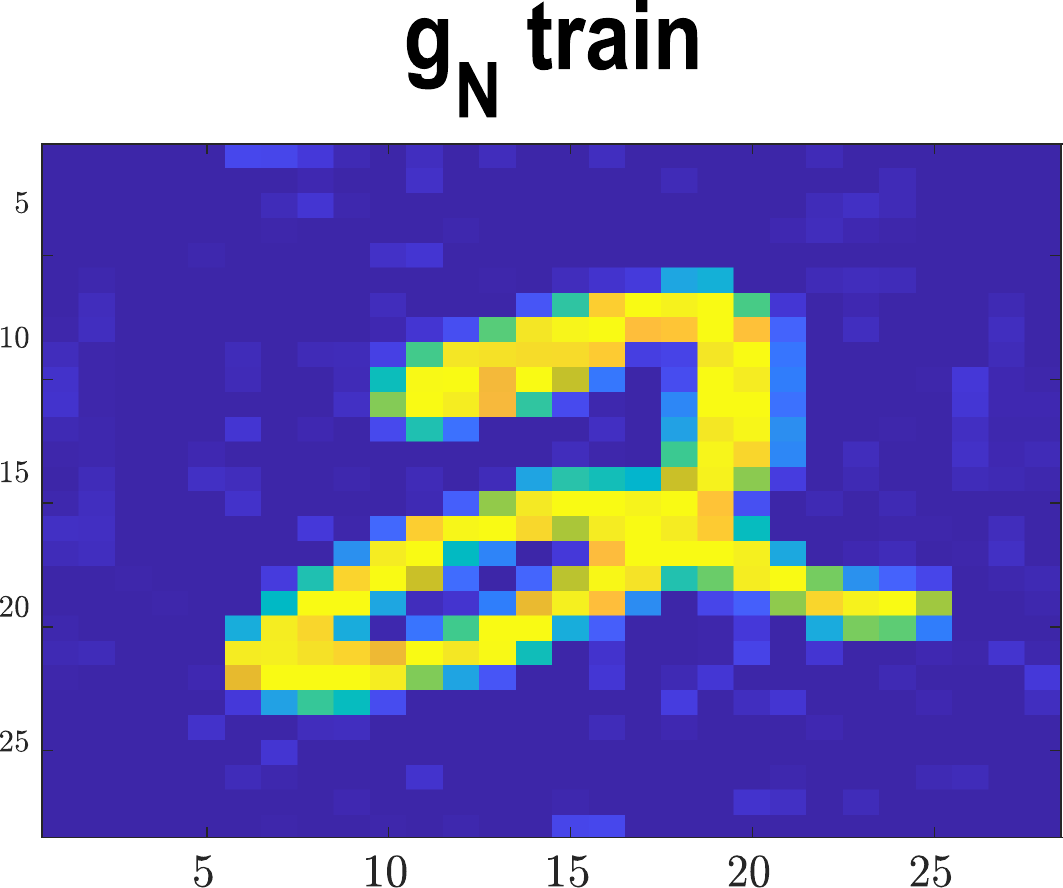}
\hspace{0.8cm}\includegraphics[width=0.14\textwidth, height=0.14\textwidth]{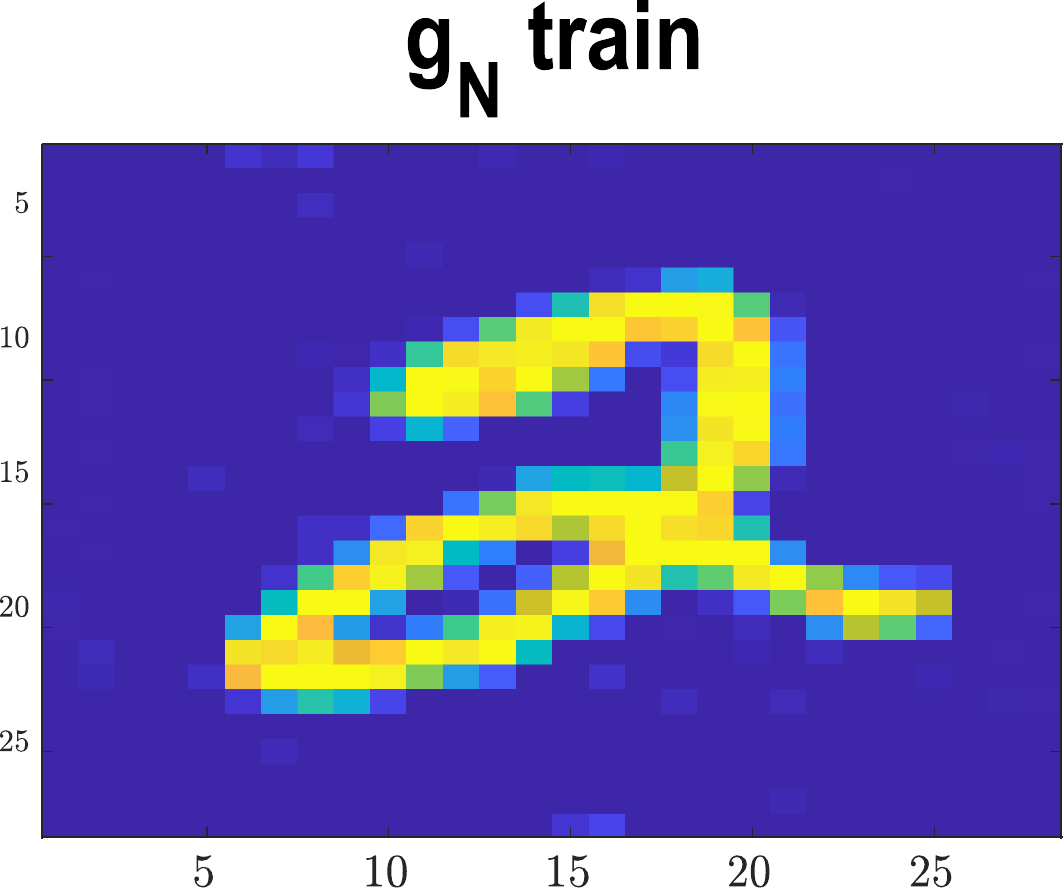}
\hspace{0.8cm}\includegraphics[width=0.14\textwidth, height=0.14\textwidth]{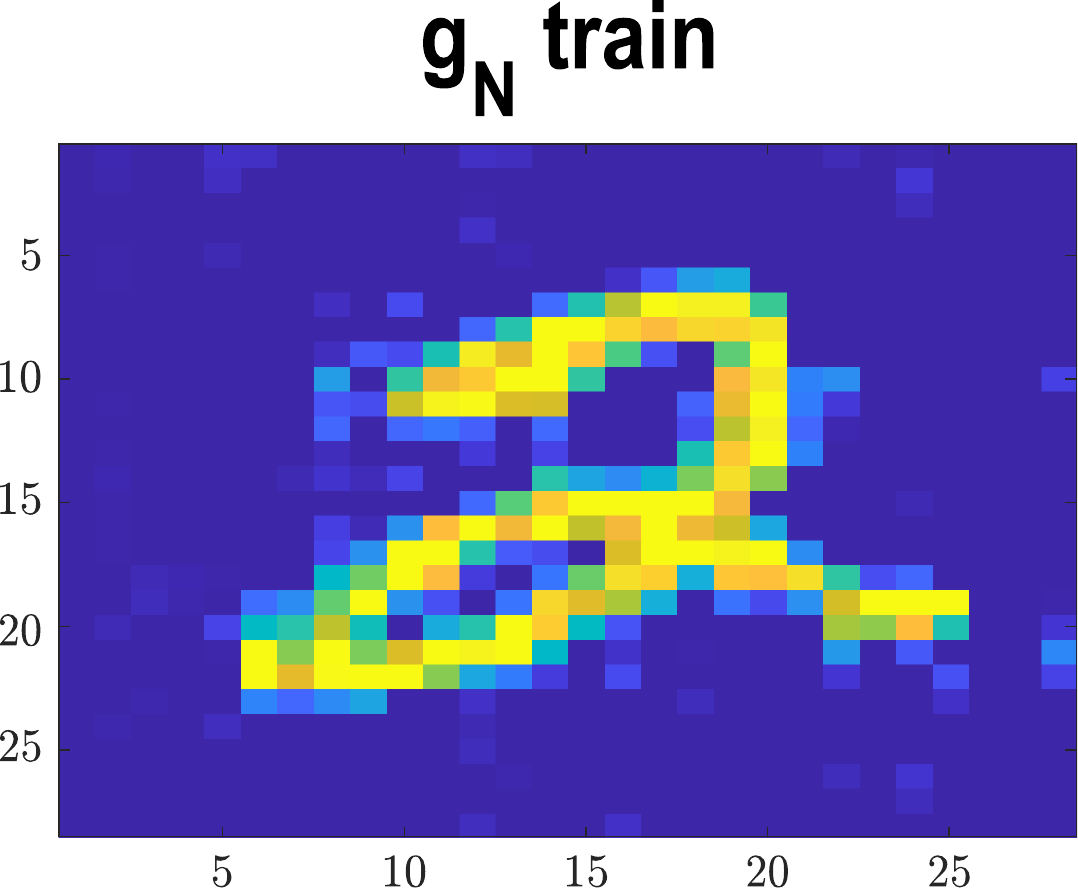}
\hspace{0.8cm}\includegraphics[width=0.15\textwidth, height=0.13\textwidth]{figures/deblurring/ref_sol_train.pdf}\\
{\tiny{
\hspace{0cm}$(N,\tau,T) = (28,0.18,5)$ \hspace{0.5cm}$(N,\tau,T) = (16,0.63,10)$ \hspace{0.5cm} $(N,\tau,T) = (4,3.75,20)$\tiny}}\\~\\
\includegraphics[width=0.14\textwidth, height=0.14\textwidth]{figures/deblurring/f0_train.pdf}
\hspace{0.8cm}\includegraphics[width=0.14\textwidth, height=0.145\textwidth]{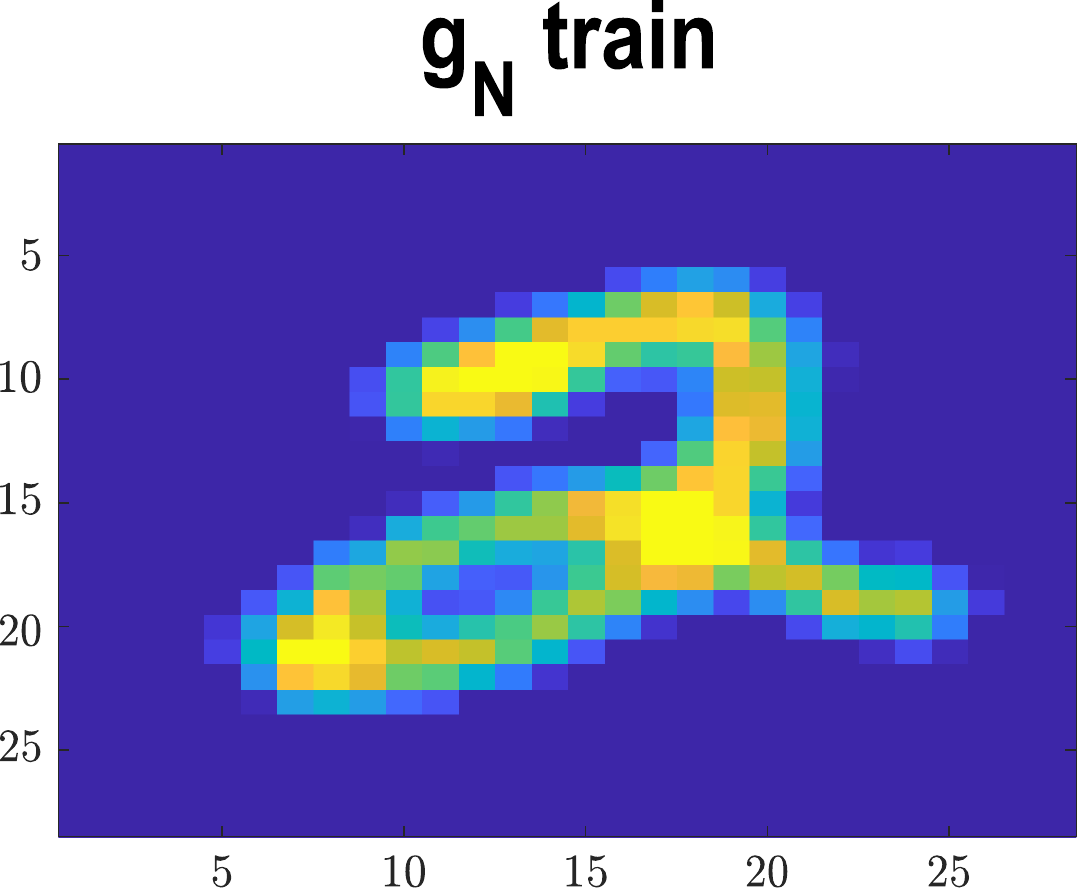}
\hspace{0.8cm}\includegraphics[width=0.14\textwidth, height=0.145\textwidth]{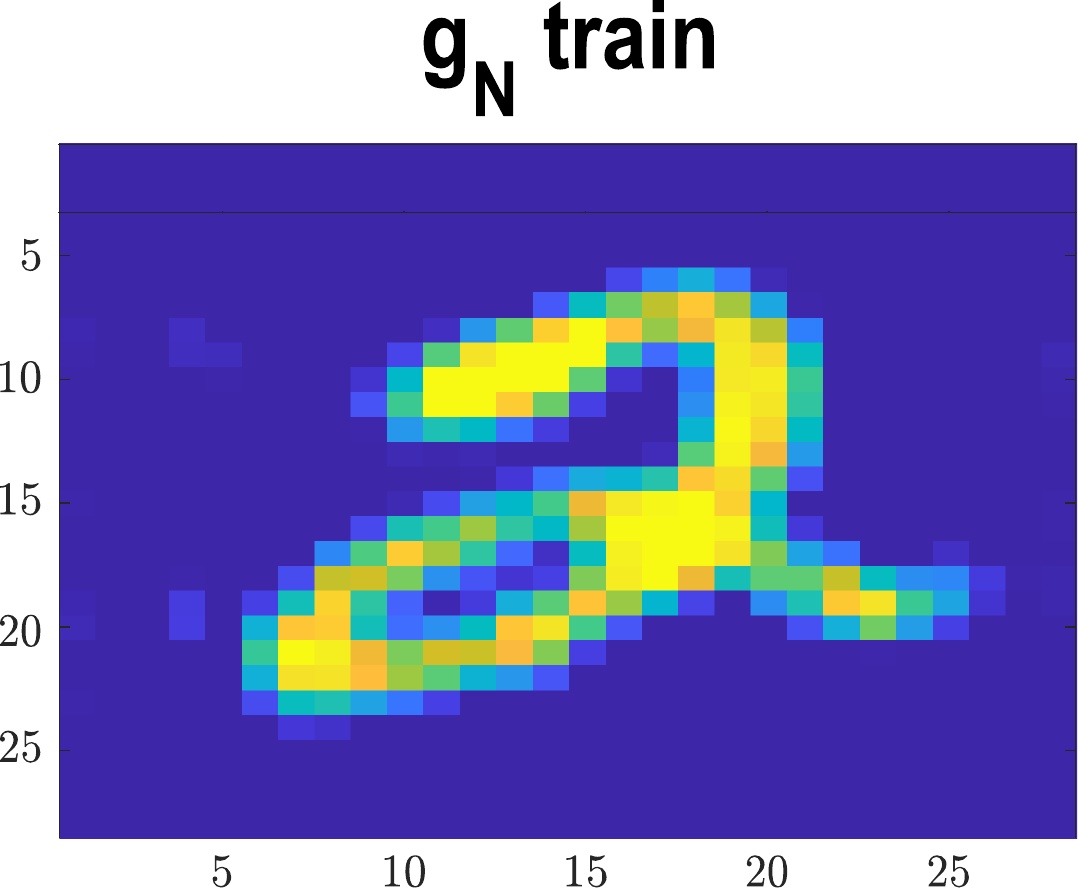}
\hspace{0.8cm}\includegraphics[width=0.14\textwidth, height=0.145\textwidth]{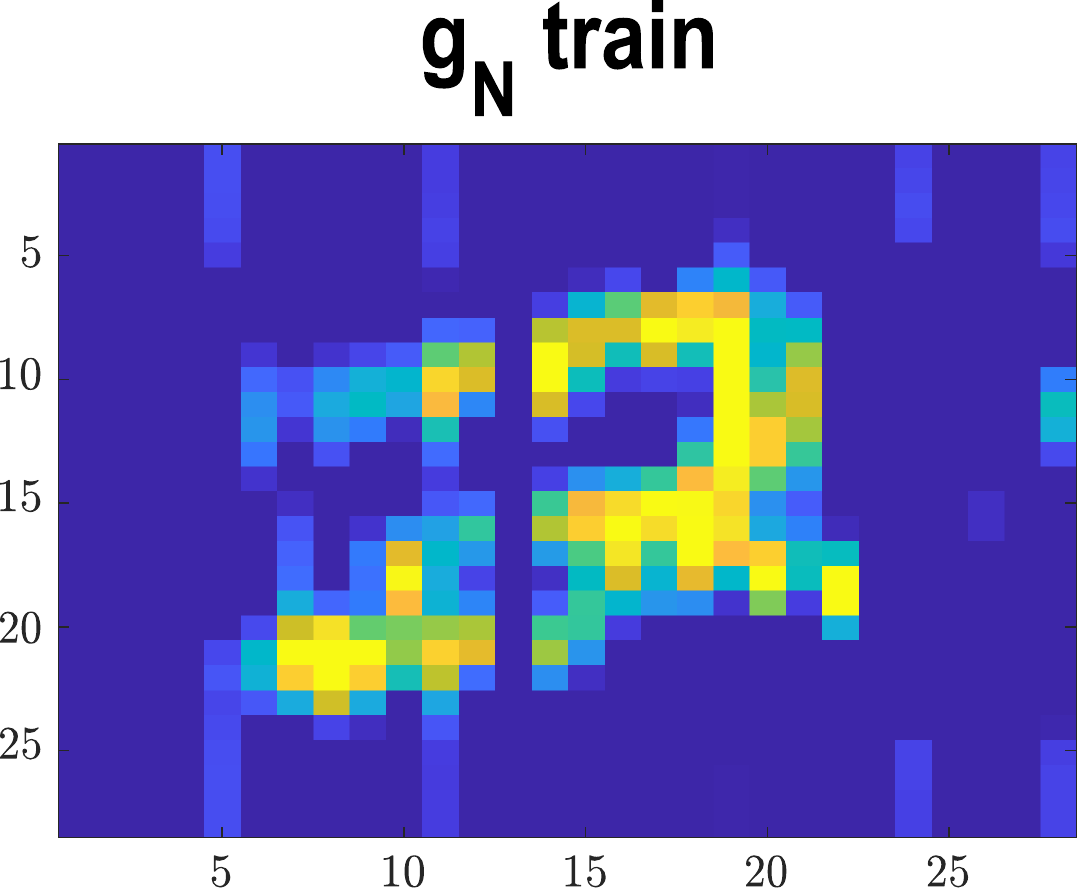}
\hspace{0.8cm}\includegraphics[width=0.15\textwidth, height=0.13\textwidth]{figures/deblurring/ref_sol_train.pdf}
{\tiny{
\hspace{0cm}$(N,\tau,T) = (28,0.18,5)$ \hspace{0.5cm}$(N,\tau,T) = (16,0.63,10)$ \hspace{0.5cm} $(N,\tau,T) = (4,3.75,20)$\tiny}}
\caption{\label{f:GoodBadRecon}Reconstructions of representative samples of digit $2$ for the image denoising (\textit{top row}) and deblurring (\textit{bottom row}) experiments, for $\tau=0.18<0.3$ (\textit{column $2$}), $\tau = 0.63 \in \tau_{proposed}$ (\textit{column $3$}) and $\tau = 3.75 >1.3$ (\textit{column $4$}), along with initial condition (\textit{column $1$}) and reference solutions (\textit{column $5$}). Observe that the reconstructions are poor for $\tau \notin \tau_{proposed}$, and successful for $\tau \in \tau_{proposed}$.}
\end{center}
\end{figure}

We confirmed similar conclusions for MNIST digit $4$, and show sample reconstructions from both within, and outside (subjacent and superjacent), our proposed range in \cref{f:GoodBadRecon}. 

\section{Discussion\label{s:Discuss}}

We introduced \textbf{OCTANE}, a novel deep autoencoder derived from an optimal control formulation involving encoder and decoder differential equations. This framework merges continuous-time dynamics with neural network architecture, enabling principled design through the discretization of state, adjoint, and design equations. A key innovation lies in the use of rank-adaptive tensor solvers, which naturally enforce the butterfly-like structure of the autoencoder through compression and decompression via rank truncation.

Our experiments demonstrate OCTANE's effectiveness on image denoising and deblurring tasks. The model learns low-rank representations while significantly reducing memory usage—up to $\textbf{16.21\%}$ for denoising, and $\textbf{57.46\%}$ for deblurring—without sacrificing reconstruction quality. Training times are modest (under 40 minutes), even with minimal data (as few as 20 samples), and minimal tuning of parameters is required.

We also offer a heuristic strategy to select the step-size $\tau$ and number of layers $N$ to balance accuracy and computational cost. A practical interval for $\tau$ and corresponding $N$ values yields consistently strong performance across tasks.

Future work includes exploring advanced tensor solvers for time integration, extending to larger and higher-dimensional datasets, implementing the method in popular deep learning frameworks (e.g., PyTorch), and embedding OCTANE into bilevel optimization pipelines for hyperparameter tuning and generalization across modalities \cite{Khatri_2019_Bonnet}.

\section*{Acknowledgments}
R. Khatri would like to thank Abram (Bram) Rodgers (NASA Ames Research Center) and Andr\'{e}s Miniguano-Trujillo (Maxwell Institute for Mathematical Sciences) for useful discussions on tensors and their numerical implementations in the earlier stages of this work. The authors are extremely grateful to Sergey Dolgov (University of Bath) for his generous help on the use of \texttt{TT-toolbox}.

\bibliographystyle{plain}
\bibliography{OCTANE_KKOA_ref}


\end{document}